\documentclass[AMA,Times1COL]{WileyNJDv5} 
\usepackage[OT2,OT1]{fontenc}
\usepackage{bm}
\usepackage{tikz}
\usepackage{flowchart}
\usepackage{esint}
\usetikzlibrary{patterns}
\tikzset{>=stealth}
\numberwithin{equation}{section}
\usepackage[final]{changes}
\usetikzlibrary{positioning, shapes.geometric}

\articletype{Article Type}%

\received{00}
\revised{00}
\accepted{00}
\journal{XXXX}
\volume{00}
\copyyear{2024}
\startpage{1}

\begin{document}

\title{\replaced{Bidirectional conformal mapping for over-break and under-break tunnelling and its application in complex variable method}{Complex variable solution on noncircular and asymmetrical tunnelling embedded by bidirectional conformal mapping incorporating Charge Simulation Method}}

\author[1]{Luo-bin Lin}

\author[2]{Fu-quan Chen}

\author[1]{Chang-jie Zheng}

\author[1]{Shang-shun Lin}

\authormark{LIN \textsc{et al.}}
\titlemark{\replaced{Bidirectional conformal mapping for over-break and under-break tunnelling and its application in complex variable method}{Complex variable solution on noncircular and asymmetrical tunnelling embedded by bidirectional conformal mapping incorporating Charge Simulation Method}}

\address[1]{\orgdiv{Fujian Provincial Key Laboratory of Advanced Technology and Informatization in Civil Engineering, College of Civil Engineering}, \orgname{Fujian University of Technology}, \orgaddress{\state{Fujian}, \country{China}}}

\address[2]{\orgdiv{College of Civil Engineering}, \orgname{Fuzhou University}, \orgaddress{\state{Fujian}, \country{China}}}

\corres{Fu-quan Chen, College of Civil Engineering, Fuzhou University, No. 2 Xueyuan Road, Shangjie University Town, Fuzhou, 350108, Fujian, China\email{phdchen@fzu.edu.cn}}



\abstract[Abstract]{\replaced{Over-break and under-break excavation is very common in practical tunnel engineering with asymmetrical cavity contour, while existing conformal mapping schemes of complex variable method generally focus on tunnelling with theoretical and symmetrical cavity contour.}{Mechanical issues of noncircular and asymmetrical tunnelling can be estimated using complex variable method with suitable conformal mapping.} \replaced{Besides, the solution strategies of existing conformal mapping schemes for noncircular tunnel generally apply optimization theory,}{Existing solution schemes of conformal mapping for noncircular tunnel generally need iteration or optimization strategy,} and are thereby mathematically complicated. This paper proposes a new bidirectional conformal mapping for \replaced{over-break and under-break}{deep and shallow} tunnels of \deleted{noncircular and} asymmetrical \replaced{contours}{shapes} by incorporating Charge Simulation Method\replaced{, which}{. The solution scheme of this new bidirectional conformal mapping} only involves a pair of \added{forward and backward} linear systems, and is therefore logically straight-forward, computationally efficient, and practically easy in coding. New numerical strategies are developed to deal with possible sharp corners of cavity by small arc simulation and densified collocation points. Several numerical examples are presented to illustrate the geometrical usage of the new bidirectional conformal mapping. Furthermore, the new bidirectional conformal mapping is embedded into two complex variable solutions of \replaced{under-break}{noncircular and asymmetrical} shallow tunnelling in gravitational geomaterial with reasonable far-field displacement. The respective result comparisons with finite element solution and existing analytical solution show good agreements, indicating \replaced{that the new bidirectional conformal mapping would extend the mechanical application range of the complex variable method in practical over-break and under-break tunnelling.}{the feasible mechanical usage of the new bidirectional conformal mapping.}}


\jnlcitation{\cname{%
\author{Lin L},
\author{Chen F},
\author{Zheng C},
\author{Lin S}}.
\ctitle{\replaced{Bidirectional conformal mapping for over-break and under-break tunnelling and its application in complex variable method}{Complex variable solution on noncircular and asymmetrical tunnelling embedded by bidirectional conformal mapping incorporating Charge Simulation Method}} \cjournal{XXXX} \cvol{2024;00(00):1--1}.}

\maketitle



\section{Introduction}
\label{sec:introduction}

\begin{figure}[htb]
  \centering
  \added{
  \begin{tabular}{cc}
    (a) Case in Ref \cite{singh2005causes} & (b) Case in Ref \cite{mohammadi2018prediction} \\
    \includegraphics[width = 0.4\textwidth]{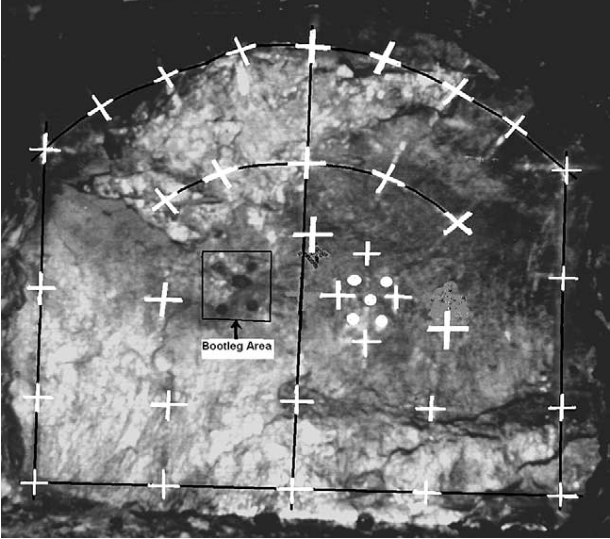} & \includegraphics[width = 0.4\textwidth]{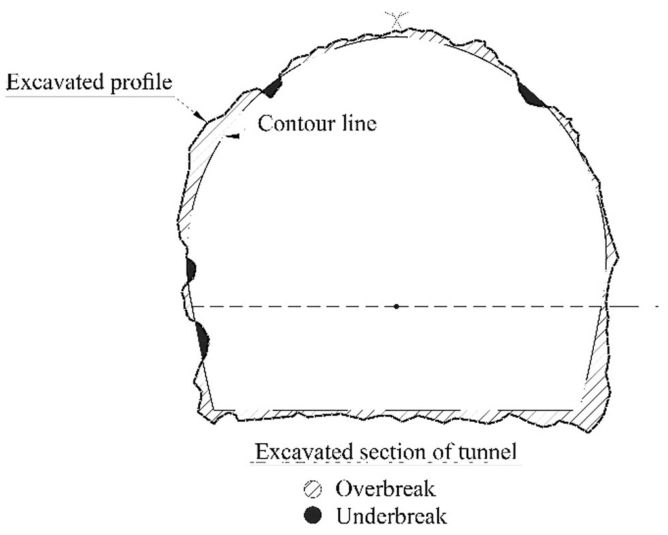} \\
    (c) Case in Ref \cite{mottahedi2018development} & (d) Case in Ref \cite{fodera2020factors} \\
    \includegraphics[width = 0.4\textwidth]{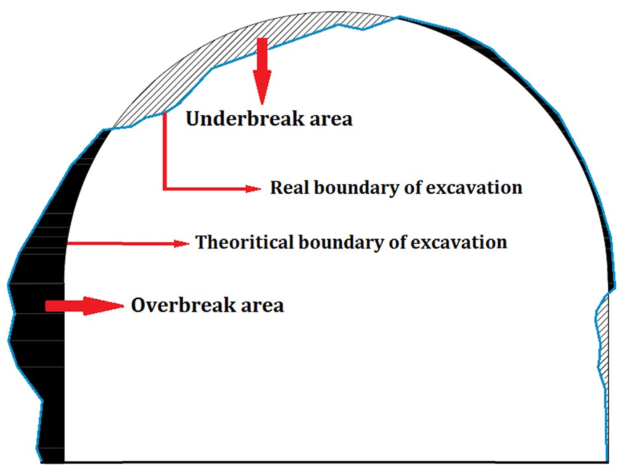} & \includegraphics[width = 0.4\textwidth]{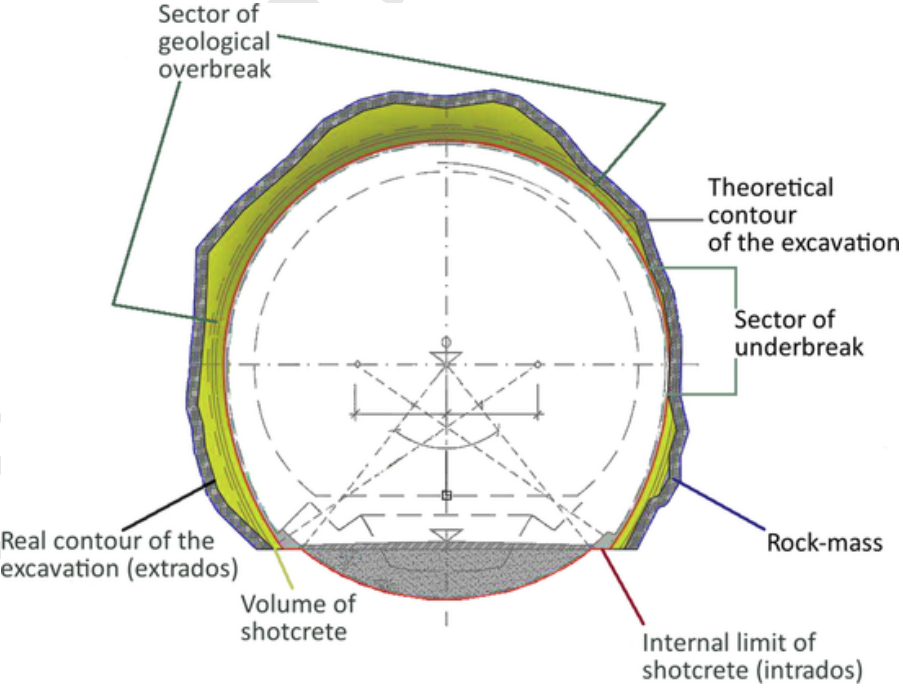} \\
  \end{tabular}
  }
  \caption{\added{Comparisons of real and theoretical excavations in tunnel engineering}}
  \label{fig:1}
\end{figure}

\replaced{Over-break and under-break excavations are common in real-world tunnel engineering, and the difference between practical and theoretical excavation contours generally exist. Fig. \ref{fig:1} shows several comparison cases of practical and theoretical excavation contours in real-world tunnel engineering \cite{singh2005causes,mandal2009evaluating,mottahedi2018development,fodera2020factors}. We can see that the cross sections of these tunnels with over-break and under-break excavations are generally asymmetrical. In other words, real-world tunnellings are always related to cavities of asymmetrical cross sections,}{Noncircular and asymmetrical excavation widely exists in tunnel engineering,} which can be studied using complex variable method \cite{Muskhelishvili1966} \added{with the aid of proper bidirectional conformal mapping}.

According \added{to} buried depth and consideration of gravitational gradient, tunnel engineering models in complex variable method can be \added{mechanically} attributed into deep and shallow categories. Geometrically, geomaterial containinng a deep tunnel is generally abstracted to an infinite two dimensional plane exterior of cavity boundary, and is considered as a simply connected region bounded by a closed Jordan curve. Meanwhile, the geomaterial containing a shallow tunnel is generally abstrated to an infinite lower half plane between an infinite ground surface and a finite cavity boundary, and is generally considered as a doubly-connected region bounded by an infinite straight line and a closed Jordan curve. The topologies of deep and shallow tunnels are different, and the conformal mapping schemes are different accordingly.

For a deep tunnel, the most widely used conformal mapping is the {\fontencoding{OT2}\selectfont Melentev} mapping \cite{Muskhelishvili1966}. The {\fontencoding{OT2}\selectfont Melentev} mapping applies an iteration strategy to gradually approximate the conformal mapping from mapping region to corresponding physical one, which can be naturally embedded into the complex variable method using both boundary condition formations of integral traction and displacement. Sufficient works on deep tunnels \cite{zhang01,zhang03:_stres,Huo2006,Kargar2014,Lu2014,Lu2015,fanghuangcheng2021,ye2024non} are based on the {\fontencoding{OT2}\selectfont Melentev} mapping. Several branches of of the {\fontencoding{OT2}\selectfont Melentev} mapping have been developed. An important branch \cite{exadaktylos2002closed} confines the upper limit of 2-norm sum of sample point coordinate differences between the physical cavity boundary and the mapping one to seek the the mapping coefficients. Technically, the coefficients should be solved by optimization or iteration strategy. Another feasible branch mapping scheme \cite{ma2022numerical} explicitly applies the least square method to solve the mapping coeffiicients using optimization strategy. Recently, a new conformal mapping for noncircular tunnels \cite{ye2023novel-aizhiyong} is proposed by combining least square method and iteration strategy. 

For \added{a} shallow tunnel, the most notable conformal mapping is the Verruijt mapping \cite{Verruijt1997traction,Verruijt1997displacement}, which is fully analytic to accurately and bidirectionally transform a lower half plane containing a circular cavity onto a circular unit annulus. Based on the Verruijt mapping \cite{Verruijt1997traction,Verruijt1997displacement}, the mechanical behaviour of shallow tunnelling is widely investigated \cite{Strack2002phdthesis,fang15:_compl,Lu2016,zhang2018complex,Wang2018shallow_surcharge,lu2019unified,ai2023stress,ye2024matrix,LIN2024appl_math_model}. Zeng et al. \cite{Zengguisen2019} make a significant extension of the Verruijt mapping by adding certain items in the mapping (see Eq. (\ref{eq:6.1})), which extends the usage of complex variable solution to noncircular shallow tunnelling issues \cite{lu2021complex,zeng2022analytical,zeng2023analytical,LIN2024comput_geotech_1,zhou2024analytical,fan2024analytical}.

\added{Though many significant complex variable solutions of deep and shallow tunnels have been conducted using corresponding conformal mappings, these studies generally focus on cavities of theoretical excavation contours, while asymmetrical cavities due to over-break and under-break excavations widely existing in real-world tunnelling are not considered. One possible reason is the lack of proper conformal mapping to simulate geomaterial containing an asymmetrical cavity, especially for a shallow geomaterial. Moreover, the conformal mapping schemes used in the above noncircular complex variable solutions are are generally computationally complicated to obtain results. }For the deep tunnel mapping schemes of {\fontencoding{OT2}\selectfont Melentev} \cite{Muskhelishvili1966} and its branches \cite{exadaktylos2002closed,ma2022numerical,ye2023novel-aizhiyong} mentioned above, the solution procedures are complicated and sometimes obscure with optimization or iteration strategy. For the noncircular shallow mapping scheme by Zeng et al. \cite{Zengguisen2019}, the nonconvex optimization strategy is applied to solved the mapping coefficients, which is laborious and generally requires computation infrastructure of high performance and possibly sophysicated coding skills (see Section \ref{sub:Shallow tunnel}). Therefore, straight-forward and simple mapping schemes for both deep and shallow tunnels \added{with over-break and under-break excavation} are necessary. 

In this paper, new bidirectional conformal mapping schemes for both deep and shallow tunnels \added{with over-break and under-break contours (asymmetrical cross sections)} are proposed by incorporating the Charge Simulation Method. The solution schemes for mapping coefficients only involve a pair of \added{forward and backward} linear systems, and are thereby straight-forward, very efficient, and also easy in programming codes. The Charge Simulation Method has been applied in our previous work \cite{lin2024over-under-excavation}, which contains defects on dealing with sharp corners (see Sections \ref{sub:Shallow tunnel} and \ref{sec:Second numerical case and comparisons with existing analytical solution}). The newly proposed bidirectional conformal mapping in this paper eliminates such defects by simulating sharp corners with corresponding rounded ones (see Fig. \ref{fig:3}).

The paper is organized as follows: Section \ref{sec:Charge Simulation Method} briefly introduces the background and outline of Charge Simulation Method; Section \ref{sec:Numerical conformal mapping} detailedly and analytically explains the formation and solution procedure of the mapping coefficients; Section \ref{sec:Numerical results} suppliments important numerical schemes to sustain the usage of the mappings and to estimate computation errors; Section \ref{sec:Numerical cases} puts forward numerical strategy to simulate sharp corners and several typical cases for both deep and shallow tunnels; Section \ref{sec:Comparisons with existing conformal mappings} analytically compares the proposed bidirectional conformal mappings with existing ones to exhibit the improvements; Section \ref{sec:Complex variable solution} embeds the proposed bidirectional conformal mapping into complex variable method to investigate two numerical cases of noncircular shallow tunnelling in gravitational geomaterial with reasonable far-field displacement, and certain comparisons are made for verification of the mapping, as well as the proposed complex variable solutions.

\section{Charge Simulation Method}%
\label{sec:Charge Simulation Method}

The Charge Simulation Method is originally proposed by Steinbigler \cite{steinbigler1969anfangsfeldstarken}, which is soon used in eletrical engineering. Then Amano found it very efficient in numerical conformal mapping \cite{amano1994charge}, which is the linear simplification of the numerical conformal mapping of boundary integration proposed by Hough and Papamichael \cite{hough1983integral}. Here, a brief introduction of Charge Simulation Method is given below.

For a Dirichlet problem of Laplace's equation:
\begin{subequations}
  \label{eq:2.1}
  \begin{equation}
    \label{eq:2.1a}
    \varDelta g(z) = 0, \quad z \in {\bm{D}}
  \end{equation}
  \begin{equation}
    \label{eq:2.1b}
    g(z) = b(z), \quad z \in \partial{\bm{D}}
  \end{equation}
\end{subequations}
where $ g(z) $ denotes the solution of the Dirichlet problem, $ \varDelta $ denotes the Laplace operator, $ z = x+{\rm{i}}y $ denotes the complex variable, $ {\bm{D}} $ is a domain with its boundary $ \partial{\bm{D}} $, $ z = x+{\rm{i}}y $ denotes the complex variable within the domain $ {\bm{D}} $ and along its boundary $ \partial{\bm{D}} $, and $ b(z) $ is a complex function defined on $ \partial{\bm{D}} $.

The solution $ g(z) $ within domain $ {\bm{D}} $ can be simulated by Charge Simulation Method as \cite{amano1994charge}
\begin{equation}
  \label{eq:2.2}
  g(z) := G(z) = \sum\limits_{k=1}^{N} Q_{k}\ln|z-Z_{k}|, \quad z \in {\bm{D}}
\end{equation}
where $ G(z) $ denotes the simulation of the original solution $ g(z) $ in Eq. (\ref{eq:2.1}); $ Z_{k} $ are called {\emph{charge points}}, which are located outside of the closure $ Z_{k} \notin \overline{\bm{D}} = {\bm{D}} \cup \partial{\bm{D}} $, and should be given beforehand; $ Q_{k} $ are called {\emph{charges}}, which are real constants to be determined to meet the boundary condition in Eq. (\ref{eq:2.1b}). Eq. (\ref{eq:2.2}) indicates that the two-dimensional Dirichlet problem in Eq. (\ref{eq:2.1a}) can be approximated by a linear combination of fundamental solutions of two-dimensional Laplace operator.

Substituting Eq. (\ref{eq:2.2}) into Eq. (\ref{eq:2.1b}) gives its {\emph{collocation conditions}} as
\begin{equation}
  \label{eq:2.3}
  \sum\limits_{k=1}^{N} Q_{k}\ln|z_{i}-Z_{k}| = b(z_{i}), \quad z_{i} \in \partial{\bm{D}}, i = 1,2,3,\cdots,N
\end{equation}
where $ z_{i} $ are called {\emph{collocation points}}, and are a series of sequentially selected points along the boundary $ \partial{\bm{D}} $. The quantity of {\emph{collocation points}} is the same to that of the {\emph{charge points}} $ Z_{k} $, so that Eq. (\ref{eq:2.3}) would form a simultaneous linear system to uniquely identify the undertemined {\emph{charges}} $ Q_{k} $. It should be addressed that as long as the boundary $ \partial{\bm{D}} $ is an analytic Jordan curve, the Charge Simulation Method is applicable.

\section{Numerical conformal mapping in tunnel engineering}
\label{sec:Numerical conformal mapping}

The Charge Simulation Method can be applied to obtain the numerical conformal mapping in tunnel engineering with certain modifications. A tunnel engineering problem can be categorized into deep or shallow tunnelling, according to tunnel depth. \deleted{For a deep tunnel, bilateral far-field stress field is applied onto the infinite geomaterial; For a shallow tunnel, gravitationl gradient of geomaterial can not be ignored, and the ground surface should be considered.} Both tunnelling problems can be simplified as plain strain one, and the typical geometries are shown in Figs. \ref{fig:2}a-1 and b-1, where a cavity of arbitrary shape is excavated within a two-dimensional geomaterial $ {\bm{D}}: z = x+{\rm{i}}y $. Fig. \ref{fig:2} geometrically indicates that the two-dimensional geomaterials of deep and shallow tunnellings have different connectivities, which involve different conformal mapping schemes. \added{To be clear, this section only focuses on conformal mappings, and no mechanical issues would be discussed.}

\subsection{Deep tunnel in infinite geomaterial}
\label{sec:Deep tunnel in infinite geomaterial}

As shown in Fig. \ref{fig:2}a-1, an infinite geomaterial weakened by a deep tunnel can be geometrically attributed to a simply-connected domain $ {\bm{D}} : z = x+{\rm{i}}y $ exterior of an arbitrary cavity boundary $ \partial{\bm{D}} $, and the corresponding mapping plane $ {\bm{d}} : \zeta = \rho \cdot {\rm{e}}^{{\rm{i}}\theta} $ is shown in Fig. \ref{fig:2}b-1 with unit circular boundary $ \partial{\bm{d}} $. The bidirectional conformal mapping can be expressed as
\begin{subequations}
  \label{eq:3.1}
  \begin{equation}
    \label{eq:3.1a}
    \zeta = \zeta(z), \quad z \in \overline{\bm{D}}
  \end{equation}
  \begin{equation}
    \label{eq:3.1b}
    z = z(\zeta), \quad \zeta \in \overline{\bm{d}}
  \end{equation}
\end{subequations}
Eqs. (\ref{eq:3.1a}) and (\ref{eq:3.1b}) are called forward and backward conformal mappings, respectively. The former maps the physical infinite plane bounded by an arbitrary cavity boundary $ \overline{\bm{D}} $ onto the mapping infinite plane bounded by a unit cavity $ \overline{\bm{d}} $, while the latter does the opposite.

The Charge Simulation Method is applied to the forward conformal mapping. The forward conformal mapping of infinite exterior domain would exist under the normalization condition $ \zeta(z_{c1}) = 0 (z_{c1} \notin \overline{\bm{D}}) $, the infinity condition $ \zeta(\infty) = \infty $, and the first-order derivative conditions $ \zeta^{\prime}(\infty) > 0 $ \cite{nehari1975conformal}. Thus, the forward conformal mapping can be expressed as
\begin{equation}
  \label{eq:3.2}
  \zeta(z) = (z-z_{c1}) \cdot \exp \left[ \varGamma + g(z) + {\rm{i}}h(z) \right], \quad z \in \overline{\bm{D}}
\end{equation}
where $ \varGamma = \lim\limits_{z\rightarrow\infty}\ln|\zeta^{\prime}(z)| $, $ g(z) $ and $ h(z) $ are conjugate harmonic functions within domain $ {\bm{D}} $. The normalization condition $ \zeta(z_{c1}) = 0 $ is satisfied in a clear manner, while the rest two conditions would rely on the formations of $ g(z) $ and $ h(z) $. We propose that
\begin{equation}
  \label{eq:3.3}
  g(z) + {\rm{i}}h(z) := G(z) + {\rm{i}}H(z) = \sum\limits_{k=1}^{N} Q_{k}\ln(z-Z_{k}), \quad z \in \overline{\bm{D}}
\end{equation}
with the following constraint
\begin{equation}
  \label{eq:3.4}
  \sum\limits_{k=1}^{N} Q_{k} = 0
\end{equation}
would meet the other two conditions, where $ G(z) $ and $ H(z) $ are approximations of $ g(z) $ and $ h(z) $, respectively; $ Q_{k} $ and $ Z_{k} $ are the same to the predefined notations. The formation in Eqs. (\ref{eq:3.3}) and (\ref{eq:3.4}) can be briefly proven as follows.

Separating the real and imaginary parts of Eq. (\ref{eq:3.3}) gives
\begin{subequations}
  \label{eq:3.5}
  \begin{equation}
    \label{eq:3.5a}
    g(z) := G(z) = \sum\limits_{k=1}^{N} Q_{k}\ln|z-Z_{k}|
  \end{equation}
  \begin{equation}
    \label{eq:3.5b}
    h(z) := H(z) = \sum\limits_{k=1}^{N} Q_{k}\arg(z-Z_{k})
  \end{equation}
\end{subequations}
Eq. (\ref{eq:3.5a}) is in the same formation as the simulation in Eq. (\ref{eq:2.2}), while Eq. (\ref{eq:3.5b}) is the conjugate imaginary part. 

Substituting Eq. (\ref{eq:3.3}) into Eq. (\ref{eq:3.2}) gives
\begin{equation}
  \label{eq:3.6}
  \zeta(z) := \tilde{\zeta}(z) = (z-z_{c1}) \cdot \exp \left[ \varGamma + \sum\limits_{k=1}^{N} Q_{k}\ln(z-Z_{k}) \right], \quad z \in \overline{\bm{D}}
\end{equation}
where $ \tilde{\zeta}(z) $ denotes the simulation of $ \zeta(z) $ with Eq. (\ref{eq:3.3}).
From Eq. (\ref{eq:3.6}), we have
\begin{subequations}
  \label{eq:3.7}
  \begin{equation}
    \label{eq:3.7a}
    \tilde{\zeta}(\infty) = \lim\limits_{z\rightarrow\infty} (z-z_{c1}) \cdot \exp\varGamma \cdot \exp\left[ \sum\limits_{k=1}^{N} Q_{k}\ln(z-Z_{k}) \right]
  \end{equation}
  \begin{equation}
    \label{eq:3.7b}
    \tilde{\zeta}^{\prime}(\infty) = \lim\limits_{z\rightarrow\infty} \exp\varGamma \cdot \exp\left[ \sum\limits_{k=1}^{N} Q_{k}\ln(z-Z_{k}) \right] \cdot \left[ 1 + \sum\limits_{k=1}^{N} Q_{k}\frac{z-z_{c}}{z-Z_{k}} \right]
  \end{equation}
\end{subequations}
From Eq. (\ref{eq:3.7}), we have
\begin{subequations}
  \label{eq:3.8}
  \begin{equation}
    \label{eq:3.8a}
    \lim\limits_{z\rightarrow\infty} \sum\limits_{k=1}^{N} Q_{k} \ln|z-Z_{k}| = \sum\limits_{k=1}^{N} Q_{k} \cdot C_{\infty}
  \end{equation}
  \begin{equation}
    \label{eq:3.8b}
    \lim\limits_{z\rightarrow\infty} \sum\limits_{k=1}^{N} Q_{k} \arg(z-Z_{k}) = \sum\limits_{k=1}^{N} Q_{k} \cdot D_{\infty}
  \end{equation}
  \begin{equation}
    \label{eq:3.8c}
    \lim\limits_{z\rightarrow\infty}\sum\limits_{k=1}^{N} Q_{k}\frac{z-z_{c1}}{z-Z_{k}} = \sum\limits_{k=1}^{N} Q_{k}
  \end{equation}
\end{subequations}
where $ C_{\infty} \leftarrow \lim\limits_{z\rightarrow\infty}\ln|z-Z_{k}| $ and $ D_{\infty} \leftarrow \lim\limits_{z\rightarrow\infty}\arg(z-Z_{k}) $ for $ k = 1,2,3,\cdots,N $.

With the constraint of $ Q_{k} $ in Eq. (\ref{eq:3.4}), Eq. (\ref{eq:3.8}) trends to zero, and Eq. (\ref{eq:3.7}) can be simplified as
\begin{equation}
  \label{eq:3.7a'}
  \tag{3.7a'}
  \tilde{\zeta}(\infty) = \lim\limits_{z\rightarrow\infty} (z-z_{c}) \cdot \exp\varGamma = \infty
\end{equation}
\begin{equation}
  \label{eq:3.7b'}
  \tag{3.7b'}
  \tilde{\zeta}^{\prime}(\infty) = \lim\limits_{z\rightarrow\infty} \exp\varGamma = \exp\varGamma
\end{equation}
It should be addressed that the values of Eqs. (\ref{eq:3.7a'}), (\ref{eq:3.7b'}), and (\ref{eq:3.8}) are simultaneously achieved with $ z\rightarrow\infty $. Eqs. (\ref{eq:3.7a'}) and (\ref{eq:3.7b'}) are the simulations of the other two conditions of conformal mapping existence. Therefore, Eqs. (\ref{eq:3.3}) and (\ref{eq:3.4}) would meet the requirements. In other words, the forward conformal mapping that maps the two-dimensional geomaterial bounded by a cavity boundary $ \overline{\bm{D}} $ in Fig. \ref{fig:2}a-1 to its image $ \overline{\bm{d}} $ in Fig. \ref{fig:2}a-2 can be expressed in a practical manner as
\begin{equation}
  \label{eq:3.9}
  \zeta(z) = (z-z_{c1}) \cdot \exp\left[ \varGamma + \sum\limits_{k=1}^{N} Q_{k}\ln(z-Z_{k}) \right], \quad z \in \overline{\bm{D}}
\end{equation}
Eq. (\ref{eq:3.9}) is first proposed in our previous study \cite{lin2024over-under-excavation}, but no mathematical proof is given.

Eq. (\ref{eq:3.9}) contains $ N+1 $ unknown real constants, which can be determined in the following scheme:

(1) Sequentially select $ N $ "suitable" {\emph{collocation points}} along the cavity boundary $ \partial{\bm{D}} $ in the positive direction (clockwisely to always keep geomaterial on the left), and these "suitable" {\emph{collocation points}} are denoted by $ z_{i} $ ($ i = 1,2,3,\cdots,N $). \added{The selection strategy of {\emph{collocation points}} are sensitive to the mapping accuracy, and are thereby very flexible. However, certain selection principles can be given in Section \ref{sub:Selection strategy of collocation points} to guarantee mapping accuracy.}

(2) For each {\emph{collocation point}}, the corresponding {\emph{charge point}} is selected in Amano's manner \cite{amano1994charge} as
\begin{equation}
  \label{eq:3.10}
  \left\{
    \begin{aligned}
      & Z_{k} = z_{k} + K_{0} \cdot H_{k} \cdot \exp ({\rm{i}}\varTheta_{k}) \\
      & H_{k} = \frac{1}{2}(|z_{k+1}-z_{k}|+|z_{k}-z_{k-1}|) \\
      & \varTheta_{k} = \arg(z_{k+1}-z_{k-1})-\frac{\pi}{2} \\
    \end{aligned}
  \right.
  , \quad k = 1,2,3,\cdots,N
\end{equation}
where $ K_{0} $ denotes an {\emph{assignment factor}} to adjust the distance between {\emph{collocation point}} and its corresponding {\emph{charge point}}. Eq. (\ref{eq:3.10}) indicates that the {\emph{charge points}} are distributed perpendicularly to and near the {\emph{collocation points}}.

(3) Eq. (\ref{eq:3.9}) can be rewritten as
\begin{equation*}
  \ln\zeta(z) = \ln(z-z_{c1}) + \varGamma + \sum\limits_{k=1}^{N} Q_{k}\ln(z-Z_{k})
\end{equation*}
With the $ N $ selected {\emph{collocation points}} and corresponding {\emph{charge points}}, the following $ N $ {\emph{collocation conditions}} can be established according to the real part of above equation as
\begin{equation}
  \label{eq:3.11}
  \sum\limits_{k=1}^{N} Q_{k}\ln|z_{i}-Z_{k}| + \varGamma = -\ln|z_{i}-z_{c1}|, \quad i = 1,2,3,\cdots,N
\end{equation}
where $ \ln|\zeta(z_{i})| = 0 $ should be theoretically satisfied, since $ \zeta(z_{i}) $ is located along the mapping boundary $ \partial{\bm{d}} $, whose radius is unity. Together with the constraint of $ Q_{k} $ in Eq. (\ref{eq:3.4}), a well-defined simultaneous linear system of $ N+1 $ equations and $ N+1 $ variables can be obtained, and the unique solution of the real contants $ Q_{k} $ \added{and $ \varGamma $} can be achieved.

(4) With the unique set of $ Q_{k} $, the forward conformal mapping in Eq. (\ref{eq:3.9}) is uniquely determined so far.

In complex variable solution, the backward conformal mapping $ z = z(\zeta) $ that maps the image in Fig. \ref{fig:2}a-2 onto its preimage in Fig. \ref{fig:2}a-1 is in need to construct boundary conditions, as well as stress and displacement solutions. The backward conformal mapping can be established using point correspondence method.

The backward conformal mapping should be in the following formation \cite{Muskhelishvili1966} as
\begin{equation}
  \label{eq:3.12}
  z(\zeta) = q_{-1}\zeta + q_{0} + \sum\limits_{k=1}^{N-2} q_{k}\zeta^{-k}, \quad \zeta \in \overline{\bm{d}}
\end{equation}
where the coefficients $ q_{k} $ are complex numbers, and can be determined by point correspondence method as
\begin{equation}
  \label{eq:3.13}
  \zeta(z_{i}) \cdot q_{-1} + 1 \cdot q_{0} + \sum\limits_{k=1}^{N-2} \zeta(z_{i})^{-k} \cdot q_{k} = z_{i}, \quad i = 1,2,3,\cdots,N
\end{equation}
Eq. (\ref{eq:3.13}) is a simultaneous linear system containing $ N $ variables and $ N $ complex equations. With Eq. (\ref{eq:3.13}), the backward conformal mapping can be obtained.

The solution method in Eqs. (\ref{eq:3.12}) and (\ref{eq:3.13}) is theoretically accurate if the modulus precisely match $ |\zeta(z_{i})| = 1 $ ($ i = 1,2,3,\cdots,N $). However, due to possible computation error of the forward conformal mapping in Eq. (\ref{eq:3.9}), we only have $ |\zeta(z_{i})| \approx 1 $. Therefore, the high-order item $ \zeta(z_{i})^{-k} $ in Eq. (\ref{eq:3.13}) would significantly enlarge the small errors to cause a huge condition number (see Section \ref{sub:Condition number and error estimates}). To avoid such significant errors, the backward conformal mapping in Eq. (\ref{eq:3.12}) can be simulated by the following one as
\begin{equation}
  \label{eq:3.12a'}
  \tag{3.12a'}
  z(\zeta) = p_{-1}\zeta + p_{0} + \sum\limits_{k=1}^{N-2} \frac{p_{k}}{\zeta-\zeta_{k}}, \quad \zeta \in \overline{\bm{d}}
\end{equation}
where $ p_{k} $ are complex coefficients to be determined, and $ \zeta_{k} $ are singular points distributed inside of boundary $ \partial{\bm{d}} $, similar to the {\emph{charge points}} in the forward conformal mapping:
\begin{equation}
  \label{eq:3.12b'}
  \tag{3.12b'}
  \left\{
    \begin{aligned}
      & \zeta_{k} = \zeta(z_{k}) + k_{0} \cdot h_{k} \cdot \exp({\rm{i}}\vartheta_{k}) \\
      & h_{k} = \frac{1}{2} \left[ \left| \zeta(z_{k+1}) - \zeta(z_{k}) \right| + \left| \zeta(z_{k}) - \zeta(z_{k-1}) \right| \right] \\
      & \vartheta_{k} = \arg\left[ \zeta(z_{k+1}) - \zeta(z_{k-1}) \right] - \frac{\pi}{2} \\
    \end{aligned}
  \right.
\end{equation}
wherein $ k_{0} $ is an {\emph{assignment factor}} to adjust the distance between points $ \zeta(z_{k}) $ and $ \zeta_{k} $. The complex coefficients $ p_{k} $ in Eq. (\ref{eq:3.12a'}) can also be determined by point correspondence as
\begin{equation}
  \label{eq:3.13'}
  \tag{3.13'}
  \zeta(z_{i}) \cdot p_{-1} + 1 \cdot p_{0} + \sum\limits_{k=1}^{N-2} \frac{p_{k}}{\zeta(z_{i})-\zeta_{k}} = z_{i}, \quad i = 1,2,3,\cdots,N
\end{equation}
Eq. (\ref{eq:3.13'}) is a well-defined simultaneous linear system to obtain unique results of $ p_{k} $. 

Therefore, Eqs. (\ref{eq:3.9}) and (\ref{eq:3.12a'}) complete the bidirectional conformal mapping between the geomaterial bounded by a cavity boundary $ \overline{\bm{D}} $ in Fig. \ref{fig:2}a-1 and its image $ \overline{\bm{d}} $ in Fig. \ref{fig:2}a-2. The fraction expression of Eq. (\ref{eq:3.12a'}) formally guarantees its analyticity.

\subsection{Shallow tunnel in infinite lower half geomaterial}
\label{sec:Shallow tunnel in infinite lower half geomaterial}

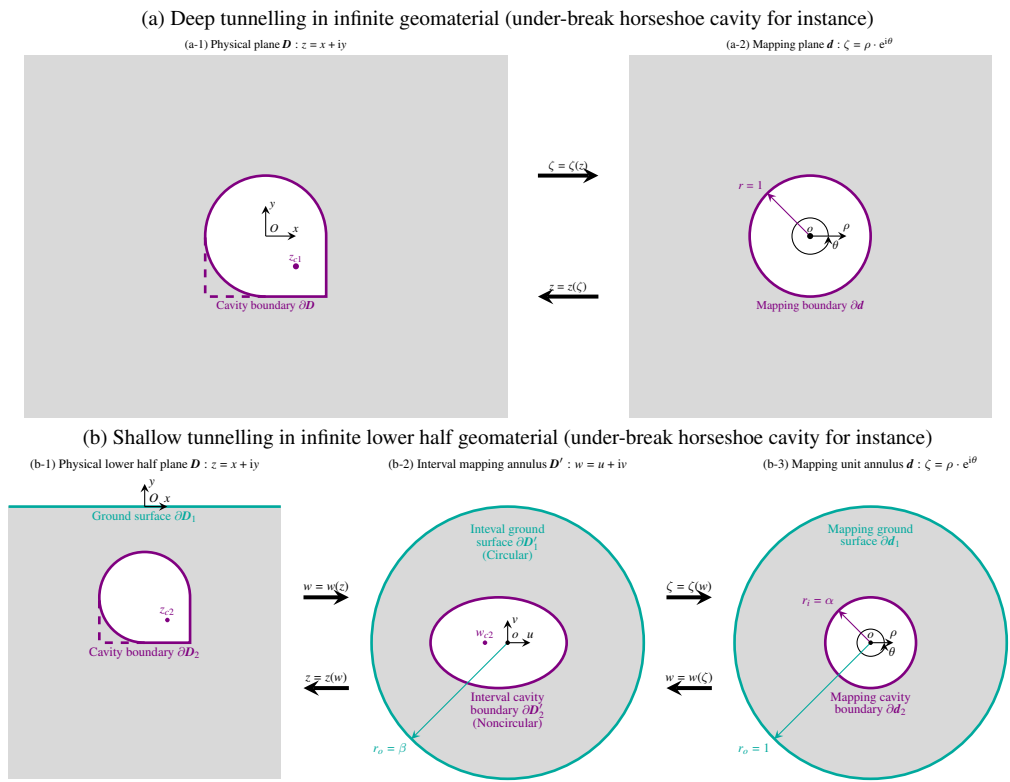
\begin{figure}[htpb]
  \begin{center}
    \small
    \begin{tabular}{c}
      (a) Deep tunnelling in infinite geomaterial (under-break horseshoe cavity for instance) \\
      \begin{tikzpicture}[scale=0.8, transform shape]
        \tikzstyle{every node} = [scale = 0.7]
        \fill [gray!30] (-4,0) rectangle (4,-6);
        \fill [white] (1,-3) arc [start angle = 0, end angle = 270, radius = 1] -- (1,-4) -- (1,-3);
        \draw [violet, line width = 1pt] (1,-3) arc [start angle = 0, end angle = 270, radius = 1] -- (1,-4) -- (1,-3);
        \draw [violet, line width = 1pt, dashed] (0,-4) -- (-1,-4) -- (-1,-3);
        \draw [->] (0,-3) -- (0.5,-3) node [above] {$ x $};
        \draw [->] (0,-3) -- (0,-2.5) node [right] {$ y $};
        \node at (0,-3) [above right] {$ O $};
        \fill [violet] (0.5,-3.5) circle [radius = 0.05];
        \node at (0.5,-3.5) [above, violet] {$ z_{c1} $};
        \node at (0,-4) [below, violet] {Cavity boundary $ \partial{\bm{D}} $};
        \node at (0,0) [above] {(a-1) Physical plane $ {\bm{D}} : z = x+{\rm{i}}y $};

        \fill [gray!30] (6,0) rectangle (12,-6);
        \fill [white] (9,-3) circle [radius = 1];
        \draw [violet, line width = 1pt] (9,-3) circle [radius = 1];
        \draw [->] (9,-3) -- (9.6,-3) node [above] {$ \rho $};
        \draw [->] (9.3,-3) arc [start angle = 0, end angle = 360, radius = 0.3];
        \node at (9.3,-3) [below right] {$ \theta $};
        \node at (9,-3) [above] {$ o $};
        \draw [->, violet] (9,-3) -- ({9+cos(135)},{-3+sin(135)}) node [above left] {$ r = 1 $};
        \fill (9,-3) circle [radius = 0.05];
        \node at (9,-4) [below, violet] {Mapping boundary $ \partial{\bm{d}} $};
        \node at (9,0) [above] {(a-2) Mapping plane $ {\bm{d}} : \zeta = \rho \cdot {\rm{e}}^{{\rm{i}}\theta} $};

        \draw [->, line width = 1.5pt] (4.5,-2) -- (5.5,-2);
        \draw [<-, line width = 1.5pt] (4.5,-4) -- (5.5,-4);
        \node at (5,-2) [above] {$ \zeta = \zeta(z) $};
        \node at (5,-4) [above] {$ z = z(\zeta) $};
      \end{tikzpicture} \\
      (b) Shallow tunnelling in infinite lower half geomaterial (under-break horseshoe cavity for instance) \\
      \begin{tikzpicture}[scale=0.6, transform shape]
        \tikzstyle{every node} = [scale = 1]
        \fill [gray!30] (-3,0) rectangle (3,-6);
        \draw [Emerald, line width = 1pt] (-3,0) -- (3,0);
        \fill [white] (1,-2) arc [start angle = 0, end angle = 270, radius = 1] -- (1,-3) -- (1,-2);
        \draw [violet, line width = 1pt] (1,-2) arc [start angle = 0, end angle = 270, radius = 1] -- (1,-3) -- (1,-2);
        \draw [violet, line width = 1pt, dashed] (0,-3) -- (-1,-3) -- (-1,-2);
        \draw [->] (0,0) -- (0.5,0) node [above] {$ x $};
        \draw [->] (0,0) -- (0,0.5) node [right] {$ y $};
        \node at (0,0) [above right] {$ O $};
        \fill [violet] (0.5,-2.5) circle [radius = 0.05];
        \node at (0.5,-2.5) [above, violet] {$ z_{c2} $};
        \node at (0,-3) [below, violet] {Cavity boundary $ \partial{\bm{D}}_{2} $};
        \node at (0,0) [below, Emerald] {Ground surface $ \partial{\bm{D}}_{1} $};
        \node at (0,0.7) [above] {(b-1) Physical lower half plane $ {\bm{D}} : z = x+{\rm{i}}y $};
        
        \fill [gray!30] (8,-3) circle [radius = 3];
        \fill [white] (7.8,-3) ellipse [x radius = 1.5, y radius = 1];
        \draw [Emerald, line width = 1pt] (8,-3) circle [radius = 3];
        \draw [violet, line width = 1pt] (7.8,-3) ellipse [x radius = 1.5, y radius = 1];
        \draw [->] (8,-3) -- (8.5,-3) node [above] {$ u $};
        \draw [->] (8,-3) -- (8,-2.5) node [right] {$ v $};
        \node at (8,-3) [above right] {$ o $};
        \fill [violet] (7.5,-3) circle [radius = 0.05];
        \node at (7.5,-3) [above, violet] {$ w_{c2} $};
        \draw [->, Emerald] (8,-3) -- ({8+3*cos(225)},{-3+3*sin(225)}) node [below left] {$ r_{o} = \beta $};
        \fill (8,-3) circle [radius = 0.05];
        \node at (8,-4) [below, violet, align = center] {Interval cavity \\ boundary $ \partial{\bm{D}}^{\prime}_{2} $ \\ (Noncircular)};
        \node at (8,-0.3) [below, Emerald, align = center] {Inteval ground \\ surface $ \partial{\bm{D}}^{\prime}_{1} $ \\ (Circular)};
        \node at (8,0.7) [above] {(b-2) Interval mapping annulus $ {\bm{D}}^{\prime} : w = u+{\rm{i}}v $};

        \draw [->, line width = 1.5pt] (3.5,-2) -- (4.5,-2);
        \draw [<-, line width = 1.5pt] (3.5,-4) -- (4.5,-4);
        \node at (4,-2) [above] {$ w = w(z) $};
        \node at (4,-4) [above] {$ z = z(w) $};
        
        \fill [gray!30] (16,-3) circle [radius = 3];
        \fill [white] (16,-3) circle [radius = 1];
        \draw [Emerald, line width = 1pt] (16,-3) circle [radius = 3];
        \draw [violet, line width = 1pt] (16,-3) circle [radius = 1];
        \draw [->] (16,-3) -- (16.5,-3) node [above] {$ \rho $};
        \draw [->] (16.3,-3) arc [start angle = 0, end angle = 360, radius = 0.3];
        \node at (16.3,-3) [below right] {$ \theta $};
        \node at (16,-3) [above] {$ o $};
        \draw [->, Emerald] (16,-3) -- ({16+3*cos(225)},{-3+3*sin(225)}) node [below left] {$ r_{o} = 1 $};
        \draw [->, violet] (16,-3) -- ({16+cos(135)},{-3+sin(135)}) node [above left] {$ r_{i} = \alpha $};
        \fill (16,-3) circle [radius = 0.05];
        \node at (16,-4) [below, violet, align = center] {Mapping cavity \\ boundary $ \partial{\bm{d}}_{2} $};
        \node at (16,-0.3) [below, Emerald, align = center] {Mapping ground \\ surface $ \partial{\bm{d}}_{1} $};
        \node at (16,0.7) [above] {(b-3) Mapping unit annulus $ {\bm{d}} : \zeta = \rho \cdot {\rm{e}}^{{\rm{i}}\theta} $};

        \draw [->, line width = 1.5pt] (11.5,-2) -- (12.5,-2);
        \draw [<-, line width = 1.5pt] (11.5,-4) -- (12.5,-4);
        \node at (12,-2) [above] {$ \zeta = \zeta(w) $};
        \node at (12,-4) [above] {$ w = w(\zeta) $};
      \end{tikzpicture} \\
    \end{tabular}
  \end{center}
  \caption{Typical conformal mappings of deep and shallow tunnelling}%
  \label{fig:2}
\end{figure}

In most complex variable solutions, a shallow tunnel problem is conducted in an infinite lower half plane. As shown in Fig. \ref{fig:2}b-1, an infinite lower half geomaterial weakened by a shallow tunnel can be geometrically categorized as a doubly-connected lower half domain $ {\bm{D}} : z = x+{\rm{i}}y $ bounded by ground surface $ \partial{\bm{D}}_{1} $ and cavity boundary of arbitrary shape $ \partial{\bm{D}}_{2} $, and $ \overline{\bm{D}} = \partial{\bm{D}}_{1} \cup {\bm{D}} \cup \partial{\bm{D}}_{2} $. 

The bidirectional conformal mapping of a lower half plane containing a cavity of arbitrary shape is decomposed into two steps, as shown in Fig. \ref{fig:2}b. The first step is to map the physical lower half plane closure $ \overline{\bm{D}} $ onto an interval annulus $ {\bm{D}}^{\prime} : w = u+{\rm{i}}v $ bounded by a circular exterior boundary $ \partial{\bm{D}}^{\prime}_{1} $ of exterior radius $ \beta $ and a noncircular interior boundary $ \partial{\bm{D}}^{\prime}_{2} $ ($ \overline{\bm{D}}^{\prime} = \partial{\bm{D}}^{\prime}_{1} \cup {\bm{D}}^{\prime} \cup \partial{\bm{D}}^{\prime}_{2} $) in bidirectional manner as
\begin{subequations}
  \label{eq:3.14}
  \begin{equation}
    \label{eq:3.14a}
    w = w(z) = \beta\frac{z-z_{c2}}{z-\overline{z}_{c2}}, \quad z \in \overline{\bm{D}}
  \end{equation}
  \begin{equation}
    \label{eq:3.14b}
    z = z(w) = \frac{w\overline{z}_{c2}-\beta z_{c2}}{w-\beta}, \quad w \in \overline{\bm{D}}^{\prime}
  \end{equation}
\end{subequations}
where $ z_{c2} $ and $ w_{c2} $ are arbitrary points within the cavity boundaries $ \partial{\bm{D}}_{2} $ and  $ \partial{\bm{D}}^{\prime}_{2} $, respectively.

The second step is to map the interval annulus closure $ \overline{\bm{D}}^{\prime} $ onto a unit annulus $ {\bm{d}} : \zeta = \rho \cdot {\rm{e}}^{{\rm{i}}\theta} $ bounded by a circular exterior boundary $ \partial{\bm{d}}_{1} $ of unity radius and a circular interior boundary $ \partial{\bm{d}}_{2} $ of $ \alpha $ radius ($ \overline{\bm{d}} = \partial{\bm{d}}_{1} \cup {\bm{d}} \cup \partial{\bm{d}}_{2} $) in a bidirectional manner as
\begin{subequations}
  \label{eq:3.15}
  \begin{equation}
    \label{eq:3.15a}
    \zeta = \zeta(w), \quad w \in \overline{\bm{D}}^{\prime}
  \end{equation}
  \begin{equation}
    \label{eq:3.15b}
    w = w(\zeta), \quad \zeta \in \overline{\bm{d}}
  \end{equation}
\end{subequations}
Via the stepwise bidirectional conformal mappings in Eqs. (\ref{eq:3.14}) and (\ref{eq:3.15}), the full bidirectional conformal mapping can be composed as
\begin{subequations}
  \label{eq:3.16}
  \begin{equation}
    \label{eq:3.16a}
    \zeta = \zeta[w(z)], \quad z \in \overline{\bm{D}}
  \end{equation}
  \begin{equation}
    \label{eq:3.16b}
    z = z[w(\zeta)], \quad \zeta \in \overline{\bm{d}}
  \end{equation}
\end{subequations}

The Charge Simulation Method is applied in the forward conformal mapping of the second-step mapping, and the details are elaborated below. According to the general formation of Charge Simulation Method in Eq. (\ref{eq:3.2}), the forward conformal mapping can be abstractly expressed as
\begin{equation}
  \label{eq:3.17}
  \zeta = \zeta(w) = (w-w_{c2})\exp\left[ \varLambda + g(w) + {\rm{i}}h(w) \right], \quad w \in \overline{\bm{D}}^{\prime}
\end{equation}
where $ \varLambda $ is a constant to be determined, $ g(w) $ and $ h(w) $ are conjugate harmonic functions within domain $ {\bm{D}}^{\prime} $. Eq. (\ref{eq:3.17}) clearly satisfies the normalization requirement $ \zeta(w_{c2}) = 0 $.

Since the doubly-connected annuli in Fig. \ref{fig:2}b-2 and b-3 are all finite, three more requirements should be imposed to the forward conformal mapping: (1) a rotation normalization should be imposed to Eq. (\ref{eq:3.17}) as $ \zeta(w_{\beta}) = 1 $, where $ w_{\beta} $ is an arbitrary point along the exterior boundary $ \partial{\bm{D}}^{\prime}_{1} $; (2) scaling invariant of the forward conformal mapping; (3) single-valuedness of the forward conformal mapping. To meet such the rotation normalization requirement (1), a possible formation may be set as
\begin{subequations}
  \label{eq:3.18}
  \begin{equation}
    \label{eq:3.18a}
    (w_{\beta}-w_{c2}) \cdot \exp\varLambda = 1
  \end{equation}
  \begin{equation}
    \label{eq:3.18b}
    g(w_{\beta}) + {\rm{i}}h(w_{\beta}) = 0
  \end{equation}
\end{subequations}

With Eq. (\ref{eq:3.18a}), Eq. (\ref{eq:3.17}) can be written into a more specific formation as
\begin{equation}
  \label{eq:3.17'}
  \tag{3.17'}
  \zeta = \zeta(w) = \frac{w-w_{c2}}{w_{\beta}-w_{c2}} \cdot \exp\left[ g(w) + {\rm{i}}h(w) \right]
\end{equation}
Similar to the deep tunnel, $ g(w)+{\rm{i}}h(w) $ can be simulated as \cite{okano2003numerical}
\begin{equation}
  \label{eq:3.19}
  g(w) + {\rm{i}}h(w) := G(w)+{\rm{i}}H(w) = \sum\limits_{k=1}^{N_{1}^{\prime}} Q_{1,k}\ln\frac{w-W_{1,k}}{w_{\beta}-W_{1,k}} + \sum\limits_{k=1}^{N_{2}} Q_{2,k}\ln\frac{w-W_{2,k}}{w_{\beta}-W_{2,k}}
\end{equation}
with the following constraints
\begin{subequations}
  \label{eq:3.20}
  \begin{equation}
    \label{eq:3.20a}
    \sum\limits_{k=1}^{N_{1}^{\prime}} Q_{1,k} = -1
  \end{equation}
  \begin{equation}
    \label{eq:3.20b}
    \sum\limits_{k=1}^{N_{2}} Q_{2,k} = 0
  \end{equation}
\end{subequations}
where $ G(w) $ and $ H(w) $ are approximations of $ g(w) $ and $ h(w) $, respectively; $ W_{1,k} $ and $ W_{2,k} $ are {\emph{charge points}} distributed outside of boundaries $ \partial{\bm{D}}^{\prime}_{1} $ and $ \partial{\bm{D}}^{\prime}_{2} $, respectively; $ Q_{1,k} $ and $ Q_{2,k} $ are {\emph{charges}} to be determined with necessary {\emph{collocation conditons}}. Clearly, Eq. (\ref{eq:3.19}) meet the other rotation normalization requirement \deleted{of requirement (1)} in Eq. (\ref{eq:3.18b}). The proof of necessity of Eq. (\ref{eq:3.20}) can be found in Ref \cite{okano2003numerical}, and would not be repeated here.

Similar to deep tunnel, the forward conformal mapping of shallow tunnel using Charge Simulation Method can be determined in the following scheme:

(1) Sequentially select $ N_{1}^{\prime} $ "suitable" {\emph{collocation points}} along the exterior boundary $ \partial{\bm{D}}^{\prime}_{1} $ (counter-clockwisely), and denote these $ N_{1}^{\prime} $ {\emph{collocation points}} as $ w_{1,i} $ ($ i = 1,2,3,\cdots,N_{1}^{\prime} $). Similarly, sequentially select $ N_{2} $ "suitable" {\emph{collocation points}} along the interior boundary $ \partial{\bm{D}}_{2} $ (clockwisely), and denote these $ N_{2} $ {\emph{collocation points}} as $ z_{2,j} $ ($ j = 1,2,3,\cdots,N_{2} $). The reason that $ N_{1}^{\prime} $ {\emph{collocation points}} are selected along boundary $ \partial{\bm{D}}^{\prime}_{1} $ of interval mapping plane $ {\bm{D}}^{\prime} $ (Fig. \ref{fig:2}b-2), instead of along boundary $ \partial{\bm{D}}_{1} $ of physical plane $ {\bm{D}} $ (Fig. \ref{fig:2}b-1), is that $ N_{1}^{\prime} $ finite {\emph{collocation points}} can be easily distributed along the former circular bounadary, but is impossible to be distributed along the infinite horizontal ground surface. In other words, finite points can not be distributed along an infinite line.

(2) For each {\emph{collocation point}} along exterior boundary $ \partial{\bm{D}}^{\prime}_{1} $ and and interior boundary $ \partial{\bm{D}}^{\prime}_{2} $, the corresponding {\emph{charge point}} is selected as
\begin{subequations}
  \label{eq:3.21}
  \begin{equation}
    \label{eq:3.21a}
    \left\{
      \begin{aligned}
        & W_{1,k} = w_{1,k} + K_{1} \cdot H_{1,k} \cdot \exp({\rm{i}}\varTheta_{1,k}) \\
        & H_{1,k} = \frac{1}{2} \left( |w_{1,k+1}-w_{1,k}| + |w_{1,k}-w_{1,k-1}| \right) \\
        & \varTheta_{1,k} = \arg\left( w_{1,k+1} - w_{1,k-1} \right) - \frac{\pi}{2} \\
      \end{aligned}
    \right.
    , \quad k = 1,2,3,\cdots,N_{1}^{\prime}
  \end{equation}
  \begin{equation}
    \label{eq:3.21b}
    \left\{
      \begin{aligned}
        & W_{2,k} = w(z_{2,k}) + K_{2} \cdot H_{2,k} \cdot \exp({\rm{i}}\varTheta_{2,k}) \\
        & H_{2,k} = \frac{1}{2} \left[ |w(z_{2,k+1})-w(z_{2,k})| + |w(z_{2,k})-w(z_{2,k-1})| \right] \\
        & \varTheta_{2,k} = \arg\left[ w(z_{2,k+1}) - w(z_{2,k-1}) \right] - \frac{\pi}{2} \\
      \end{aligned}
    \right.
    , \quad k = 1,2,3,\cdots,N_{2}
  \end{equation}
\end{subequations}
where $ K_{1} $ and $ K_{2} $ denote {\emph{assigment factors}} to adjust the distance between {\emph{collocation point}} and its corresponding {\emph{charge point}}.

(3) With Eq. (\ref{eq:3.19}), Eq. (\ref{eq:3.17'}) can be rewritten as
\begin{equation}
  \label{eq:3.22}
  \sum\limits_{k=1}^{N_{1}^{\prime}} Q_{1,k}\ln\frac{w-W_{1,k}}{w_{\beta}-W_{1,k}} + \sum\limits_{k=1}^{N_{2}} Q_{2,k}\ln\frac{w-W_{2,k}}{w_{\beta}-W_{2,k}} - \ln\zeta = \ln\frac{w-w_{c2}}{w_{\beta}-w_{c2}}
\end{equation}
With the $ N_{1}^{\prime} + N_{2} $ {\emph{collocation points}} and corresponding {\emph{charge points}}, the following $ N_{1}^{\prime} + N_{2} $ {\emph{collocation conditions}} can be established according the real part of Eq. (\ref{eq:3.22}) as
\begin{subequations}
  \label{eq:3.23}
  \begin{equation}
    \label{eq:3.23a}
    \sum\limits_{k=1}^{N_{1}^{\prime}} Q_{1,k}\ln\left|\frac{w_{1,i}-W_{1,k}}{w_{\beta}-W_{1,k}}\right| + \sum\limits_{k=1}^{N_{2}} Q_{2,k}\ln\left|\frac{w_{1,i}-W_{2,k}}{w_{\beta}-W_{2,k}}\right| - \ln r_{o} = -\ln\left|\frac{w_{1,i}-w_{c2}}{w_{\beta}-w_{c2}}\right|, \quad i = 1,2,3,\cdots,N_{1}^{\prime}
  \end{equation}
  \begin{equation}
    \label{eq:3.23b}
    \sum\limits_{k=1}^{N_{1}^{\prime}} Q_{1,k}\ln\left|\frac{w(z_{2,j})-W_{1,k}}{w_{\beta}-W_{1,k}}\right| + \sum\limits_{k=1}^{N_{2}} Q_{2,k}\ln\left|\frac{w(z_{2,j})-W_{2,k}}{w_{\beta}-W_{2,k}}\right| - \ln r_{i} = -\ln\left|\frac{w(z_{2,j})-w_{c2}}{w_{\beta}-w_{c2}}\right|, \quad j = 1,2,3,\cdots,N_{2}
  \end{equation}
\end{subequations}
where $ r_{o} $ and $ r_{i} $ denote the outer and inner radii of the mapping unit annulus $ \overline{\bm{d}} $ in Fig. \ref{fig:2}a-3, and $ r_{o} $ should be approximate to unity. Eqs. (\ref{eq:3.20}) and (\ref{eq:3.23}) form a well-defined simultaneous linear system containing $ N_{1}^{\prime} + N_{2} + 2 $ real variables and linear equations.

(4) With the unique set of  $ Q_{1,k} $, $ Q_{2,k} $, $ \ln r_{o} $, and $ \ln r_{i} $, the forward conformal mapping in Eq. (\ref{eq:3.17'}) can be obtained as
\begin{equation}
  \label{eq:3.24}
  \zeta = \zeta(w) = \frac{w-w_{c2}}{w_{\beta}-w_{c2}} \cdot \exp\left( \sum\limits_{k=1}^{N_{1}^{\prime}} Q_{1,k}\ln\frac{w-W_{1,k}}{w_{\beta}-W_{1,k}} + \sum\limits_{k=1}^{N_{2}} Q_{2,k}\ln\frac{w-W_{2,k}}{w_{\beta}-W_{2,k}} \right), \quad w \in \overline{\bm{D}}^{\prime}
\end{equation}

With the forward conformal mapping, the backward conformal mapping can be subsequently solved using the Complex Dipole Simulation Method (CDSM) \cite{sakakibara2020bidirectional} in the following formation as
\begin{equation}
  \label{eq:3.25}
  w = w(\zeta) = \sum\limits_{k=1}^{N_{1}^{\prime}} \frac{q_{1,k}}{\zeta-\zeta_{1,k}} + \sum\limits_{k=1}^{N_{2}} \frac{q_{2,k}}{\zeta-\zeta_{2,k}}, \quad \zeta \in \overline{\bm{d}}
\end{equation}
where $ \zeta_{1,k} $ and $ \zeta_{2,k} $ denote the {\emph{charge points}} distributed outside of boundary $ \partial{\bm{d}}_{1} $ and $ \partial{\bm{d}}_{2} $ of Fig. \ref{fig:2}a-3, respectively. The {\emph{charge points}} are determined using the Amano's manner \cite{sakakibara2020bidirectional} as
\begin{subequations}
  \label{eq:3.26}
  \begin{equation}
    \label{eq:3.26a}
    \left\{
      \begin{aligned}
        & \zeta_{1,k} = \zeta(w_{1,k}) + k_{1} \cdot h_{1,k} \cdot \exp({\rm{i}}\vartheta_{1,k}) \\
        & h_{1,k} = \frac{1}{2}\left( |\zeta(w_{1,k+1})-\zeta(w_{1,k})| + |\zeta_{1,k}-\zeta_{1,k-1}| \right) \\
        & \vartheta_{1,k} = \arg\left[ \zeta(w_{1,k+1}) - \zeta(w_{1,k-1}) \right] - \frac{\pi}{2} \\
      \end{aligned}
    \right.
    , \quad k = 1,2,3,\cdots,N_{1}
  \end{equation}
  \begin{equation}
    \label{eq:3.26b}
    \left\{
      \begin{aligned}
        & \zeta_{2,k} = \zeta(w_{2,k}) + k_{2} \cdot h_{2,k} \cdot \exp({\rm{i}}\vartheta_{2,k}) \\
        & h_{2,k} = \frac{1}{2}\left( |\zeta(w_{2,k+1})-\zeta(w_{2,k})| + |\zeta_{2,k}-\zeta_{2,k-1}| \right) \\
        & \vartheta_{2,k} = \arg\left[ \zeta(w_{2,k+1}) - \zeta(w_{2,k-1}) \right] - \frac{\pi}{2} \\
      \end{aligned}
    \right.
    , \quad k = 1,2,3,\cdots,N_{2}
  \end{equation}
\end{subequations}
where $ k_{1} $ and $ k_{2} $ are also called {\emph{assignment factors}} to adjust the distance between {\emph{collocation point}} and corresponding {\emph{charge point}}. Eq. (\ref{eq:3.26}) indicates that the {\emph{charge points}} $ \zeta_{1,k} $ and $ \zeta_{2,k} $ are determined by the {\emph{collocation points}} $ \zeta(w_{1,k}) $ and $ \zeta(w_{2,k}) $, respectively, which are computed via the forward conformal maping in Eq. (\ref{eq:3.24}).

With the {\emph{collocation points}} computed via Eq. (\ref{eq:3.24}) and the {\emph{charge points}} computed via Eq. (\ref{eq:3.26}), the {\emph{charges}} in Eq. (\ref{eq:3.25}) can be solved by
\begin{subequations}
  \label{eq:3.27}
  \begin{equation}
    \label{eq:3.27a}
    \sum\limits_{k=1}^{N_{1}^{\prime}} \frac{q_{1,k}}{\zeta(w_{1,i})-\zeta_{1,k}} + \sum\limits_{k=1}^{N_{2}} \frac{q_{2,k}}{\zeta(w_{1,i})-\zeta_{2,k}} = w_{1,i}, \quad i = 1,2,3,\cdots,N_{1}
  \end{equation}
  \begin{equation}
    \label{eq:3.27b}
    \sum\limits_{k=1}^{N_{1}^{\prime}} \frac{q_{1,k}}{\zeta(w_{2,j})-\zeta_{1,k}} + \sum\limits_{k=1}^{N_{2}} \frac{q_{2,k}}{\zeta(w_{2,j})-\zeta_{2,k}} = w_{2,j}, \quad j = 1,2,3,\cdots,N_{2}
  \end{equation}
\end{subequations}
Eq.(\ref{eq:3.27}) contains $ N_{1}^{\prime} + N_{2} $ complex variables and complex linear equations, and would lead to unique solution of $ q_{1,k} $ and $ q_{2,k} $. The backward conformal mapping can be obtained.

Till now, both steps of forward and backward conformal mappings connecting the physical lower half plane containg a cavity of arbitrary shape in Fig. \ref{fig:2}b-1 and the mapping unit annulus in Fig. \ref{fig:2}b-3 are established. The fraction expressions in Eqs. (\ref{eq:3.14b}) and (\ref{eq:3.25}) formally guarantee the analyticity of the stepwise backward conformal mapping in Eq. (\ref{eq:3.16b}).

\section{Numerical strategies}
\label{sec:Numerical results}

\subsection{Mapping single-valuedness \added{and solution flowcharts}}
\label{sub:Mapping single-valuedness}

\begin{figure}[htb]
  \centering
  \begin{tabular}{c}
    (a) Deep tunnel (One-step mapping) \\
    \begin{tikzpicture}[node distance = 1cm]
      \tikzstyle{every node} = [scale = 0.9]
      \node[draw, rounded corners, align = center] (start) {Forward mapping starts \\ Eqs. (\ref{eq:4.1}) and (\ref{eq:4.2})};
      \node[draw, right = of start, align = center] (step 1) {{\emph{Collocation points}} \\ $ z_{i} (i = 1,2,3,\cdots,N) $};
      \node[draw, right = 1.2cm of step 1, align = center] (step 2) {{\emph{Charge points}} \\ $ Z_{k} (k = 1,2,3,\cdots,N) $};
      \node[draw, right = 1.2cm of step 2, align = center] (step 3) {Coefficients $ \varGamma $ and \\ $ Q_{k} (k = 1,2,3,\cdots,N) $};
      \node[draw, right = of step 3, align = center] (step 4) {Eq. \ref{eq:3.9'}};
      \node[below = 0.5cm of step 1] (step mid) {$ \; $};
      \node[draw, below = 1.2cm of step 1, align = center] (step 5) {Mapping points \\ $ \zeta(z_{i}) (i = 1,2,3,\cdots,N) $};
      \node[draw, right = 1.2cm of step 5, align = center] (step 6) {{\emph{Charge points}} \\ $ \zeta_{k} (k = 1,2,3,\cdots,N) $};
      \node[draw, right = 1.2cm of step 6, align = center] (step 7) {Complex coefficients \\ $ q_{k} (k = 1,2,3,\cdots,N) $};
      \node[draw, rounded corners, right = of step 7, align = center] (step 8) {Eq. (\ref{eq:3.12a'})};
      \draw [->] (start) -- (step 1);
      \draw [->] (step 1) -- node [align = center] {Eq. (\ref{eq:3.10}) \\ $ K_{0} $} (step 2);
      \draw [->] (step 2) -- node [align = center] {Eq. (\ref{eq:3.11}) \\ $ z_{c1} $} (step 3);
      \draw [->] (step 3) -- node [below, align = center] {Forward \\ mapping \\ ends} (step 4);
      \draw [->, dashed] (step 4) -- (step 4 |- step mid) -- (step mid);
      \draw [->] (step 1) -- node [left, align = right] {Backward \\ mapping \\ starts} (step 5);
      \draw [->] (step 5) -- node [above] {Eq. (\ref{eq:3.12b'})} (step 6);
      \draw [->] (step 6) -- node [above] {Eq. (\ref{eq:3.13'})} (step 7);
      \draw [->] (step 7) -- node [above, align = center] {Backward \\ mapping \\ ends} (step 8);
    \end{tikzpicture}\\
    \\
    (b) Shallow tunnel (Two-step mapping) \\
    \begin{tikzpicture}[node distance = 1cm]
      \tikzstyle {every node} = [scale = 0.9]
      \node [draw, rounded corners, align = center, minimum width = 3cm, minimum height = 1cm] (start) {Forward mapping starts \\ using selection strategy \\ Eqs. (\ref{eq:4.1}) and (\ref{eq:4.2})};
      \node [draw, right = of start, align = center, minimum width = 3cm, minimum height = 1cm] (step 0-2) {{\emph{Collocation points}} \\ in physical plane \\ $ z_{2,j} (j = 1,2,3,\cdots,N_{2}) $};
      \node [draw, below = of start, align = center, minimum width = 3cm, minimum height = 1cm] (step 1-1) {{\emph{Collocation points}} \\ in mid-mapping plane \\ $ w_{1,i} (i = 1,2,3,\cdots,N_{1}^{\prime}) $ \\ \textcolor{red}{Block 1}};
      \node [draw, below = of step 0-2, align = center, minimum width = 3cm, minimum height = 1cm] (step 1-2) {Mapping points \\ in mid-mapping plane \\ $ w_{2,j} (j = 1,2,3,\cdots,N_{2}) $ \\ \textcolor{red}{Block 2}};
      \node [draw, below = of step 1-1, align = center, minimum width = 3cm, minimum height = 1cm] (step 2-1) {{\emph{Charge points}} \\ $ W_{1,k} (k = 1,2,3,\cdots,N_{1}^{\prime}) $};
      \node [draw, right = of step 2-1, align = center, minimum width = 3cm, minimum height = 1cm] (step 2-2) {{\emph{Charge points}} \\ $ W_{2,k} (k = 1,2,3,\cdots,N_{2}) $};
      \node [right = 0.5cm of step 2-1] (step mid1) {\;};
      \node [draw, below = 1.5cm of step mid1, align = center, minimum width = 3cm, minimum height = 1cm] (step 3) {Coefficients $ \ln r_{o} $, $ \ln r_{i} $ \\ $ Q_{1,k} (k = 1,2,3,\cdots,N_{1}^{\prime}) $ \\ $ Q_{2,k} (k = 1,2,3,\cdots,N_{2}) $};
      \node [draw, below = of step 3, align = center, minimum width = 3cm, minimum height = 1cm] (step 4) {Step 2 of forward mapping \\ Eq. (\ref{eq:3.24'})};
      \node [draw, right = 2cm of step 4, align = center, minimum width = 3cm, minimum height = 1cm] (step 5-1) {Mapping points \\ $ \zeta(w_{1,i}) (i = 1,2,3,\cdots,N_{1}^{\prime}) $};
      \node [draw, right = of step 5-1, align = center, minimum width = 3cm, minimum height = 1cm] (step 5-2) {Mapping points \\ $ \zeta(w_{2,j}) (j = 1,2,3,\cdots,N_{2}) $};
      \node [below = 0.3cm of step 5-2] (step mid2) {\;};
      \node [draw, above = of step 5-1, align = center, minimum width = 3cm, minimum height = 1cm] (step 6-1) {{\emph{Charge points}} \\ $ \zeta_{1,k} (k = 1,2,3,\cdots,N_{1}^{\prime}) $};
      \node [draw, above = of step 5-2, align = center, minimum width = 3cm, minimum height = 1cm] (step 6-2) {{\emph{Charge points}} \\ $ \zeta_{2,k} (k = 1,2,3,\cdots,N_{2}) $};
      \node [right = 0.5cm of step 6-1] (step mid3) {\;};
      \node [draw, above = 1.5cm of step mid3, align = center, minimum width = 3cm, minimum height = 1cm] (step 7) {Complex coefficients \\ $ q_{1,k} (k = 1,2,3,\cdots,N_{1}^{\prime}) $ \\ $ q_{2,k} (k = 1,2,3,\cdots,N_{2}) $};
      \node [draw, above = of step 7, align = center, minimum width = 3cm, minimum height = 1cm] (step 8) {Step 2 of backward mapping \\ Eq. (\ref{eq:3.25})};
      \node [draw, rounded corners, above = of step 8, align = center, minimum width = 3cm, minimum height = 1cm] (step 9) {Step 1 of backward mapping \\ Eq. (\ref{eq:3.14b})};
      \node [right = 0.5cm of start] (mid1) {\;};
      \node [above = 0.5cm of mid1] (forward) {Forward conformal mapping};
      \node [right = 3.7cm of forward] (backward) {Backward conformal mapping};
      \draw [->] (start) -- (step 0-2);
      \draw [->] (start) -- (step 1-1);
      \draw [->] (step 0-2) -- node [left, align = right] {Step 1 of \\ forward mapping \\ Eq. (\ref{eq:3.14a})} (step 1-2);
      \draw [->] (step 1-1) -- node [left, align = right] {Eq. (\ref{eq:3.21a})} (step 2-1);
      \draw [->] (step 1-1) -- node [right, align = right] {$ K_{1} $} (step 2-1);
      \draw [->] (step 1-2) -- node [left, align = right] {Eq. (\ref{eq:3.21b})} (step 2-2);
      \draw [->] (step 1-2) -- node [right, align = right] {$ K_{2} $} (step 2-2);
      \draw [->] (step 2-1) -- (step 2-1 -| step 3) -- node [left, align = right] {Eqs. (\ref{eq:3.20}) \\ and (\ref{eq:3.23})} (step 3);
      \draw [->] (step 2-2) -- (step 2-2 -| step 3) -- node [right, align = left] {$ w_{\beta} $ \\ $ w_{c2} $} (step 3);
      \draw [->] (step 3) -- (step 4);
      \draw [->, dashed, red] (step 4) -- node [above, red] {+Block 1} (step 5-1);
      \draw [->, dashed, red] (step 4) -- (step 4 |- step mid2) -- node [above, red] {+Block 2} (step mid2) -- (step 5-2);
      \draw [->] (step 5-1) -- node [left, align = right] {Eq. (\ref{eq:3.26a})} (step 6-1);
      \draw [->] (step 5-1) -- node [right, align = left] {$ k_{1} $} (step 6-1);
      \draw [->] (step 5-2) -- node [left, align = right] {Eq. (\ref{eq:3.26b})} (step 6-2);
      \draw [->] (step 5-2) -- node [right, align = left] {$ k_{2} $} (step 6-2);
      \draw [->] (step 6-1) -- (step 6-1 -| step 7) -- node [left] {Eq. (\ref{eq:3.27})} (step 7);
      \draw [->] (step 6-2) -- (step 6-2 -| step 7) -- (step 7);
      \draw [->] (step 7) -- (step 8);
      \draw [->] (step 8) -- (step 9);
    \end{tikzpicture}
  \end{tabular}
  \caption{\added{Flowcharts of two bidirectional conformal mappings}}
  \label{fig:4}
\end{figure}
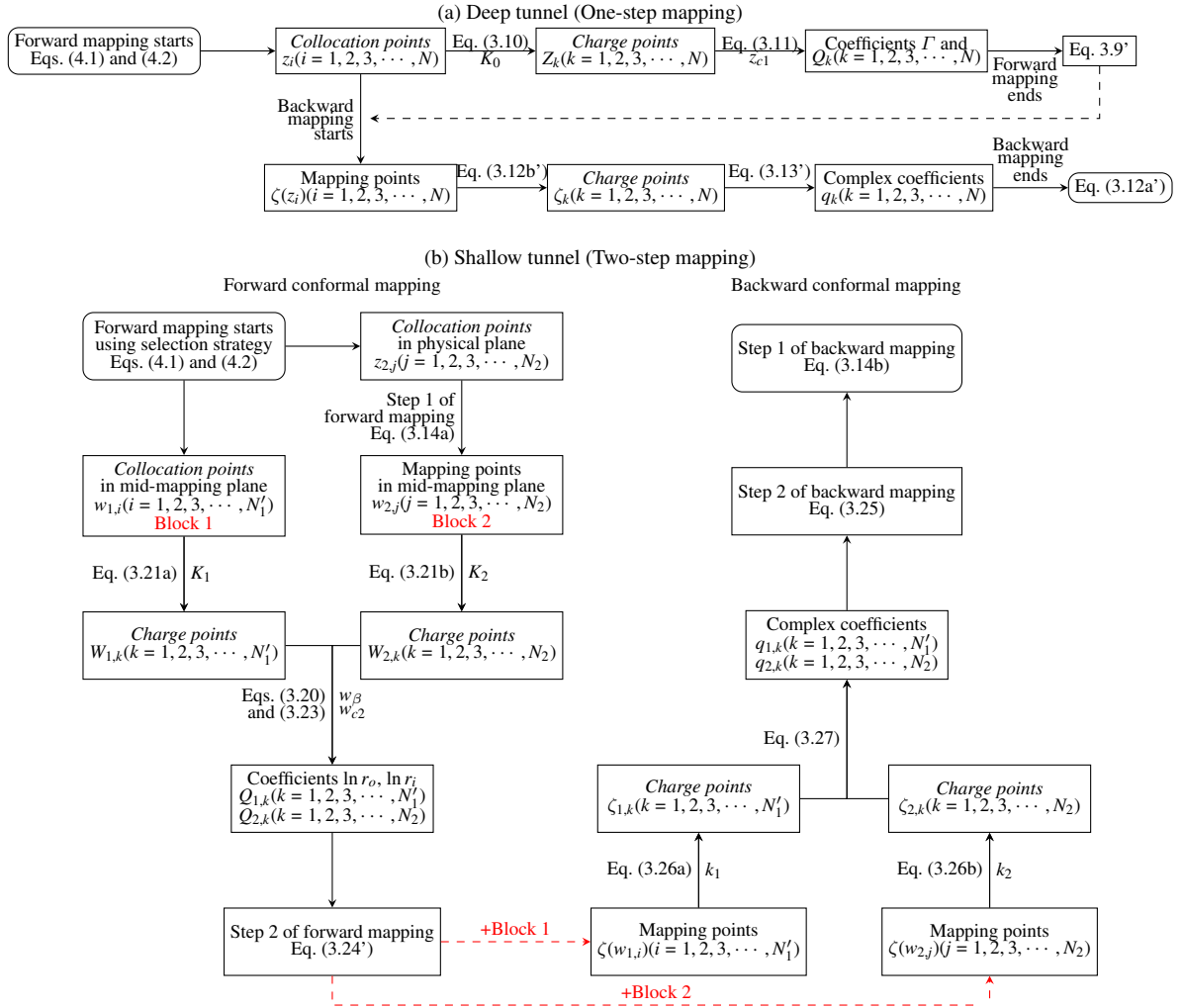

In numerical computation, the bidirectional conformal mapping should be performed on coding languages, such as \texttt{Python}, \texttt{Fortran}, \texttt{Julia}, \texttt{C/C++}, \texttt{Matlab}, \texttt{Mathematica}, and so on. All these coding languages share a similar stipulation that the principle augument ranges within $ (-\pi,\pi] $ in complex variable under natural logarithmic computation \texttt{ln}, which is not yet considered in our previous theoretical deduction. Therefore, in numerical computation, all the natural logarithmic items should be modified to guarantee that the argument should be within $ (-\pi,\pi] $ to avoid possible multi-valuedness. To achieve such a goal, we only need to modify Eqs. (\ref{eq:3.9}) and (\ref{eq:3.24}) as
\begin{equation}
  \label{eq:3.9'}
  \tag{3.9'}
  \zeta(z) = (z-z_{c1}) \cdot \exp\left[ \varGamma + \sum\limits_{k=1}^{N} Q_{k}\ln\frac{z-Z_{k}}{z-z_{c1}} \right], \quad z \in \overline{\bm{D}}
\end{equation}
\begin{equation}
  \label{eq:3.24'}
  \tag{3.24'}
  \zeta(w) = \frac{w-w_{c2}}{w_{\beta}-w_{c2}} \cdot \exp\left[ 
    \begin{aligned}
      & \sum\limits_{k=1}^{N_{1}} Q_{1,k}\left(\ln\frac{w-W_{1,k}}{w_{1,0}-W_{1,k}}-\ln\frac{w_{\beta}-W_{1,k}}{w_{1,0}-W_{1,k}}\right) \\
      & + \sum\limits_{k=1}^{N_{2}} Q_{2,k}\left(\ln\frac{w-W_{2,k}}{w-w_{c2}}-\ln\frac{w_{\beta}-W_{2,k}}{w_{\beta}-w_{c2}}\right) 
    \end{aligned}
  \right], \quad w \in \overline{\bm{D}}^{\prime}
\end{equation}
where $ w_{1,0} $ is an arbitrary point within domain $ {\bm{D}}^{\prime} $. The modifications above use Eq. (\ref{eq:3.4}) and (\ref{eq:3.20b}), respectively. With the modifications of Eqs. (\ref{eq:3.9'}) and (\ref{eq:3.24'}), the arguments of the natural logarithmic items would be always constrained within the principle range $ (-\pi,\pi] $, so that the mapping results are always single-valued. \added{With Eqs. (\ref{eq:3.9'}) and (\ref{eq:3.24'}), the solution procedures of bidirectional conformal mapping for deep and shallow tunnels are illustrated in the flowcharts in Figs. \ref{fig:4}a and \ref{fig:4}b, respectively.}

\subsection{Selection strategy of collocation points}
\label{sub:Selection strategy of collocation points}

The bidirectional conformal mappings for both deep and shallow tunnels are derived in formalized schemes, except for the "suitable" {\emph{collocation points}}, which are mathematically vague in our previous description, and should be explained explicitly. Similar to other numerical method, such as boundary element method, finite difference method, finite element method, and finite volume method, the selection of collocation points of Charge Simulation Method is more or less emprical, and may require trial computation for specific cases. However, there are several recommended fundamental schemes to follow.

Eqs. (\ref{eq:3.10}), (\ref{eq:3.21a}), and (\ref{eq:3.21b}) show that $ N $, $ N_{1} $ and $ N_{2} $ {\emph{collocation points}} are extracted to simulate boundaries $ \partial{\bm{D}} $ of deep tunnel in physical plane $ {\bm{D}} $ (Fig. \ref{fig:2}a-1), $ \partial{\bm{D}}^{\prime}_{1} $ of shallow tunnel in interval mapping plane $ {\bm{D}}^{\prime} $ (Fig. \ref{fig:2}b-2), and $ \partial{\bm{D}}_{2} $ of shallow tunnel in physical plane $ {\bm{D}} $ (Fig. \ref{fig:2}a-1), respectively. The angle variations and relative distances of adjacent {\emph{collocation points}} for all three boundaries can be defined as
\begin{subequations}
  \label{eq:4.1}
  \begin{equation}
    \label{eq:4.1a}
    \left\{
      \begin{aligned}
        & \vartheta_{k} = \arg(z_{k+1}-z_{k}), \quad k = 1,2,3,\cdots,N \\
        & \vartheta_{1,k} = \arg(w_{1,k+1}-w_{1,k}), \quad k = 1,2,3,\cdots,N_{1} \\
        & \vartheta_{2,k} = \arg(z_{2,k+1}-z_{2,k}), \quad k = 1,2,3,\cdots,N_{2} \\
      \end{aligned}
    \right.
  \end{equation}
  \begin{equation}
    \label{eq:4.1b}
    \left\{
      \begin{aligned}
        & \delta_{k} = \frac{| z_{k+1} - z_{k} |}{\sum\limits_{l=1}^{N} | z_{l} - z_{l-1} |}, \quad k = 1,2,3,\cdots,N \\
        & \delta_{1,k} = \frac{| w_{1,k+1} - w_{1,k} |}{\sum\limits_{l=1}^{N} | w_{1,l} - w_{1,l-1} |}, \quad k = 1,2,3,\cdots,N_{1} \\
        & \delta_{2,k} = \frac{| z_{2,k+1} - z_{2,k} |}{\sum\limits_{l=1}^{N} | z_{2,l} - z_{2,l-1} |}, \quad k = 1,2,3,\cdots,N_{2} \\
      \end{aligned}
    \right.
  \end{equation}
\end{subequations}
where $ \vartheta_{k} $, $ \vartheta_{1,k} $, and $ \vartheta_{2,k} $ denote the angle variations of adjacent {\emph{collocation points}} along boundaries $ \partial{\bm{D}} $, $ \partial{\bm{D}}^{\prime}_{1} $, and $ \partial{\bm{D}}_{2} $, respectively; $ \delta_{k} $, $ \delta_{1,k} $, and $ \delta_{2,k} $ denote relative distances of adjacent {\emph{collocation points}} along boundaries  $ \partial{\bm{D}} $, $ \partial{\bm{D}}^{\prime}_{1} $, and $ \partial{\bm{D}}_{2} $, respectively.

These three boundaries consisting of {\emph{collocation points}} should be analytic, so that the Charge Simulation Method can be applicable. To guarantee the analyticity of all three boundaries, the following constraints of angle variations and relative distances of adjacent {\emph{collocation points}} are recommended as:
\begin{subequations}
  \label{eq:4.2}
  \begin{equation}
    \label{eq:4.2a}
    \max\left\{ \max_{1 \leq k \leq N}|\vartheta_{k}|, \max_{1 \leq k \leq N_{1}} |\vartheta_{1,k}|, \max_{1 \leq k \leq N_{2}} |\vartheta_{2,k}| \right\} \leq 10^{\circ}
  \end{equation}
  \begin{equation}
    \label{eq:4.2b}
    \max\left\{ \max_{1 \leq k \leq N} \delta_{k}, \max_{1 \leq k \leq N_{1}} \delta_{1,k}, \max_{1 \leq k \leq N_{2}} \delta_{2,k} \right\} \leq 10^{-2}
  \end{equation}
\end{subequations}
Generally, if the curvature of cavity boundary is more stable, in other word, more geometrically approximate to a circle, the quantity of {\emph{collocation points}} can be correspondingly smaller. Eq. (\ref{eq:4.2}) is summarized via many attemts of trial computation, and is suitable to tunnel boundaries without sharp corners. 

Eq. (\ref{eq:4.2a}) would no longer hold for boundaries with sharp corners (angle variations of $ 90^{\circ} $ in Figs. \ref{fig:2}a-1 and b-1 for instance), which lead to boundary singularities \cite{hough1981use,hough1983integral} and further cause failure of the Charge Simulation Method. Such a failure would not exist in real-world tunnel engineering, since tunnel boundaries with sharp corners are generally not applied to avoid typical and serious stress concentrations, which provides a strong realistic base to simulate sharp corners with smooth and analytic curves.

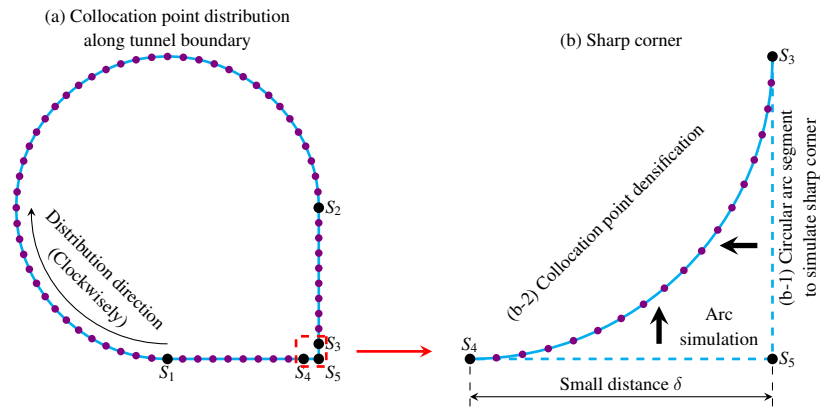
\begin{figure}[htb]
  \centering
  \begin{tikzpicture}
    \tikzstyle{every node} = [scale = 0.7]
    \draw [cyan, line width = 1pt] (2,0) arc [start angle = 0, end angle = 270, radius = 2];
    \draw [cyan, line width = 1pt] (2,0) -- (2,-2) -- (0,-2);
    \foreach \x in {0,1,2,3,...,45} \fill [violet] ({2*cos(\x/45*270)},{2*sin(\x/45*270)}) circle [radius = 0.05];
    \foreach \x in {0,1,2,3,...,10} \fill [violet] ({2},{-\x/10*2}) circle [radius = 0.05];
    \foreach \x in {0,1,2,3,...,10} \fill [violet] ({\x/10*2},{-2}) circle [radius = 0.05];
    \fill (2,0) circle [radius = 0.07];
    \fill (0,-2) circle [radius = 0.07];
    \fill (2,-2) circle [radius = 0.07];
    \fill (1.8,-2) circle [radius = 0.07];
    \fill (2,-1.8) circle [radius = 0.07];
    \node at (0,-2) [below] {$ S_{1} $};
    \node at (2,0) [right] {$ S_{2} $};
    \node at (2,-1.8) [right] {$ S_{3} $};
    \node at (1.8,-2) [below] {$ S_{4} $};
    \node at (2,-2) [below right] {$ S_{5} $};
    \draw [dashed, red, line width = 1pt] (2.1, -1.7) rectangle (1.7,-2.1);
    \draw [->, line width = 1pt, red] (2.5,-1.9) -- (3.5,-1.9);
    \draw [->] (0,-1.8) arc [start angle = -90, end angle = -180, radius = 1.8];
    \node at (-1.2,-1.2) [above, rotate = -45, align = center] {Distribution direction \\ (Clockwisely)};
    \node at (0,2) [above, align = center] {(a) Collocation point distribution \\ along tunnel boundary};
    
    \draw [cyan, line width = 1pt, dashed] (8,2) -- (8,-2) -- (4,-2);
    \draw [cyan, line width = 1pt] (4,-2) arc [start angle = -90, end angle = 0, radius = 4];
    \fill (8,2) circle [radius = 0.07];
    \fill (8,-2) circle [radius = 0.07];
    \fill (4,-2) circle [radius = 0.07];
    \node at (8,2) [right] {$ S_{3} $};
    \node at (4,-2) [above] {$ S_{4} $};
    \node at (8,-2) [right] {$ S_{5} $};
    \foreach \x in {1,2,3,...,17} \fill [violet] ({4+4*cos(-90+\x/18*90)},{2+4*sin(-90+\x/18*90)}) circle [radius = 0.05];
    \draw [->, line width = 2pt] (6.5,-1.8) -- (6.5,-1.3);
    \draw [->, line width = 2pt] (7.8,-0.5) -- (7.3,-0.5);
    \node at (7.9,-1.9) [above left, align = center] {Arc \\ simulation};
    \draw [dashed] (8,-2.1) -- (8,-2.7);
    \draw [dashed] (4,-2.1) -- (4,-2.7);
    \draw [<->] (8,-2.5) -- (4,-2.5);
    \node at (6,-2.5) [above] {Small distance $ \delta $};
    \node at (8,0) [below, rotate = 90, align = center] {(b-1) Circular arc segment \\ to simulate sharp corner};
    \node at ({6-0.2},-0.2) [rotate = 45, align = center] {(b-2) Collocation point densification};
    \node at (6,2) [above] {(b) Sharp corner};
  \end{tikzpicture}
  \caption{Selection scheme of collocation points near sharp corners}
  \label{fig:3}
\end{figure}

Fig. \ref{fig:3} presents one possible sharp corner simulation strategy, and provide a feasible scheme to select collocation points near sharp corners in a densified manner. Fig. \ref{fig:3}a illustrates a clockwise {\emph{collocation point}} distribution (to always keep the remaining geomaterial on left side) along the boundary of the under-break horseshoe cavities for deep and shallow tunnelling in Fig. \ref{fig:2}. The {\emph{collocation points}} are distributed in a relatively uniform manner, as recommended in Eq. (\ref{eq:4.2}), and the curvature of the boundary is relatively stable to 0 or 1, except for the bottom right corner.

The boundary near the sharp corner is enlarged in Fig. \ref{fig:3}b, and a circular arc segment is applied to replace the original straight lines, which is identical to the engineering demand of no sharp corners. The circular arc segment is further divided into smaller straight line segments by densifying more collocation points within, so that the angle variations in Eq. (\ref{eq:4.2a}) would hold to guarantee the analyticity of the circular arc segment. In other words, a sharp corner should always be simulated by a small circular arc segment, in which {\emph{collocation points}} should be densified.

\subsection{Condition number and error estimates}
\label{sub:Condition number and error estimates}

The forward and backward conformal mappings for deep tunnel are computed via real linear system of Eq. (\ref{eq:3.11}) and complex linear system in Eqs. (\ref{eq:3.13}) and (\ref{eq:3.13'}), respectively; while the second-step forward and backward conformal mappings for shallow tunnel are computed via real linear system of Eqs. (\ref{eq:3.20}) and (\ref{eq:3.23}) and complex linear system of Eq. (\ref{eq:3.27}), respectively. All the computations of linear systems should consider numerical stability measured by condition number $ C_{N} $, which denotes the sensitivity of solution vector to perturbations of constant vector, and is an intrinsic property of coefficient matrix. The condition number of coefficient matrix is measured by its 2-norm as
\begin{equation}
  \label{eq:4.3}
  C_{N} = \frac{S_{1}}{S_{-1}}
\end{equation}
where $ S_{1} $ and $ S_{-1} $ denote the largest and smallest singular values of the coefficient matrix, respectively. The condition number of a coefficient matrix should be small for good numerical stability. To be specific, the condition numbers of the coefficient matrices in Eqs. (\ref{eq:3.11}), (\ref{eq:3.13'}), (\ref{eq:3.20}) and (\ref{eq:3.23}), (\ref{eq:3.27}) are denoted by $ C_{Nf}^{d} $, $ C_{Nb}^{d} $, $ C_{Nf}^{s} $, $ C_{Nb}^{s} $, respectively.

Error estimates of the bidirectional conformal mappings of deep and shallow tunnels should be conducted. For deep tunnel, error estimates of forward conformal mapping in Eq. (\ref{eq:3.9}) and backward conformal mapping in Eq. (\ref{eq:3.12}) can be evaluated by
\begin{equation}
  \label{eq:4.4}
  \left\{
    \begin{aligned}
      & \varepsilon_{f} = \max_{1 \leq k \leq N} \left| \zeta(z_{k+\frac{1}{2}}) - 1 \right| \\
      & \varepsilon_{b} = \max_{1 \leq k \leq N} \left| z[\zeta(z_{k+\frac{1}{2}})] - z_{k+\frac{1}{2}} \right| \\
    \end{aligned}
  \right.
  , \quad z_{k+\frac{1}{2}} = \frac{1}{2}(z_{k}+z_{k+1})
\end{equation}
where $ \varepsilon_{f} $ and $ \varepsilon_{b} $ denote the errors of forward and backward conformal mappings, respectively; $ z_{k+\frac{1}{2}} $ denotes the midpoint of {\emph{collocation points}} $ z_{k} $ and $ z_{k+1} $.

For shallow tunnel, two boundaries should apply different error estimate schemes. For the exterior boundary (ground surface), error estimates can only be conducted to the second-step forward conformal mapping in Eq. (\ref{eq:3.24}) and backward conformal mapping in Eq. (\ref{eq:3.25}), since the first-step bidirectional conformal mapping is accurate and {\emph{collocation points}} are selected in the interval mapping annulus in Fig. \ref{fig:2}b-2. The error estimates can be expressed as
\begin{subequations}
  \label{eq:4.5}
  \begin{equation}
    \label{eq:4.5a}
    \left\{
      \begin{aligned}
        & \varepsilon_{fe} = \max_{1 \leq k \leq N_{1}} \left| \zeta(w_{1,k+\frac{1}{2}}) - r_{o} \right| \\
        & \varepsilon_{be} = \max_{1 \leq k \leq N_{1}} \left| w[\zeta(w_{1,k+\frac{1}{2}})] - w_{1,k+\frac{1}{2}} \right| \\
      \end{aligned}
    \right.
    , \quad w_{1,k+\frac{1}{2}} = \frac{1}{2}(w_{1,k}+w_{1,k+1})
  \end{equation}
  where $ \varepsilon_{fe} $ and $ \varepsilon_{be} $ denote the errors of second-step forward and backward conformal mappings of exterior boundary of ground surface, respectively; $ w_{1,k+\frac{1}{2}} $ denotes the midpoint of {\emph{collocation points}} $ w_{1,k} $ and $ w_{1,k+1} $.

For the interior boundary (cavity periphery), the {\emph{collocation points}} should be selected in the physical plane in Fig. \ref{fig:2}b-1, and the error estimates should be conducted to the whole bidirectional conformal mapping. The error estimates of forward conformal mapping in Eq. (\ref{eq:3.16a}) and backward conformal mapping in Eq. (\ref{eq:3.16b}) can be expressed as
  \begin{equation}
    \label{eq:4.5b}
    \left\{
      \begin{aligned}
        & \varepsilon_{fi} = \max_{1 \leq k \leq N_{2}} \left| \zeta[w(z_{2,k+\frac{1}{2}})] - r_{i} \right| \\
        & \varepsilon_{bi} = \max_{1 \leq k \leq N_{2}} \left| z[w(\zeta_{2,k+\frac{1}{2}})] - z_{2,k+\frac{1}{2}} \right| \\
      \end{aligned}
    \right.
    , \quad z_{2,k+\frac{1}{2}} = \frac{1}{2}(z_{2,k}+z_{2,k+1})
  \end{equation}
\end{subequations}
where $ \varepsilon_{fi} $ and $ \varepsilon_{bi} $ denote the errors of forward and backward conformal mappings of interior boundary (tunnel cavity boundary), respectively; $ z_{2,k+\frac{1}{2}} $ denotes the midpoint of {\emph{collocation points}} $ z_{2,k} $ and $ z_{2,k+1} $; $ \zeta_{2,k+\frac{1}{2}} = \zeta[w(z_{2,k+\frac{1}{2}})] $.

\section{Numerical results of bidirectional conformal mapping}
\label{sec:Numerical cases}

In this section, several numerical cases will be conducted to illustrate the bidirectional conformal mapping incorporating Charge Simulation Method in this paper. The benchmark cases for deep and shallow tunnels are the under-break horseshoe cavity shown in Fig. \ref{fig:3}a, in which the cavities consist of three quadrants of a circle ($ \overset{\LARGE{\frown}}{S_{1}S_{2}} $), two straight lines ($ S_{2}S_{3} $ and $ S_{4}S_{1} $), and a sharp right corner ($ \measuredangle S_{3}S_{5}S_{4} $), which is simulated by a small quadrant $ \overset{\LARGE{\frown}}{S_{3}S_{4}} $ in Fig. \ref{fig:3}b. According to Eq. (\ref{eq:4.2}), the three quadrants ($ \overset{\LARGE{\frown}}{S_{1}S_{2}} $) are all uniformly discreted and simulated by $ 3 \cdot n_{1} = 3 \cdot 40 $ collocation points, the two lines ($ S_{2}S_{3} $ and $ S_{4}S_{1} $) are both uniformly discreted and simulated by $ n_{2} = 30 $ collocation points, the small quadrant ($ \overset{\LARGE{\frown}}{S_{3}S_{4}} $) is uniformly discreted and simulated by $ n_{3} = 20 $ collocation points. Therefore, we have $ N = 3n_{1}+2n_{2}+n_{3} $ collocation points for deep tunnel boundary in Eqs. (\ref{eq:3.9'}) and (\ref{eq:3.12a'}), while $ N_{2} = 3n_{1}+2n_{2}+n_{3} $ collocation points for shallow tunnel boundary in Eqs. (\ref{eq:3.24'}) and (\ref{eq:3.25}). Since the exterior boundary of shallow tunnel for Charge Simulation Method (see Fig. \ref{fig:2}b-2) is a unit circle, we can simply select $ N_{1} = 90 $ uniform collocation points. With the selected collocation points, the bidirectional conformal mappings of deep and shallow tunnels incorporating Charge Simulation Method can be obtained.

The numerical cases in this section and in Sections \ref{sub:Results and comparisons with finite element method} and \ref{sec:Second numerical case and comparisons with existing analytical solution} are all computed using programming language \texttt{Fortran} of compiler \texttt{GCC 14.1.1}. The condition number in Eq. (\ref{eq:4.3}) is computed using \texttt{DGESVD} and \texttt{ZGESVD} packages of \texttt{LAPACK 3.12.0} for real and complex coefficient matrix, respectively. The real linear systems in Eqs. (\ref{eq:3.11}), (\ref{eq:3.20}) and (\ref{eq:3.23}) are solved using \texttt{DGESV} package, while the complex linear systems in Eqs. (\ref{eq:3.13}), (\ref{eq:3.27}), and (\ref{eq:7.16}) are solved using \texttt{ZGESV} package. The figures are plotted using \texttt{Gnuplot-6.0}. All source codes of numerical cases in this section and Sections \ref{sub:Results and comparisons with finite element method} and \ref{sec:Second numerical case and comparisons with existing analytical solution} are released in author \textsf{Luobin Lin}'s github repository {\emph{github.com/luobinlin987/conformal-mapping-incorporating-Charge-Simulation-Method}}.

\subsection{Deep tunnel cases}
\label{sec:Deep tunnel cases}

\begin{figure}[htb]
  \centering
  \includegraphics[width = 0.85\textwidth]{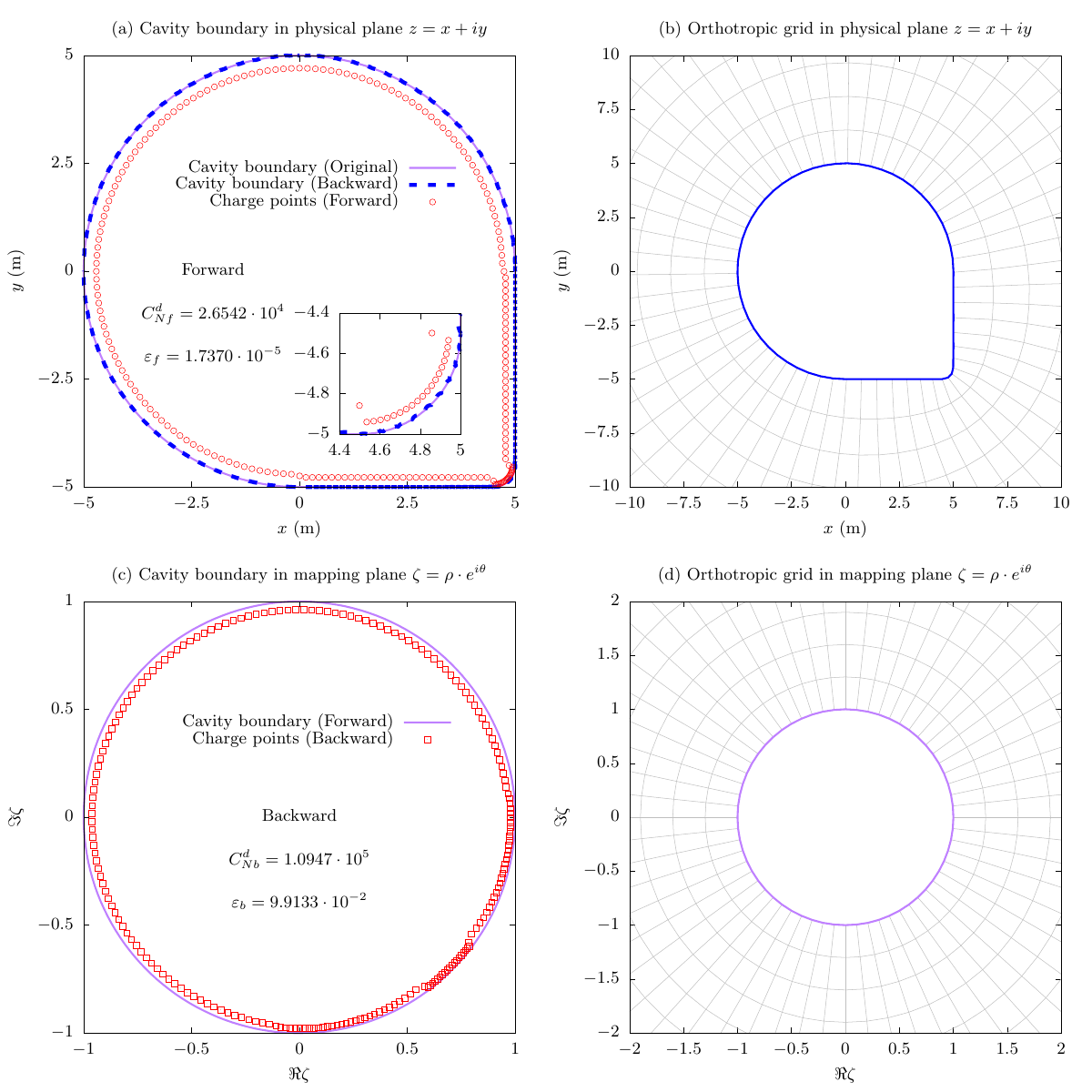}
  \caption{Deep tunnel of under-break horseshoe cavity}
  \label{fig:5}
\end{figure}

The cavity size and original boundary of the deep tunnel are shown in Fig. \ref{fig:5}a, where the coordinates of the {\emph{collocation points}} along the boundary $ \partial{\bm{D}} $ can be obtained analytically as
\begin{equation}
  \label{eq:5.1}
  z_{i} = 
  \left\{
    \begin{aligned}
      & 5 \cdot \exp\left( {\rm{i}} \frac{3}{2}\pi - \frac{i-1}{3n_{1}}\cdot \frac{3}{2}\pi \right), \quad 1 \leq i \leq 3n_{1} \\
      & 5 - {\rm{i}} \frac{i-3n_{1}-1}{n_{2}} \cdot 4.5, \quad 3n_{1} \leq i \leq 3n_{1}+n_{2} \\
      & 4.5-4.5{\rm{i}} + 0.5 \cdot \exp\left( -{\rm{i}}\frac{i-3n_{1}-n_{2}-1}{n_{3}} \cdot \frac{\pi}{2} \right), \quad 3n_{1}+n_{2}+1 \leq i \leq 3n_{1}+n_{2}+n_{3} \\
      & 4.5-5{\rm{i}}-\frac{i-3n_{1}-n_{2}-n_{3}-1}{n_{2}} \cdot 4.5, \quad 3n_{1}+n_{2}+n_{3}+1 \leq i \leq 3n_{1}+2n_{2}+n_{3} \\
    \end{aligned}
  \right.
\end{equation}
Eq. (\ref{eq:5.1}) determines the "suitable" {\emph{collocation points}}, which are the most flexible and empirical part of the Charge Simulation Method, while the rest procedure in Eq. (\ref{eq:3.9})-(\ref{eq:3.13}) is relatively routine.

Substituting the $ N = 3n_{1}+2n_{2}+n_{3} $ collocation points into the bidirectional procedure in Eqs. (\ref{eq:3.9})-(\ref{eq:3.13}) gives Fig. \ref{fig:5}, where $ z_{c1} = 0 $. Fig. \ref{fig:5}a shows the comparison between the original cavity boundary and cavity boundary computed via the backward conformal mapping, and the {\emph{charge points}} are also marked. The small quadrant section is magnified for more details. The condition number $ C_{Nf}^{d} $ and error estimate $ \varepsilon_{f} $ both suggest that the forward conformal mapping is numerically stable and accurate. Fig. \ref{fig:5}c shows the mapping boundary after forward conformal mapping, and {\emph{charge points}} preparing for the backward conformal mapping. The condition number $ C_{Nb}^{d} $ and $ \varepsilon_{b} $ are much larger than those of the forward conformal mapping, but the cavity boundary comparison in Fig. \ref{fig:5}a indicates good backward mapping. Therefore, the backward conformal mapping is still accurate enough. Substituting the orthotropic one in Fig. \ref{fig:5}d into the backward conformal mapping gives the smooth grid in Fig. \ref{fig:5}b, which numerically indicates the analyticity of backward conformal mapping.

To better illustrate the bidirectional conformal mapping of deep tunnel, three more deep tunnel cases of complicated cavity shapes are conducted in Fig. \ref{fig:4}. The boundaries of horseshoe cavity, over-break horseshoe cavity, and square cavity apply $ N = 2n_{1}+4n_{2}+2n_{3} $, $ N = n_{1}+6n_{2}+3n_{3} $, and $ N = 8n_{2}+4n_{3} $ collocation points, respectively. The coordinates of these {\emph{collocation points}} can be analytically given in the same manner of Eq. (\ref{eq:5.1}). Since the coordinate expressions are simple, they are not presented here. The condition numbers and error estimates of forward and backward conformal mappings are all listed within corresponding figures, which are relatively small. The comparisons between original cavity boundaries and the ones computed via corresponding backward conformal mappings show good agreements.

\subsection{Shallow tunnel cases}
\label{sec:Shallow tunnel cases}

For consistency, the size and shape of the shallow tunnel are the same to that in previous section, and the "suitable" {\emph{collocation points}} along cavity boundary $ \partial{\bm{D}}_{2} $ in physical plane $ z = x+{\rm{i}}y $ can be determined by simply shifting the coordinates in Eq. (\ref{eq:5.1}) as
\begin{equation}
  \label{eq:5.2}
  z_{2,j} = z_{j} -z_{c1} + z_{c2}, \quad j = 1,2,3,\cdots,N_{2}
\end{equation}
where $ z_{c2} = -10{\rm{i}} $, $ N_{2} = N $, and the subscript $ i $ is replaced by $ j $ without loss of generality. The {\emph{collocation points}} along the finite ground surface $ {\bm{D}}^{\prime}_{1} $ in the interval mapping plane $ w = u+{\rm{i}}v $ can be analytically given as
\begin{equation}
  \label{eq:5.3}
  w_{1,i} = \exp\left( {\rm{i}} \frac{i-1}{N_{1}} \cdot 2\pi \right), \quad i = 1,2,3,\cdots,N_{1}
\end{equation}

\begin{figure}[htb]
  \centering
  \includegraphics[width = 0.85\textwidth]{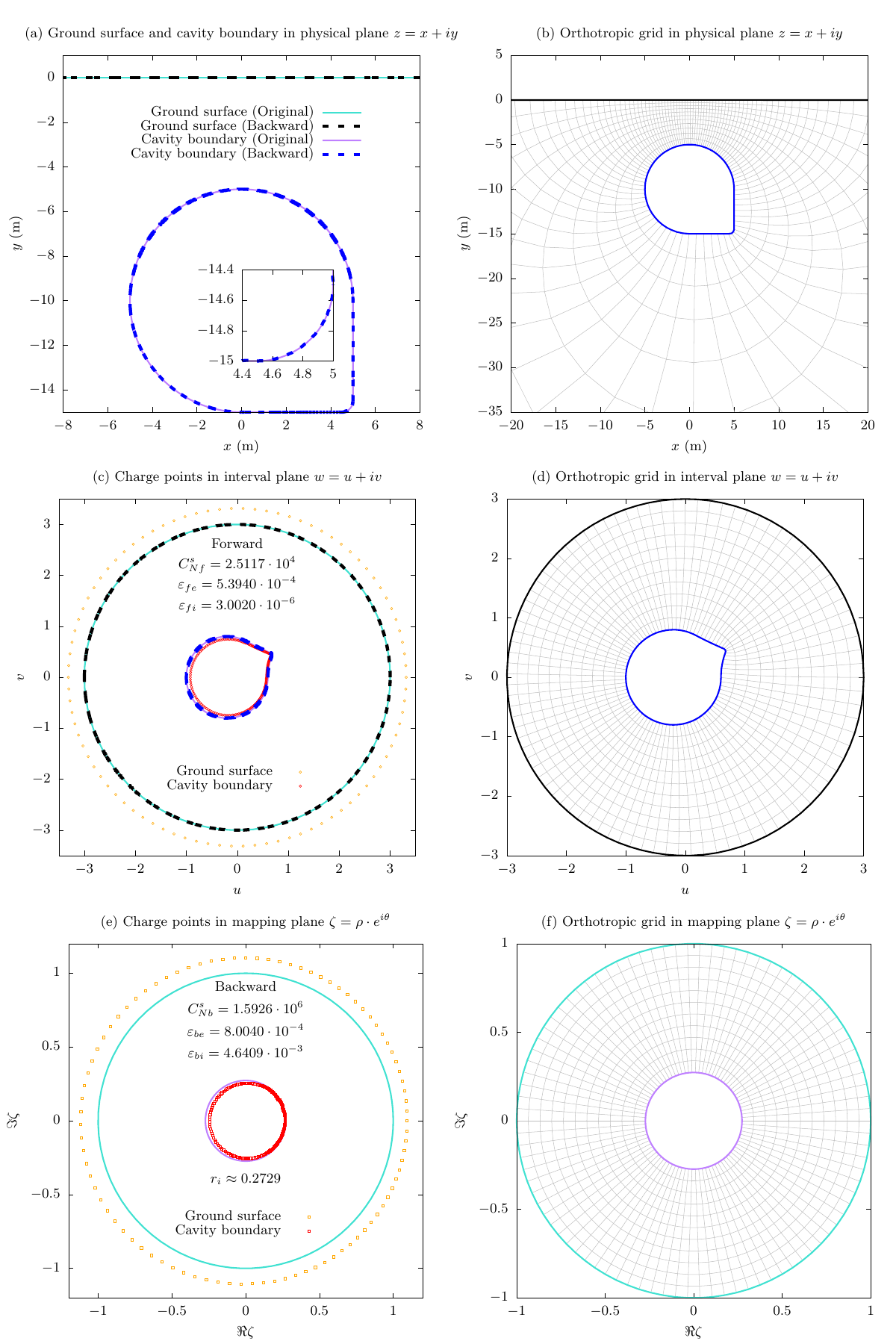}
  \caption{Shallow tunnel of under-break horseshoe cavity}
  \label{fig:7}
\end{figure}

Substituting the {\emph{collocation points}} into the bidirectional conformal mapping of lower half geomaterial containing a cavity of arbitrary shape in Eqs. (\ref{eq:3.21})-(\ref{eq:3.27}) gives Fig. \ref{fig:7}. Figs. \ref{fig:7}a and c show good agreements between both original boundaries and both boundaries computed via backward conformal mapping. The condition numbers and error estimates for forward and backward conformal mappings in Figs. \ref{fig:7}c and e are relatively small, indicating good numerical stability and numerical accuracy. Figs. \ref{fig:7}f, d, and b sequentially show the orthotropic grid of stepwise backward conformal mappings in mapping plane $ \zeta = \rho\cdot{\rm{e}}^{{\rm{i}\theta}} $, interval mapping plane $ w = u+{\rm{i}}v $, and physical plane $ z = x+{\rm{i}}y $, respectively. All the grids are smooth, which also numerically indicates the analycity of the backward conformal mapping.

\begin{figure}[htb]
  \centering
  \includegraphics[width = 0.85\textwidth]{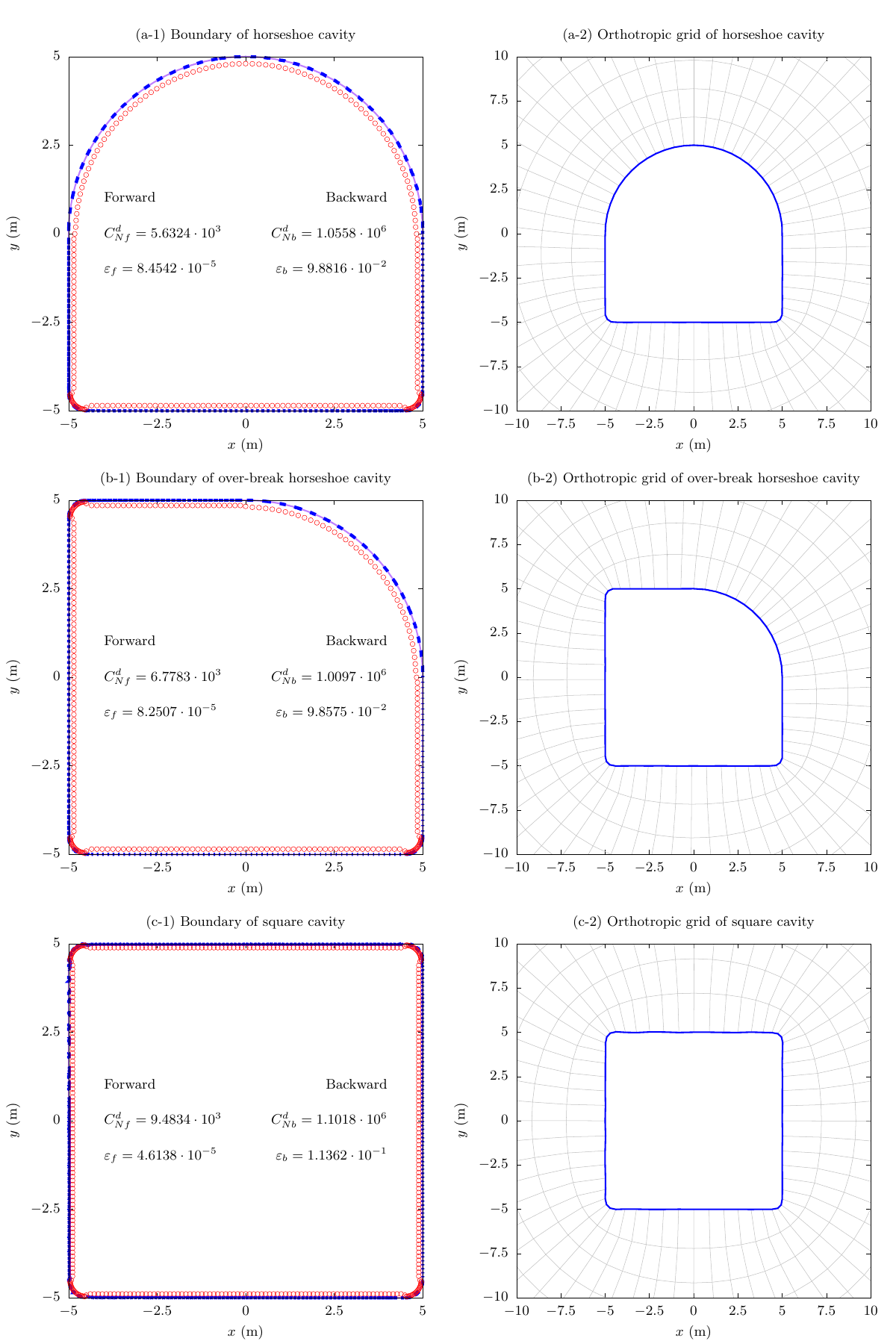}
  \caption{Deep tunnel of three different cavities}
  \label{fig:6}
\end{figure}

\begin{figure}[htb]
  \centering
  \includegraphics[width = 0.85\textwidth]{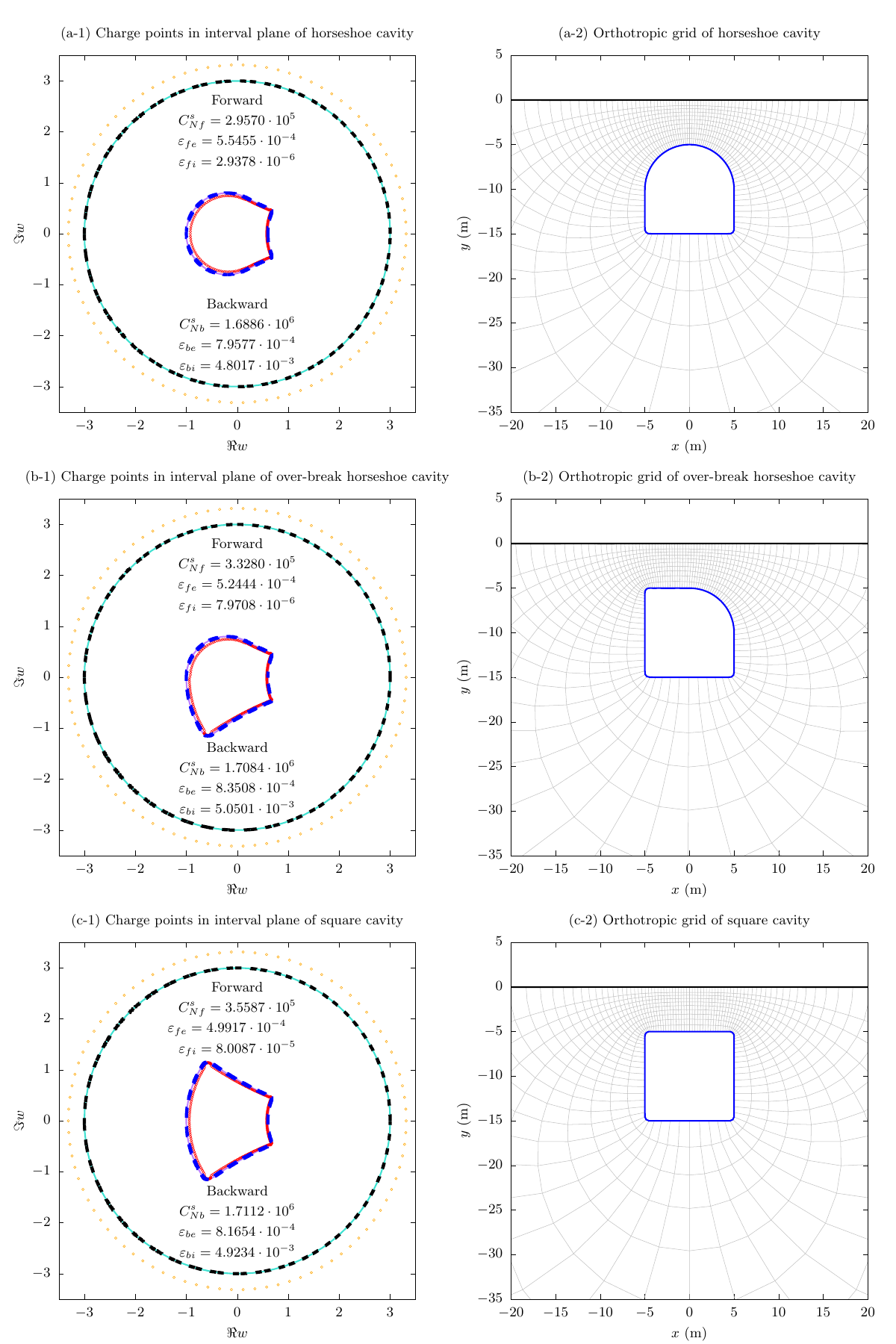}
  \caption{Shallow tunnel of three different cavities}
  \label{fig:8}
\end{figure}

Similar to the deep tunnel cases, three more shallow tunnel cases are shown in Fig. \ref{fig:8}. The boundaries are the same of deep tunnel cases, and the {\emph{collocation point}} quantities of inner boundaries $ \partial{\bm{D}_{2}} $ in physical plane $ z = x+{\rm{i}}y $ take $ N_{2} = 2n_{1}+4n_{2}+2n_{3} $, $ N_{2} = n_{1}+6n_{2}+3n_{3} $, and $ N_{2} = 8n_{2}+4n_{3} $, respectively. These {\emph{collocation point}} coordinates are also simply shifted from the $ N $ {\emph{collocation points}} in Fig. \ref{fig:6} via the same manner in Eq. (\ref{eq:5.2}). The {\emph{collocation points}} of exterior boundaries $ \partial{\bm{D}}^{\prime}_{1} $ in interval mapping plane $ w = u+{\rm{i}}v $ take the same values in Eq. (\ref{eq:5.3}). The small values of condition numbers and error estimates in Figs. \ref{fig:8}a-1, b-1, and c-1, as well as the good agreements between original boundaries and the boundaries computed via backward conformal mappings, all suggest that the stepwise bidirectional conformal mapping of shallow tunnel is numerically stable and accurate. The smooth grids in Figs. \ref{fig:8}a-2, b-2, and c-2 further fortify the analyticty of the backward conformal mapping.

\section{Comparisons with existing conformal mappings}%
\label{sec:Comparisons with existing conformal mappings}

\subsection{Deep tunnel}%
\label{sub:Deep tunnel}

Comparing to exsting conformal mappings, the bidirectional conformal mapping for deep tunnel in Eq. (\ref{eq:3.9'}) makes little accuracy improvement, thus, no comparison of deep tunnel cases should be necessary. However, an extra and minor advantage of the present bidirectional conformal mapping of deep tunnel cases against existing ones should be pointed out that the present conformal mapping is very fast in solving procedure, and takes very low computation cost. To be specific, all the deep and shallow tunnel cases are computed on an old \texttt{thinkpad} laptop \texttt{x61s} (CPU type: \texttt{Core2 DUO L7300} of 1.40 GHz), and the elapsing time of all these cases is alway within 8 seconds. 

In existing conformal mapping of deep tunnels, the classical {\fontencoding{OT2}\selectfont Melentev} method \cite{Muskhelishvili1966} and the improving method by Ye and Ai \cite{ye2023novel-aizhiyong} need iterations; while the least square method by Ma et al. \cite{ma2022numerical} is based on mixed penalty function for nonconvex optimization. Comparing to these methods, the present conformal mapping is straight-forward, and only a pair of well-defined simultaneous linear systems in Eqs. (\ref{eq:3.11}) and (\ref{eq:3.13'}) should be solved. Therefore, the present conformal mapping of deep tunnel is highly efficient and accurate.

\subsection{Shallow tunnel}%
\label{sub:Shallow tunnel}

The most classical conformal mapping is the Verruijt mapping \cite{Verruijt1997traction,Verruijt1997displacement}, which is fully analytic and would not be considered as a comparison with the present numerical one. 

The second conformal mapping is the extensions of the Verruijt mapping proposed by Zeng et al. \cite{Zengguisen2019}, which is semi-analytical with extra items to adapt noncircular cavity boundary:
\begin{equation}
  \label{eq:6.1}
  z = z(\zeta) = -{\rm{i}}a\frac{1+\zeta}{1-\zeta} + {\rm{i}}\sum\limits_{k=1}^{n} b_{k}\left( \zeta^{k} - \zeta^{-k} \right)
\end{equation}
where $ b_{k} $ are real constants to be determined. The equality of $ \overline{z}(\zeta) = -z(\zeta) $ of Eq. (\ref{eq:6.1}) indicates that such an extending mapping by Zeng et al. \cite{Zengguisen2019} is only suitable for symmetrical cavity, and can not handle asymmetrical cavity as shown in Figs. \ref{fig:7} and \ref{fig:8}. Additionally, the real constants $ b_{k} $ can be obtained using nonconvex optimization \cite{Zengguisen2019}. A feasible solution of the extending mapping \cite{Zengguisen2019} is proposed by the authors \cite{lin2022solution} using the Partical Swarm Optimization \cite{kennedy1995particle}, which requires parallel computation infrastructure (\texttt{Nvidia} GPUs for instance) and corresponding parallel computation suite (\texttt{CUDA} or \texttt{OpenACC}). Apparently, the requirement of infrastructure cost and sophisticated coding skills would limit its usage. 

The third conformal mapping for shallow tunnel is recently proposed by the authors \cite{lin2024over-under-excavation}, which is also bidirectional and stepwise by incorporating Charge Simulation Method. The first step of the forward conformal mapping applies the Charge Simulation Method to the cavity boundary in the physical plane alone, so that the cavity in the interval mapping plane is a unit circle, but the ground surface in the interval mapping plane is laterally affected to a curve instead of a straight line. The second step of forward conformal mapping applies the Verruijt conformal mapping, and the exterior boundary of the mapping annulus is no longer circular, which is corresponding to the non-straight ground surface in the interval mapping plane. In theory, the complex potentials can only be expanded into bilateral Laurent series formation, as long as the mapping region is a fully circular annulus. However, with the third conformal mapping \cite{lin2024over-under-excavation}, the mapping region is just a quasi circular annulus with noncircular exterior boundary, instead of a fully circular one. To maintain the usage of complex potentials in Laurent series formation, the noncircular annulus has to be simulated by a circular one (the geometry is changed), and certain errors would exist in the simulating backward conformal mapping, and the ground surface will no longer hold a straight line (see numerical example in Figs. \ref{fig:15}a and b), as well as subsequent complex variable solution. In other words, the third bidirectional conformal mapping \cite{lin2024over-under-excavation} is obscure and indirect with complicated and extra simulation.

In contrast, though the present conformal mapping of shallow tunnel is stepwise by incorporating Charge Simulation Method, the mapping sequence is reorganized. The first step applies the Verruijt mapping to transform the infinite lower half plane containing an arbitrary cavity to a finite unit annulus with noncircular interior boundary in the interval mapping plane. The second step applies the Charge Simulation Method to the unit annulus to further transform the unit annulus in the interval mapping plane to a full unit annulus with both circular boundaries, so that the complex potentials can be spontaneously expanded into Larent series formation without any further simulation. In other words, the present conformal mapping is cleaner and more accurate, comparing to the third existing conformal mapping \cite{lin2024over-under-excavation}. Moreover, the third bidirectional conformal mapping \cite{lin2024over-under-excavation} claims that the Charge Simulaion Method can not deal with compicated cavity shape and sharp corners, which has been proven possible in Figs. \ref{fig:7} and \ref{fig:8} by rounding corners with small arcs and densifying {\emph{collocation points}}.

\section{Complex variable solution using present conformal mapping and verification}
\label{sec:Complex variable solution}

The aim of this paper is to extend the usage of complex variable solution \cite{Muskhelishvili1966} by using present bidirectional conformal mapping. Therefore, in this section, two numerical cases of noncircular and asymmetrical shallow tunnelling in gravitational geomaterial with reasonable far-field displacement is investigated for verification.

\subsection{Mechanical model and brief solution procedure}%
\label{sub:Mechanical model and brief solution procedure}

\begin{figure}[htb]
  \centering
  \begin{tikzpicture}[scale=0.8, transform shape]
    \tikzstyle{every node} = [scale = 0.7]
    \fill [gray!30] (-4,0) rectangle (4,-6);
    \fill [gray!50] (1,-2) arc [start angle = 0, end angle = 270, radius = 1] -- (1,-3) -- (1,-2);
    \draw [dashed, violet, line width = 1pt] (1,-2) arc [start angle = 0, end angle = 270, radius = 1] -- (1,-3) -- (1,-2);
    \draw [line width = 1pt, Emerald] (-4,0) -- (4,0);
    \foreach \x in {1,2,3,...,12} \draw [->] ({-4-\x/12*0.5},{-\x/12*6}) -- (-4,{-\x/12*6});
    \draw (-4,0) -- (-4.5,-6) node [below] {$ k_{x}\gamma y $};
    \foreach \x in {1,2,3,...,12} \draw [->] ({4+\x/12*0.5},{-\x/12*6}) -- (4,{-\x/12*6});
    \draw (4,0) -- (4.5,-6) node [below] {$ k_{x}\gamma y $};
    \foreach \x in {0,1,2,...,16} \draw [->] ({-4+\x/16*8},-6.5) -- ({-4+\x/16*8},-6);
    \draw (-4,-6.5) -- (4,-6.5);
    \node at (0,-6.5) [below] {$ \gamma y $};
    \draw [->, line width = 1.5pt] (3,-4) -- (3,-5) node [below right] {$ \gamma $};
    \foreach \x in {0,1,2,3} \draw [->] ({-0.2+cos(0+22.5*\x)},{-2+sin(0+22.5*\x)}) -- ({0+cos(0+22.5*\x)},{-2+sin(0+22.5*\x)});
    \foreach \x in {1,2,3,4,5,6,7} \draw [->] ({0.2+cos(90+22.5*\x)},{-2+sin(90+22.5*\x)}) -- ({0+cos(90+22.5*\x)},{-2+sin(90+22.5*\x)});
    \foreach \x in {1,2,3} \draw [->] (-0.2+1,{-2-0.25*\x}) -- (1,{-2-0.25*\x});
    \foreach \x in {1,2,3,4,5,6,7} \draw [->] ({0+cos(0+22.5*\x)},{-0.2-2+sin(0+22.5*\x)}) -- ({0+cos(0+22.5*\x)},{-2+sin(0+22.5*\x)});
    \foreach \x in {1,2,3,4} \draw [->] ({0+cos(180+22.5*\x)},{0.2-2+sin(180+22.5*\x)}) -- ({0+cos(180+22.5*\x)},{-2+sin(180+22.5*\x)});
    \foreach \x in {1,2,3} \draw [->] ({0+0.25*\x},-3+0.2) -- ({0+0.25*\x},-3);
    \node at (1,-2) [right] {$ X_{0}(S) $};
    \node at (0,-2.8) [above] {$ Y_{0}(S) $};
    \node at (-4,-6) [above right] {$ \overline{\bm{\varOmega}}_{0} = \overline{\bm{D}} \cup {\bm{C}} $};
    \draw [->] (0,0) -- (0.5,0) node [above right] {$ x $};
    \draw [->] (0,0) -- (0,0.5) node [above right] {$ y $};
    \node at (0,0) [above right] {$ O $};
    \node at (-2,0) [above, Emerald] {Ground surface $ \partial{\bm{D}}_{1} $};
    \node at (0,-3) [below, violet] {Cavity boundary $ \partial{\bm{D}}_{2} $};
    \node at (0,-2) {$ {\bm{C}} $};
    \node at (0,1) {(a) Geomaterial and initial stress field};

    \node at (5,-3) {\LARGE{+}};
    
    \fill [gray!30] (6,0) rectangle (14,-6);
    \fill [white] (11,-2) arc [start angle = 0, end angle = 270, radius = 1] -- (11,-3) -- (11,-2);
    \draw [violet, line width = 1pt] (11,-2) arc [start angle = 0, end angle = 270, radius = 1] -- (11,-3) -- (11,-2);
    \draw [Emerald, line width = 1pt] (7,0) -- (13,0);
    \fill [pattern = north east lines] (6,0) rectangle (7,0.2);
    \fill [pattern = north east lines] (13,0) rectangle (14,0.2);
    \draw [line width = 1pt] (6,0) -- (7,0);
    \draw [line width = 1pt] (13,0) -- (14,0);
    \foreach \x in {0,1,2,3} \draw [<-] ({10-0.2+cos(0+22.5*\x)},{-2+sin(0+22.5*\x)}) -- ({10+cos(0+22.5*\x)},{-2+sin(0+22.5*\x)});
    \foreach \x in {1,2,3,4,5,6,7} \draw [<-] ({10+0.2+cos(90+22.5*\x)},{-2+sin(90+22.5*\x)}) -- ({10+cos(90+22.5*\x)},{-2+sin(90+22.5*\x)});
    \foreach \x in {1,2,3} \draw [<-] (10-0.2+1,{-2-0.25*\x}) -- (11,{-2-0.25*\x});
    \foreach \x in {1,2,3,4,5,6,7} \draw [<-] ({10+cos(0+22.5*\x)},{-0.2-2+sin(0+22.5*\x)}) -- ({10+cos(0+22.5*\x)},{-2+sin(0+22.5*\x)});
    \foreach \x in {1,2,3,4} \draw [<-] ({10+cos(180+22.5*\x)},{0.2-2+sin(180+22.5*\x)}) -- ({10+cos(180+22.5*\x)},{-2+sin(180+22.5*\x)});
    \foreach \x in {1,2,3} \draw [<-] ({10+0.25*\x},-3+0.2) -- ({10+0.25*\x},-3);
    \node at (10.8,-2) [left] {$ X(S) $};
    \node at (10,-2.8) [above] {$ Y(S) $};
    \draw [->] (10,0) -- (10.5,0) node [above right] {$ x $};
    \draw [->] (10,0) -- (10,0.5) node [above right] {$ y $};
    \node at (10,0) [above right] {$ O $};
    \node at (10,0) [below, Emerald] {Free ground surface $ \partial{\bm{D}}_{11} $};
    \node at (6,0) [below right] {$ \partial{\bm{D}}_{12} $};
    \node at (14,0) [below left] {$ \partial{\bm{D}}_{12} $};
    \node at (10,-0.5) [below] {$ \partial{\bm{D}}_{1} = \partial{\bm{D}}_{11} \cup \partial{\bm{D}}_{12} $};
    \fill [red] (7,0) circle [radius = 0.05];
    \fill [red] (13,0) circle [radius = 0.05];
    \node at (7,0) [below, red] {$ T_{1} $};
    \node at (13,0) [below, red] {$ T_{2} $};
    \node at (10,-3) [below, violet] {Cavity boundary $ \partial{\bm{D}}_{2} $};
    \node at (6,-6) [above right] {$ \overline{\bm{D}} $};
    \node at (10,1) {(b) Mixed boundaries of reansonable far-field displacement and cavity excavation};
  \end{tikzpicture}
  \caption{Mechanical model of shallow tunnelling in gravitational geomaterial with reasonable far-field displacement}
  \label{fig:9}
\end{figure}
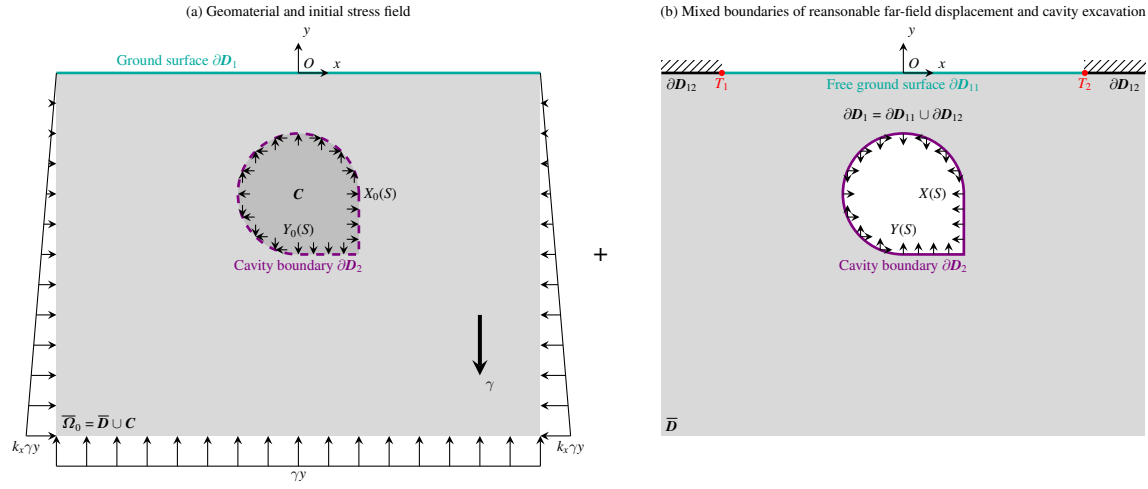

Fig. \ref{fig:9} shows the elastic mechanical model of the numerical case for verification, which is decomposed into one sub-model of prior to excavation in Fig. \ref{fig:9}a and the other one of excavation in Fig. \ref{fig:9}b.

In Fig. \ref{fig:9}a, geomaterial of elastic modulus $ E $ and Poisson's ratio $ \nu $ ocuppies the lower half plane in a rectangular complex coordinate system $ z = x+{\rm{i}}y $. The ground surface and cavity boundary are denoted by $ \partial{\bm{D}}_{1} $ and $ \partial{\bm{D}}_{2} $, respectively, just as in previous sections. The geomaterial within cavity boundary $ \partial{\bm{D}}_{2} $ is denoted by $ {\bm{C}} $; the geomaterial between $ \partial{\bm{D}}_{1} $ and $ \partial{\bm{D}}_{2} $ is denoted by $ {\bm{D}} $, and $ \overline{\bm{D}} = \partial{\bm{D}}_{1} \cup {\bm{D}} \cup \partial{\bm{D}}_{2} $, just as in previous sections. The whole geomaterial $ \overline{\bm{\varOmega}}_{0} = \overline{\bm{D}} \cup {\bm{C}} $. The ground surface $ \partial{\bm{D}}_{1} $ is fully free, and the geomaterial in the lower half plane $ \overline{\bm{\varOmega}}_{0} $ is subjected to the initial stress field as
\begin{equation}
  \label{eq:7.1}
  \left\{
    \begin{aligned}
      & \sigma_{x}^{0}(z) = k_{x}\gamma y \\
      & \sigma_{y}^{0}(z) = \gamma y \\
      & \tau_{xy}^{0}(z) = 0 \\
    \end{aligned}
  \right.
  , \quad z = x+{\rm{i}}y \in \overline{\bm{\varOmega}}_{0}
\end{equation}
where $ \sigma_{x}^{0} $, $ \sigma_{y}^{0} $, and $ \tau_{xy}^{0} $ denote horizontal, vertical, and shear components of initial stress field, respectively; $ \gamma $ denotes unit weight of geomaterial; $ k_{x} $ denotes lateral stress coefficient. 

The initial stress field causes surface traction along cavity boundary $ \partial{\bm{D}}_{2} $ as
\begin{equation}
  \label{eq:7.2}
  \left\{
    \begin{aligned}
      & X_{0}(S) = \sigma_{x}^{0}(S)\cdot\cos\langle \vec{n},\vec{x} \rangle + \tau_{xy}^{0}(S)\cdot\cos\langle \vec{n},\vec{y} \rangle \\
      & Y_{0}(S) = \tau_{xy}^{0}(S)\cdot\cos\langle \vec{n},\vec{x} \rangle + \sigma_{y}^{0}(S)\cdot\cos\langle \vec{n},\vec{y} \rangle
    \end{aligned}
  \right.
  , \quad S \in \partial{\bm{D}}_{2}
\end{equation}
where $ X_{0}(S) $ and $ Y_{0}(S) $ denote horizontal and vertical surface traction of arbitrary point along cavity boundary $ \partial{\bm{D}}_{2} $, respectively; $ \vec{n} $ denotes outward normal of cavity boundary $ \partial{\bm{D}}_{2} $ from region $ {\bm{D}} $.

In Fig. \ref{fig:9}b, the ground surface is divided into two segments. The far-field ground surface segment $ \partial{\bm{D}}_{12} $ is fixed, while the rest finite one $ \partial{\bm{D}}_{11} $ is left free, and the whole ground surface has $ \partial{\bm{D}}_{1} = \partial{\bm{D}}_{11} \cup \partial{\bm{D}}_{12} $. The joint points of segments $ \partial{\bm{D}}_{11} $ and $ \partial{\bm{D}}_{12} $ are denoted by $ T_{1} $ and $ T_{2} $, respectively. \added{In practical tunnel engineering, excavation would generally affect geomaterial near the tunnel, while the far-field ground should remain intact with no deformation. If the ground surface is set to be fully free, displacement singularity at infinity would be laterally raised due to the nonzero resultant caused by excavation, as reported in our previous study \cite{Self2020JEM}. Therefore, more reasonable mechanical models \cite{LIN2024appl_math_model,LIN2024comput_geotech_1,lin2024over-under-excavation} should be applied to cancel such displacement singularity at infinity. In this section, we would use such a reasonable mechanical model to illustrate the possible application of the present bidirectional conformal mapping.}

The stress and displacement within the geomaterial region $ \overline{\bm{D}} $ using complex potentials as can be expressed as \cite{Muskhelishvili1966}
\begin{subequations}
  \label{eq:7.3}
  \begin{equation}
    \label{eq:7.3a}
    \left\{
      \begin{aligned}
        & \sigma_{y}(z)+\sigma_{x}(z) = 2\left[ \varphi^{\prime}(z) + \overline{\varphi^{\prime}(z)} \right] \\
        & \sigma_{y}(z)-\sigma_{x}(z) + 2{\rm{i}}\tau_{xy}(z) = 2\left[ \overline{z}\varphi^{\prime\prime}(z) + \psi^{\prime}(z) \right] \\
      \end{aligned}
    \right.
    , \quad z \in \overline{\bm{D}}
  \end{equation}
  \begin{equation}
    \label{eq:7.3b}
    2G[u_{x}(z)+{\rm{i}}u_{y}(z)] = \kappa \varphi(z) - z \overline{\varphi^{\prime}(z)} - \overline{\psi(z)}, \quad z \in \overline{\bm{D}}
  \end{equation}
\end{subequations}
where $ \sigma_{x} $, $ \sigma_{y} $, and $ \tau_{xy} $ denote horizontal, vertical, and shear stress components due to excavation, respectively; $ u_{x} $ and $ u_{y} $ denote horizontal and vertical displacement components due to excavation, respectively; $ \varphi(z) $ and $ \psi(z) $ denote complex potentials to be uiquely determined by the well-defined mixed boundary conditions in Eq. (\ref{eq:7.5}) below \cite{LIN2024appl_math_model,LIN2024comput_geotech_1,lin2024over-under-excavation}; $ \bullet^{\prime} $ and $ \bullet^{\prime\prime} $ denote first and second derivatives of function $ \bullet $, respectively; $ G = \frac{E}{2(1+\nu)} $ denotes shear modulus; $ \kappa = 3-4\nu $ is the Kolosov parameter for plane strain condition. The final stress within geomaterial can be computed as the sum of Eq. (\ref{eq:7.1}) and (\ref{eq:7.3a}):
\begin{equation}
  \label{eq:7.3a'}
  \tag{7.3a'}
  \left\{
    \begin{aligned}
      & \sigma_{x}^{\ast}(z) = \sigma_{x}^{0}(z) + \sigma_{x}(z) \\
      & \sigma_{y}^{\ast}(z) = \sigma_{y}^{0}(z) + \sigma_{y}(z) \\
      & \tau_{xy}^{\ast}(z) = \tau_{xy}^{0}(z) + \tau_{xy}(z) \\
    \end{aligned}
  \right.
  , \quad z \in \overline{\bm{D}}
\end{equation}
where $ \sigma_{x}^{\ast} $, $ \sigma_{y}^{\ast} $, and $ \tau_{xy}^{\ast} $ denote horizontal, vertical, and shear components of final stress field.

In Fig. \ref{fig:9}b, the inverse surface traction along cavity boundary $ \partial{\bm{D}}_{2} $ is applied to cancel the surface traction in Eq. (\ref{eq:7.2}) to mechanically simulate excavation as
\begin{equation}
  \label{eq:7.4}
  \left\{
    \begin{aligned}
      & X(S) = -X_{0}(S) \\
      & Y(S) = -Y_{0}(S) \\
    \end{aligned}
  \right.
  , \quad S \in \partial{\bm{D}}_{2}
\end{equation}
where $ X $ and $ Y $ denote horizontal and vertical component of the inverse surface traction along arbitrary boundary in the remaing geomaterial region $ \overline{\bm{D}} $, respectively. 

Eq. (\ref{eq:7.4}) can be rewritten into an alternative integral formation as
\begin{subequations}
  \label{eq:7.5}
  \begin{equation}
    \label{eq:7.5a}
    \int_{\partial{\bm{D}}_{2}} \left[ X(S) + {\rm{i}}Y(S) \right]{\rm{d}}S = - \int_{\partial{\bm{D}}_{2}} \left[ X_{0}(S) + {\rm{i}} Y_{0}(S) \right] {\rm{d}}S
  \end{equation}	
  The fixed far-field ground surface should maintain zero displacement as
  \begin{equation}
    \label{eq:7.5b}
    u_{x}(T) + {\rm{i}}u_{y}(T) = 0, \quad T \in \partial{\bm{D}}_{12}
  \end{equation}
  where $ u_{x} $ and $ u_{y} $ denote horizontal and vertical displacement components of geomaterial in physical plane $ z = x+{\rm{i}}y $. The rest finite and free ground surface above cavity is denoted by $ \partial{\bm{D}}_{11} $, and the traction can be expressed as
  \begin{equation}
    \label{eq:7.5c}
    X(T)+{\rm{i}}Y(T) = 0, \quad T \in \partial{\bm{D}}_{11}
  \end{equation}
\end{subequations}

The inverse surface traction in Eq. (\ref{eq:7.5a}) leads to an upward resultant $ R_{y} $ as \cite{Strack2002phdthesis,fang2023analytical,LIN2024appl_math_model}	
\begin{equation*}
  R_{y} = \gamma \iint_{\bm{C}} {\rm{d}}x{\rm{d}}y
\end{equation*}
The nonzero resultant $ R_{y} $ should be equilibriated by the fixed far-field ground surface $ \partial{\bm{D}}_{12} $ in Fig. \ref{fig:9}b \cite{LIN2024appl_math_model,LIN2024comput_geotech_1,lin2024over-under-excavation} as
\begin{equation}
  \label{eq:7.6}
  \int_{\partial{\bm{D}}_{12}} \left[ X(T) + {\rm{i}}Y(T) \right]{\rm{d}}T = -{\rm{i}}R_{y}
\end{equation}
\added{Eq. (\ref{eq:7.6}) shows that the fixed far-field ground surface $ \partial{\bm{D}}_{12} $ would generate necessary constraining force to equilibrate the nonzero resultant caused by excavation of the geomaterial within tunnel periphery $ \partial{\bm{D}}_{2} $, and the mechanical model in Fig. \ref{fig:9}b is an equilibrium one.} Furthermore, the displacement in Eq. (\ref{eq:7.3b}) should be single-valued.

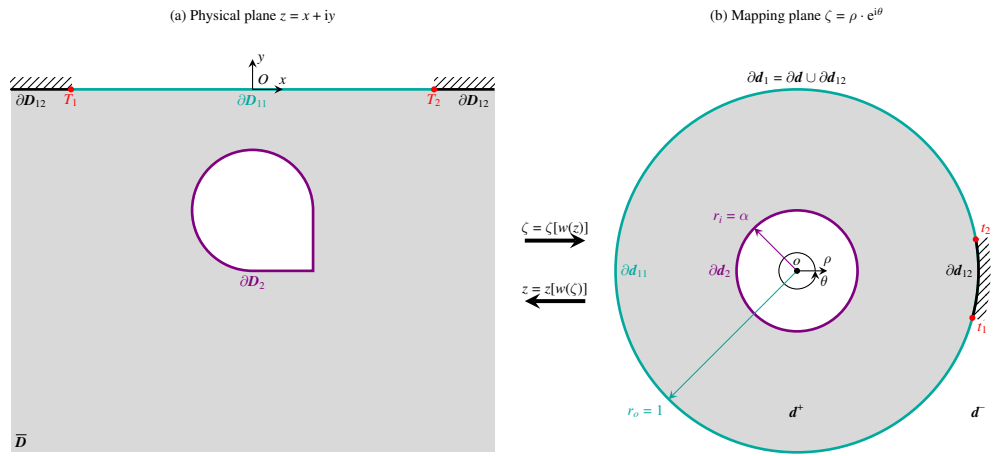
\begin{figure}[htpb]
  \centering
  \begin{tikzpicture}[scale=0.8, transform shape]
    \tikzstyle{every node} = [scale = 0.7]
    \fill [gray!30] (-4,0) rectangle (4,-6);
    \fill [white] (1,-2) arc [start angle = 0, end angle = 270, radius = 1] -- (1,-3) -- (1,-2);
    \draw [violet, line width = 1pt]  (1,-2) arc [start angle = 0, end angle = 270, radius = 1] -- (1,-3) -- (1,-2);
    \fill [pattern = north east lines] (-4,0) rectangle (-3,0.2);
    \fill [pattern = north east lines] (4,0) rectangle (3,0.2);
    \draw [line width = 1pt] (-4,0) -- (-3,0);
    \draw [line width = 1pt] (3,0) -- (4,0);
    \draw [line width = 1pt, Emerald] (-3,0) -- (3,0);
    \fill [red] (-3,0) circle [radius = 0.05];
    \fill [red] (3,0) circle [radius = 0.05];
    \node at (-4,0) [below right] {$ \partial{\bm{D}}_{12} $};
    \node at (4,0) [below left] {$ \partial{\bm{D}}_{12} $};
    \node at (-3,0) [below,red] {$ T_{1} $};
    \node at (3,0) [below,red] {$ T_{2} $};
    \node at (0,0) [below,Emerald] {$ \partial{\bm{D}}_{11} $};
    \node at (0,-3) [below,violet] {$ \partial{\bm{D}}_{2} $};
    \node at (-4,-6) [above right] {$ \overline{\bm{D}} $};
    \draw [->] (0,0) -- (0.5,0) node [above] {$ x $};
    \draw [->] (0,0) -- (0,0.5) node [right] {$ y $};
    \node at (0,0) [above right] {$ O $};
    \node at (0,1) [above] {(a) Physical plane $ z = x + {\rm{i}}y $};

    \fill [gray!30] (9,-3) circle [radius = 3];
    \fill [white] (9,-3) circle [radius = 1];
    \draw [line width = 1pt, Emerald] (9,-3) circle [radius = 3];
    \fill [pattern = north east lines] ({9+3.2*cos(-15)},{-3+3.2*sin(-15)}) arc [start angle = -15, end angle = 10, radius = 3.2] -- ({9+3*cos(10)},{-3+3*sin(10)}) arc [start angle = 10, end angle = -15, radius = 3] -- ({9+3.2*cos(-15)},{-3+3.2*sin(-15)});
    \draw [line width = 1pt] ({9+3*cos(-15)},{-3+3*sin(-15)}) arc [start angle = -15, end angle = 10, radius = 3];
    \draw [line width = 1pt, violet] (9,-3) circle [radius = 1];
    \fill [red] ({9+3*cos(10)},{-3+3*sin(10)}) circle [radius = 0.05];
    \fill [red] ({9+3*cos(-15)},{-3+3*sin(-15)}) circle [radius = 0.05];
    \node at (6,-3) [right, Emerald] {$ \partial{\bm{d}}_{11} $};
    \node at (12,-3) [left] {$ \partial{\bm{d}}_{12} $};
    \node at (8,-3) [left,violet] {$ \partial{\bm{d}}_{2} $};
    \node at (9,-5.5) [above] {$ {\bm{d}}^{+} $};
    \node at (12,-5.5) [above] {$ {\bm{d}}^{-} $};
    \node at (9,0) [above] {$ \partial{\bm{d}}_{1} = \partial{\bm{d}} \cup \partial{\bm{d}}_{12} $};
    \node at ({9+3*cos(-15)},{-3+3*sin(-15)}) [below right,red] {$ t_{1} $};
    \node at ({9+3*cos(10)},{-3+3*sin(10)}) [above right,red] {$ t_{2} $};
    \draw [->] (9,-3) -- (9.5,-3) node [above] {$ \rho $};
    \draw [->] (9.3,-3) arc [start angle = 0, end angle = 360, radius = 0.3];
    \node at (9.3,-3) [below right] {$ \theta $};
    \node at (9,-3) [above] {$ o $};
    \draw [->,violet] (9,-3) -- ({9+cos(135)},{-3+sin(135)}) node [above left] {$ r_{i} = \alpha $};
    \draw [->,Emerald] (9,-3) -- ({9+3*cos(225)},{-3+3*sin(225)}) node [below left] {$ r_{o} = 1 $};
    \fill (9,-3) circle [radius = 0.05];
    \node at (9,1) [above] {(b) Mapping plane $ \zeta = \rho \cdot {\rm{e}}^{{\rm{i}}\theta} $};

    \draw [->, line width = 1.5pt] (4.5,-2.5) -- (5.5,-2.5);
    \draw [->, line width = 1.5pt] (5.5,-3.5) -- (4.5,-3.5);
    \node at (5,-2.5) [above] {$ \zeta = \zeta[w(z)] $};
    \node at (5,-3.5) [above] {$ z = z[w(\zeta)] $};
  \end{tikzpicture}
  \caption{Schematic diagram of bidirectional conformal mapping of mechanical model for verification}
  \label{fig:10}
\end{figure}

The geomaterial in Fig. \ref{fig:9}b can be mapped onto a circular unit annulus $ \overline{\bm{d}} $ as shown in Fig. \ref{fig:10}. The mapping pattern is the same to that in Fig. \ref{fig:2}b with more specific details that the finite ground surface $ \partial{\bm{D}}_{11} $, the fixed far-field ground surface $ \partial{\bm{D}}_{12} $, and the joint points $ T_{1} $ and $ T_{2} $ in physical plane ($ \overline{\bm{D}} : z = x+{\rm{i}}y $) are bidirectionally mapped onto arc of unit radius $ \partial{\bm{d}}_{11} $, arc of unit radius $ \partial{\bm{d}}_{12} $, and joint points $ t_{1} $ and $ t_{2} $ in mapping plane ($ \overline{\bm{d}} : \zeta = \rho\cdot{\rm{e}}^{{\rm{i}}\theta} $), respectively. Additionally, arbitrary points $ T $ along ground surface $ \partial{\bm{D}}_{1} $ and $ S $ along cavity boundary $ \partial{\bm{D}}_{2} $ in the physical plane are bidirectionally mapped onto corresponding points $ t $ along arc $ \partial{\bm{d}}_{1} $ and $ \partial{\bm{d}}_{2} $, respectively. Apparently, we have $ \partial{\bm{d}}_{1} = \partial{\bm{d}}_{11} \cup \partial{\bm{d}}_{12} $. The mapping annulus without exterior boundary is denoted by $ {\bm{d}}^{+} = {\bm{d}} \cup \partial{\bm{d}}_{2}  $, and the region outside of the exterior boundary of the annulus is denoted by $ {\bm{d}}^{-} $.

With the bidirectional conformal mapping in Fig. \ref{fig:10}, the stress and displacement in physical plane in Eq. (\ref{eq:7.3}) can be bidirectionally mapped onto their curvilinear ones in the mapping plane as
\begin{subequations}
  \label{eq:7.7}
  \begin{equation}
    \label{eq:7.7a}
    \sigma_{\theta}(\zeta) + \sigma_{\rho}(\zeta) = 2\left[ \varPhi(\zeta) + \overline{\varPhi(\zeta)} \right], \quad \zeta \in \overline{\bm{d}}
  \end{equation}
  \begin{equation}
    \label{eq:7.7b}
    \sigma_{\rho}(\zeta) + {\rm{i}}\tau_{\rho\theta}(\zeta) = \left[ \varPhi(\zeta) + \overline{\varPhi(\zeta)} \right] - \frac{\overline{\zeta}}{\zeta} \left[ \frac{z(\zeta)}{z^{\prime}(\zeta)}\overline{\varPhi^{\prime}(\zeta)} + \frac{\overline{z^{\prime}(\zeta)}}{z^{\prime}(\zeta)}\overline{\varPsi(\zeta)} \right], \quad \zeta \in \overline{\bm{d}}
  \end{equation}
  \begin{equation}
    \label{eq:7.7c}
    g(\zeta) = 2G\left[ u_{x}(\zeta)+{\rm{i}}u_{y}(\zeta) \right] = \kappa \varphi(\zeta) - z(\zeta) \overline{\varPhi(\zeta)} - \overline{\psi(\zeta)}, \quad \zeta \in \overline{\bm{d}}
  \end{equation}
  with
  \begin{equation}
    \label{eq:7.7d}
    \left\{
      \begin{aligned}
        & \varPhi(\zeta) = \frac{\varphi^{\prime}(\zeta)}{z^{\prime}(\zeta)} \\
        & \varPsi(\zeta) = \frac{\psi^{\prime}(\zeta)}{z^{\prime}(\zeta)} \\
      \end{aligned}
    \right.
  \end{equation}
\end{subequations}
where $ \sigma_{\theta}(\zeta) $, $ \sigma_{\rho}(\zeta) $, and $ \tau_{\rho\theta}(\zeta) $ denote hoop, radia, and tangential components of the curvlinear stress components mapped onto the mapping plane, respectively; $ g(\zeta) $ denotes generalized rectangular displacement mapped onto the mapping plane. It should be addressed that $ z(\zeta) = z[w(\zeta)] $ with composite mappings and that $ z^{\prime}(\zeta) = \frac{{\rm{d}}z[w(\zeta)]}{{\rm{d}}\zeta} $ should take chain law of composite function. 

With the bidirectional conformal mapping in Fig. \ref{fig:10}, the three mixed boundary conditions in Eqs. (\ref{eq:7.5a}), (\ref{eq:7.5b}), and (\ref{eq:7.5c}) can be also bidirectionally mapping onto their corresponding mapping formations, respectively, as
\cite{LIN2024appl_math_model,LIN2024comput_geotech_1,lin2024over-under-excavation}
\begin{subequations}
  \label{eq:7.8}
  \begin{equation}
    \label{eq:7.8a}
    \int_{\partial{\bm{d}}_{2}} z^{\prime}(s) \left[ \sigma_{\rho}(s) + {\rm{i}}\tau_{\rho\theta}(s) \right] {\rm{d}}s = -\gamma \int_{\partial{\bm{d}}_{2}} y(s)\left[ k_{x}\frac{{\rm{i}{\rm{d}}}y(s)}{{\rm{d}}\sigma} + \frac{{\rm{d}}x(s)}{{\rm{d}}\sigma} \right] {\rm{d}}\sigma
  \end{equation}
  \begin{equation}
    \label{eq:7.8b}
    g(t) = 2G\left[ u_{x}(t) + u_{y}(t) \right] = 0, \quad t \in \partial{\bm{d}}_{12}
  \end{equation}
  \begin{equation}
    \label{eq:7.8c}
    \sigma_{\rho}(t) + {\rm{i}}\tau_{\rho\theta}(t) = 0, \quad t \in \partial{\bm{d}}_{11}
  \end{equation}
  where
  \begin{equation}
    \label{eq:7.9d}
    \left\{
      \begin{aligned}
        & x(s) = \frac{1}{2}\left\{ \overline{z[w(s)]} + z[w(s)] \right\} \\
        & y(s) = \frac{\rm{i}}{2}\left\{ \overline{z[w(s)]} - z[w(s)] \right\} \\
      \end{aligned}
    \right.
  \end{equation}
\end{subequations}

Eqs. (\ref{eq:7.8b}) and (\ref{eq:7.8c}) can be transformed to a homogenerous Riemann-Hilbert problem as \cite{LIN2024appl_math_model,LIN2024comput_geotech_1,lin2024over-under-excavation}
\begin{subequations}
  \label{eq:7.9}
  \begin{equation}
    \label{eq:7.9a}
    \left. z^{\prime}(\zeta)\left[ \sigma_{\rho}(\zeta) + {\rm{i}}\tau_{\rho\theta}(\zeta) \right] \right|_{\rho \rightarrow 1} = \varphi^{\prime+}(\sigma) - \varphi^{\prime-}(\sigma) = 0, \quad \sigma \in \partial{\bm{d}}_{11}
  \end{equation}
  \begin{equation}
    \label{eq:7.9b}
    g^{\prime}(\zeta)|_{\rho \rightarrow 1} = \kappa \varphi^{\prime+}(\sigma) + \varphi^{\prime-}(\sigma) = 0, \quad \sigma \in \partial{\bm{d}}_{12}
  \end{equation}
  where $ \sigma = {\rm{e}}^{{\rm{i}}\theta} $; $ \varphi^{\prime+}(\sigma) $ and $ \varphi^{\prime-}(\sigma) $ denote the boundary values of $ \varphi^{\prime}(\zeta) $ approaching boundary $ \partial{\bm{d}}_{1} $ from the inside region $ {\bm{d}}^{+} $ and outside region $ {\bm{d}}^{-} $ in Fig. \ref{fig:10}b, respectively. Eq. (\ref{eq:7.9}) implicitly indicates that $ \varphi^{\prime}(\zeta) $ is analytic within the whole mapping region $ {\bm{d}}^{+} \cup \partial{\bm{d}}_{1} \cup {\bm{d}}^{-} $.
\end{subequations}

The general solution of Eq. (\ref{eq:7.9}) can be expressed as
\begin{equation}
  \label{eq:7.10}
  \varphi^{\prime}(\zeta) = X(\zeta) \sum\limits_{n=-\infty}^{\infty} f_{n}\zeta^{n}
\end{equation}
where $ f_{n} $ denote complex coefficients to be determined, and $ X(\zeta) $ can be expanded as
\begin{equation}
  \label{eq:7.11}
  X(\zeta) = 
  \left\{
    \begin{aligned}
      & \sum\limits_{k=0}^{\infty} \alpha_{k}\zeta^{k}, \quad \zeta \in {\bm{d}}^{+} \\
      & \sum\limits_{k=0}^{\infty} \beta_{k}\zeta^{k}, \quad \zeta \in {\bm{d}}^{-} \\
    \end{aligned}
  \right.
\end{equation}
wherein the expressions of $ \alpha_{k} $ and $ \beta_{k} $ can be found in Ref \cite{lin2024over-under-excavation}, and would not be repeated here.

Since the complex potentials $ \varphi^{\prime}(\zeta) $ is analytic within the whole mapping region according to Eq. (\ref{eq:7.10}), it can be expanded as
\begin{subequations}
  \label{eq:7.12}
  \begin{equation}
    \label{eq:7.12a}
    \varphi^{\prime}(\zeta) = \sum\limits_{k=-\infty}^{\infty} A_{k}\zeta^{k}, \quad \zeta \in {\bm{d}}^{+}
  \end{equation}
  \begin{equation}
    \label{eq:7.12b}
    \varphi^{\prime}(\zeta) = \sum\limits_{k=-\infty}^{\infty} B_{k}\zeta^{k}, \quad \zeta \in {\bm{d}}^{-}
  \end{equation}
\end{subequations}
where $ A_{k} $ and $ B_{k} $ are complex constants to be determined, and can be expanded as \cite{LIN2024appl_math_model,LIN2024comput_geotech_1,lin2024over-under-excavation}
\begin{equation}
  \label{eq:7.13}
  \left\{
    \begin{aligned}
      & A_{k} = \sum\limits_{n=-\infty}^{k} \alpha_{k-n}f_{n} \\
      & B_{k} = \sum\limits_{n=k}^{\infty} \beta_{-k+n}f_{n} \\
    \end{aligned}
  \right.
\end{equation}
The other complex potential $ \replaced{\psi}{\varphi}^{\prime}(\zeta) $ can be expressed using $ \varphi^{\prime}(\zeta) $ as \cite{lin2024non}
\begin{equation}
  \label{eq:7.14}
  \replaced{\psi}{\varphi}^{\prime}(\zeta) = \sum\limits_{k=-\infty}^{\infty} \overline{B}_{-k-2}\zeta^{k} - \left[ \frac{\replaced{\overline{z(\zeta)}}{\overline{z}(\zeta^{-1})}}{z^{\prime}(\zeta)}\varphi^{\prime}(\zeta) \right]\added{^{\prime}}
\end{equation}

Substituting Eqs. (\ref{eq:7.12a}), (\ref{eq:7.14}), (\ref{eq:7.7d}), and (\ref{eq:7.7b}) into the left-hand side of Eq. (\ref{eq:7.8a}) yields
\begin{equation}
  \label{eq:7.15}
  \begin{aligned}
    & \sum\limits_{\substack{k=-\infty \\ k \neq 0}}^{\infty} \left( A_{k-1}\frac{\alpha^{k}}{k}\sigma^{k} + B_{-k-1}\frac{\alpha^{k}}{k}\sigma^{-k} \right) + \frac{z(\alpha\sigma)-\replaced{\overline{z(\alpha\sigma)}}{z(\alpha^{-1}\sigma)}}{\overline{z^{\prime}(\alpha\sigma)}}\sum\limits_{k=-\infty}^{\infty} \overline{A}_{k}\alpha^{k}\sigma^{-k} \\
    & + (A_{-1}+B_{-1})\ln\alpha + C_{a} + (A_{-1}-B_{-1}){\rm{Ln}}\sigma = -\gamma \int_{\partial{\bm{d}}_{2}} y(s)\left[ k_{x}\frac{{\rm{i}{\rm{d}}}y(s)}{{\rm{d}}\sigma} + \frac{{\rm{d}}x(s)}{{\rm{d}}\sigma} \right] {\rm{d}}\sigma
  \end{aligned}
\end{equation}
where $ C_{a} $ denotes an integral constant. Using sample point simulation \cite{LIN2024appl_math_model}, the following coefficient equilibriums can be established as \cite{lin2024over-under-excavation}
\begin{subequations}
  \label{eq:7.16}
  \begin{equation}
    \label{eq:7.16a}
    A_{-k-1} = \sum\limits_{n=k+1}^{\infty} \alpha_{-k-1+n}f_{n} = -k\alpha^{k}I_{k} + \alpha^{2k}B_{-k-1} + k\alpha^{2k}\sum\limits_{l=-\infty}^{\infty}\alpha^{l}d_{l}(\alpha)\overline{A}_{l+k}, \quad k \geq 1
  \end{equation}
  \begin{equation}
    \label{eq:7.16b}
    B_{k-1} = \sum\limits_{n=k}^{\infty} \beta_{-k+1+n}f_{n} = -k\alpha^{k}J_{k} + \alpha^{2k}A_{k-1} + k\sum\limits_{l=-\infty}^{\infty} \alpha^{l}d_{l}(\alpha)\overline{A}_{l-k}, k \geq 1
  \end{equation}
  where
  \begin{equation*}
    \sum\limits_{k=-\infty}^{\infty} d_{l}(\rho)\sigma^{k} = \frac{z(\rho\sigma)-\replaced{\overline{z(\rho\sigma)}}{z(\rho^{-1}\sigma)}}{\overline{z^{\prime}(\rho\sigma)}}
  \end{equation*}
  \begin{equation*}
    \sum\limits_{k=1}^{\infty} I_{k}\sigma^{k} + \sum\limits_{k=1}^{\infty} J_{k}\sigma^{k} + K_{0}{\rm{Ln}}\sigma = -\gamma \int_{\partial{\bm{d}}_{2}} y(s)\left[ k_{x}\frac{{\rm{i}{\rm{d}}}y(s)}{{\rm{d}}\sigma} + \frac{{\rm{d}}x(s)}{{\rm{d}}\sigma} \right] {\rm{d}}\sigma
  \end{equation*}
  The resultant equilibrium in Eq. (\ref{eq:7.6}) and displacement single-valuedness would give $ A_{-1} $ and $ B_{1} $ free from conformal mapping as \cite{LIN2024appl_math_model,LIN2024comput_geotech_1,lin2024over-under-excavation} 
  \begin{equation}
    \label{eq:7.16c}
    A_{-1} = \sum\limits_{n=1}^{\infty} \alpha_{-1+n}f_{-n} = \frac{K_{0}}{1+\kappa}
  \end{equation}
  \begin{equation}
    \label{eq:7.16d}
    B_{-1} = \sum\limits_{n=0}^{\infty} \beta_{1+n}f_{n} = \frac{-\kappa K_{0}}{1+\kappa}
  \end{equation}
\end{subequations}
The coefficient $ K_{0} $ should sastisfy
\begin{equation*}
  K_{0} = - {\rm{i}} \frac{R_{y}}{2\pi}
\end{equation*}

Eq. (\ref{eq:7.16}) is a well-defined simultaneous linear system to uniquely determine $ f_{n} $, which can be iteratively solved in the same manner in Refs \cite{LIN2024appl_math_model,LIN2024comput_geotech_1,lin2024over-under-excavation}. Then, the curvilinear stress and displacement components in mapping plane $ \overline{\bm{d}} $ can be uniquely determined with truncation of Eq. (\ref{eq:7.10}) into $ 2N_{0}+1 $ items and Lanczos filtering \cite{LIN2024appl_math_model,LIN2024comput_geotech_1,lin2024over-under-excavation} as
\begin{subequations}
  \label{eq:7.17}
  \begin{equation}
    \label{eq:7.17a}
    \sigma_{\theta}(\rho\sigma)+\sigma_{\rho}(\rho\sigma) = 4\Re\left[ \frac{1}{z^{\prime}(\rho\sigma)} \sum\limits_{k=-N_{0}}^{N_{0}} A_{k}\rho^{k}\sigma^{k} \added{\cdot L_{k}} \right] \deleted{\cdot L_{k}}
  \end{equation}
  \begin{equation}
    \label{eq:7.17b}
    \sigma_{\rho}(\rho\sigma)+{\rm{i}}\tau_{\rho\theta}(\rho\sigma) = \frac{1}{z^{\prime}(\rho\sigma)}\sum\limits_{k=-N_{0}}^{N_{0}} \left[ A_{k}\rho^{k} - B_{k}\rho^{-k-2} + (k+1) \sum\limits_{l=-N_{0}+k+1}^{N_{0}+k+1} d_{l}(\rho)\overline{A}_{l-k-1}\rho^{l-k-2} \right]\sigma^{k} \cdot L_{k}
  \end{equation}
  \begin{equation}
    \label{eq:7.17c}
    \begin{aligned}
      2G[u_{x}(\rho\sigma)+{\rm{i}}u_{y}(\rho\sigma)]
      = & \sum\limits_{\substack{k=-N_{0} \\ k\neq 0}}^{N_{0}} \left( \kappa A_{k-1} \frac{\rho^{k}\sigma^{k}-1^{k}}{k} + B_{k-1}\frac{\rho^{-k}\sigma^{k}-1^{k}}{k} \right) \cdot L_{k} \\
      & - \sum\limits_{k=-N_{0}}^{N_{0}}\sum\limits_{l=-M+k}^{M+k} d_{l}(\rho)\overline{A}_{l-k}\rho^{l-k}\sigma^{k} \cdot L_{k}+(\kappa A_{-1}-B_{-1})\ln\rho
    \end{aligned}
  \end{equation}
\end{subequations}
where
\begin{equation*}
  L_{k} = 
  \left\{
    \begin{aligned}
      & 1, \quad k = 0 \\
      & \sin\left( \frac{k}{N_{0}}\pi \right) / \left( \frac{k}{N_{0}}\pi \right), \quad {\rm{otherwise}}
    \end{aligned}
  \right.
\end{equation*}

The stress and displacement due to excavation in the physical plane $ \overline{\bm{D}} $ can be obtained by transformation as
\begin{subequations}
  \label{eq:7.18}
  \begin{equation}
    \label{eq:7.18a}
    \left\{
      \begin{aligned}
        & \sigma_{y}(z) + \sigma_{x}(z) = \sigma_{\theta}[\zeta(z)] + \sigma_{\rho}[\zeta(z)], z \in {\bm{D}} \\
        & \sigma_{y}(z) - \sigma_{x}(z) + 2{\rm{i}}\tau_{xy}(z) = \left\{ \sigma_{\theta}[\zeta(z)] - \sigma_{\rho}[\zeta(z)] + 2{\rm{i}}\tau_{\rho\theta}[\zeta(z)] \right\} \cdot \frac{\overline{\zeta(z)}}{\zeta(z)} \left. \frac{\overline{z^{\prime}(\zeta)}}{z^{\prime}(\zeta)} \right|_{\zeta\rightarrow\zeta(z)}, \quad z \in \overline{\bm{D}}
      \end{aligned}
    \right.
  \end{equation}
  \begin{equation}
    \label{eq:7.18b}
    u_{x}(z) + {\rm{i}}u_{y}(z) =  u_{x}[\zeta(z)] + {\rm{i}}u_{y}[\zeta(z)], \quad z \in \overline{\bm{D}}
  \end{equation}
\end{subequations}
Eq. (\ref{eq:7.18}) is the replacement of Eq. (\ref{eq:7.3}), and the final stress field can be obtained in Eq. (\ref{eq:7.3a'}).

\subsection{First numerical case and comparisons with finite element method}%
\label{sub:Results and comparisons with finite element method}

\begin{figure}[htb]
  \centering
  \begin{tabular}{c}
    (a) Model geometry and seed distribution \\
    \begin{tikzpicture}[scale=0.8, transform shape]
      \tikzstyle{every node} = [scale = 0.8]
      \fill [gray!30] (-8,0) rectangle (8,-8);
      \draw [line width = 1pt, Emerald] (-3,0) -- (3,0);
      \fill [orange!60] (0,-4.5) -- (1.35,-4.5) arc [start angle = -90, end angle = 0, radius = 0.15] -- (1.5,-3) arc [start angle = 0, end angle = 270, radius = 1.5];
      \draw [red, line width = 1pt] (1.5,-3) arc [start angle = 0, end angle = 270, radius = 1.5];
      \draw [violet, line width = 1pt] (0,-4.5) -- (1.35,-4.5);
      \draw [blue, line width = 1pt] (1.35,-4.5) arc [start angle = -90, end angle = 0, radius = 0.15];
      \draw [violet, line width = 1pt] (1.5,-4.35) -- (1.5,-3);
      \fill (0,-4.5) circle [radius = 0.05];
      \fill (1.35,-4.5) circle [radius = 0.05];
      \fill (1.5,-4.35) circle [radius = 0.05];
      \fill (1.5,-3) circle [radius = 0.05];
      \fill [pattern = north east lines] (-3,0) -- (-8,0) -- (-8,-8) -- (8,-8) -- (8,0) -- (3,0) -- (3,0.2) -- (8.2,0.2) -- (8.2,-8.2) -- (-8.2,-8.2) -- (-8.2,0.2) -- (-3,0.2) -- (-3,0);
      \draw [line width = 1pt] (-3,0) -- (-8,0) -- (-8,-8) -- (8,-8) -- (8,0) -- (3,0);
      \draw [dashed] (-3,0) -- (-3,-8);
      \draw [dashed] (-5,0) -- (-5,-8);
      \draw [dashed] (3,0) -- (3,-8);
      \draw [dashed] (5,0) -- (5,-8);
      \draw [dashed] (-5,-6) -- (5,-6);
      \node at (-8,-4) [below, rotate = 90, align = center] {100m \\ \textcolor{cyan}{200 seeds}};
      \node at (8,-4) [above, rotate = 90, align = center] {\textcolor{cyan}{200 seeds} \\ 100m};
      \node at (-6.5,0) [below] [align = center] {80m \\ \textcolor{cyan}{160 seeds}};
      \node at (-6.5,-8) [above] [align = center] {\textcolor{cyan}{160 seeds} \\ 80m};
      \node at (6.5,0) [below] [align = center] {80m \\ \textcolor{cyan}{160 seeds}};
      \node at (6.5,-8) [above] [align = center] {\textcolor{cyan}{160 seeds} \\ 80m};
      \node at (-5,-3) [above, rotate = 90, align = center] {\textcolor{cyan}{40 seeds} \\ 20m};
      \node at (5,-3) [below, rotate = 90, align = center] {20m \\ \textcolor{cyan}{40 seeds}};
      \node at (-5,-7) [above, rotate = 90, align = center] {\textcolor{cyan}{160 seeds} \\ 80m};
      \node at (-3,-7) [below, rotate = 90, align = center] {80m \\ \textcolor{cyan}{160 seeds}};
      \node at (5,-7) [below, rotate = 90, align = center] {80m \\ \textcolor{cyan}{160 seeds}};
      \node at (3,-7) [above, rotate = 90, align = center] {\textcolor{cyan}{160 seeds} \\ 80m};
      \node at (-4,-6) [above, align = center] {\textcolor{cyan}{40 seeds} \\ 10m};
      \node at (-4,-8) [above, align = center] {\textcolor{cyan}{40 seeds} \\ 10m};
      \node at (4,-6) [above, align = center] {\textcolor{cyan}{40 seeds} \\ 10m};
      \node at (4,-8) [above, align = center] {\textcolor{cyan}{40 seeds} \\ 10m};
      \node at (-4,0) [below, align = center] {10m \\ \textcolor{cyan}{120 seeds}};
      \node at (4,0) [below, align = center] {10m \\ \textcolor{cyan}{120 seeds}};
      \node at (-3,-3) [above, rotate = 90, align = center] {\textcolor{cyan}{240 seeds} \\ 20m};
      \node at (3,-3) [below, rotate = 90, align = center] {20m \\ \textcolor{cyan}{240 seeds}};
      \node at (0,-8) [above, align = center] {\textcolor{cyan}{80 seeds} \\ 20m};
      \node at (0,-6) [below, align = center] {20m \\ \textcolor{cyan}{80 seeds}};
      \node at (0,0) [above, align = center, Emerald] {20m $ \quad $ 240 seeds};
      \fill [red] (-3,0) circle [radius = 0.05];
      \fill [red] (3,0) circle [radius = 0.05];
      \node at (-3,0) [above right, red] {$ T_{1} $};
      \node at (3,0) [above left, red] {$ T_{2} $};
      \draw [<->] (0,0) -- (0,-1.5);
      \node at (0,-0.75) [above, rotate = 90] {5m};
      \draw [->,red] (0,-3) -- ({1.5*cos(45)},{-3+1.5*sin(45)}) node [above right] {5m};
      \fill [red] (0,-3) circle [radius = 0.05];
      \node at (0.7,-4.5) [above, violet] {4.5m};
      \node at (1.5,-3.7) [above, rotate = 90, violet] {4.5m};
      \node at (0.,-4.5) [below] {$ S_{1} $};
      \node at (1.5,-3) [right] {$ S_{2} $};
      \node at (1.5,-4.35) [right] {$ S_{3} $};
      \node at (1.35,-4.5) [below] {$ S_{4} $};
      \node at ({1.5*cos(135)},{-3+1.5*sin(135)}) [above, rotate = 45, align = center, red] {$ \overset{\LARGE{\frown}}{S_{1}S_{2}} $ \\ 270 seeds};
      \node at (0.7,-4.5) [below, align = center, violet] {$ \overline{S_{4}S_{1}} $ \\ 100 seeds};
      \node at (1.5,-3.7) [below, align = center, rotate = 90, violet] {$ \overline{S_{2}S_{3}} $ \\ 100 seeds};
      \node at ({2.4*cos(-45)},{-3+2.4*sin(-45)}) [below, rotate = 45, align = center, blue] {$ \overset{\LARGE{\frown}}{S_{3}S_{4}} $ \\ 45 seeds};
    \end{tikzpicture} \\
    (b) Meshing near cavity \\
    \includegraphics[width = 0.8\textwidth]{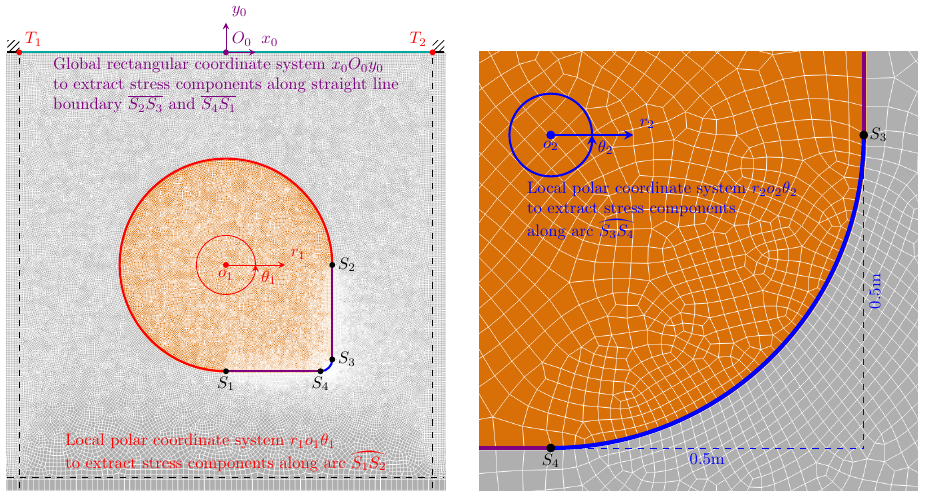}
  \end{tabular}
  \caption{Mechanical model for finite element solution and meshing near cavity}
  \label{fig:11}
\end{figure}

In the first numerical case, the cavity geometry takes the same \added{under-break} one in Fig. \ref{fig:7}a, and the mechanical parameters in Table \ref{tab:1} are used. To verify the feasible usage of present conformal mapping in complex variable method introduced in Section \ref{sub:Mechanical model and brief solution procedure}, a corresponding finite element solution is put forward, and the mechanical model is shown in Fig. \ref{fig:11}. The finite element solution is constructed and computed using finite element software \texttt{ABAQUS 2020}. \added{In Table \ref{tab:1}, the coordinates of the joint points $ T_{1} $ and $ T_{2} $ are chosen on purpose to illustrate good comparisons of the very acute changes of both stress and displacement along ground surface between the complex variable and finite element results, as shown in Eqs. (\ref{eq:7.8b}) and (\ref{eq:7.8c}).}

In Fig. \ref{fig:11}a, the model geometry and seed distribution is illustrated. The cavity boundary is divided into four segments by four points: the first segment is a $ 270^{\circ} $ arc $ \overset{\LARGE{\frown}}{S_{1}S_{2}} $ with radius of $ {\rm{5m}} $; the second segment is a straight line $ \overline{S_{2}S_{3}} $ with length of $ {\rm{4.5m}} $, the third segment is a small $ 90^{\circ} $ arc $ \overset{\LARGE{\frown}}{S_{3}S_{4}} $ with radius $ {\rm{0.5m}} $, which is the rounded corner to be identical to cavity geometry in Fig. \ref{fig:7}a; the last segment is another straight line $ \overline{S_{4}S_{1}} $ with length of $ {\rm{4.5m}} $.

In Fig. \ref{fig:11}b, the meshing near cavity is shown. Moreover, three coordinate systems are established in the finite element model to extract stress component datum. The first one is the global rectangular coordinate system $ x_{0}O_{0}y_{0} $ established in the \texttt{Assemble} submodule to extract stress datum along ground surface and straight cavity boundary segments $ \overline{S_{2}S_{3}} $ and $ \overline{S_{4}S_{1}} $. The second one is a local polar coordinate system $ r_{1}o_{1}\theta_{1} $ located at the center of cavity (coordinate $ (0,-10) $ in global coordinate system $ x_{0}O_{0}y_{0} $) to extract stress datum along arc cavity boundary $ \overset{\LARGE{\frown}}{S_{1}S_{2}} $. The last one is another local poloar coordinate system $ r_{2}o_{2}\theta_{2} $ located at coordinate $ (4.5,-14.5) $ in global coordinate system $ x_{0}O_{0}y_{0} $ to extract stress datum along arc cavity boundary $ \overset{\LARGE{\frown}}{S_{3}S_{4}} $.

\begin{figure}[htb]
  \centering
  \includegraphics[width = 0.8\textwidth]{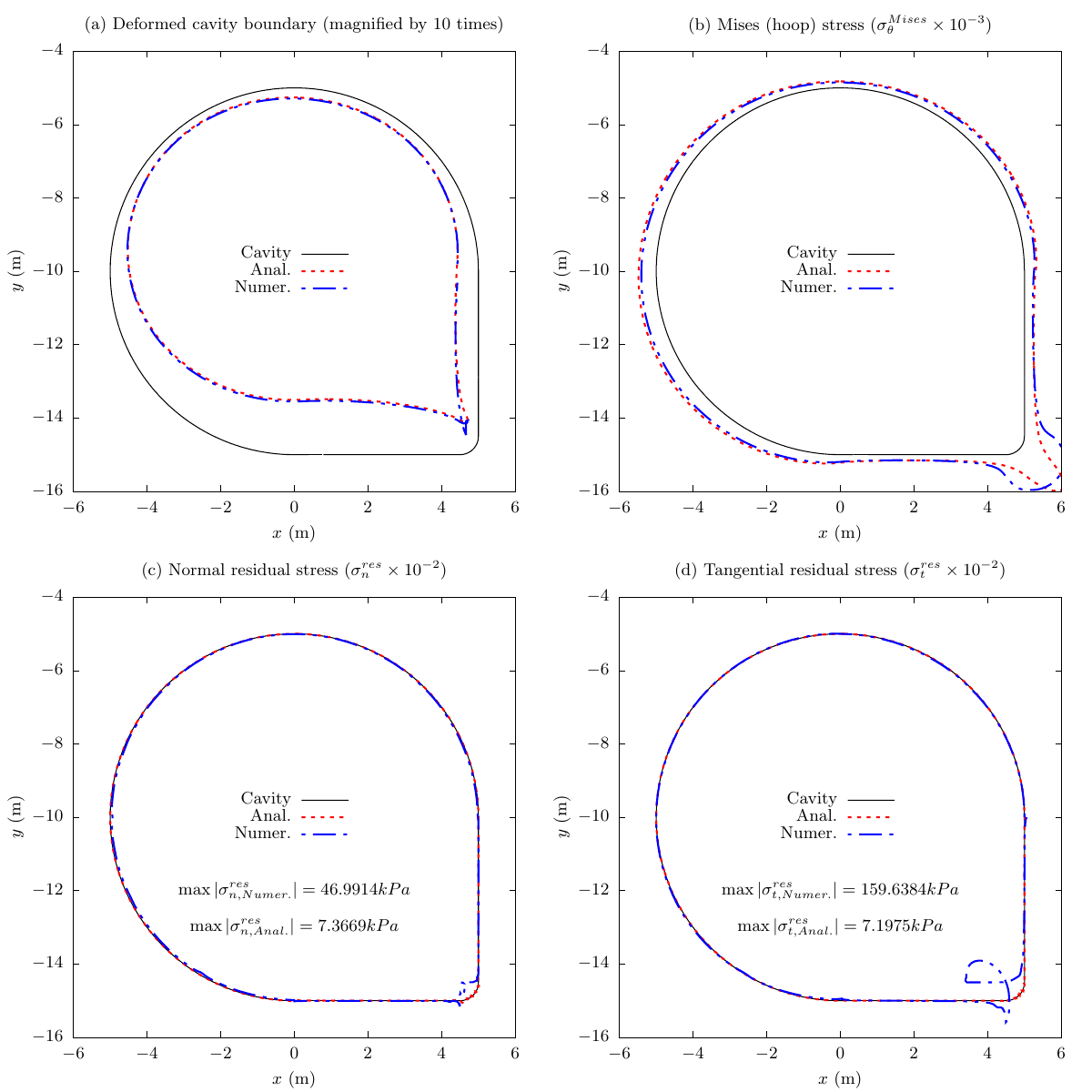}
  \caption{Deformation and stress components along cavity boundary between analytical and numerical solutions}
  \label{fig:12}
\end{figure}

\begin{figure}[htb]
  \centering
  \includegraphics[width = 0.8\textwidth]{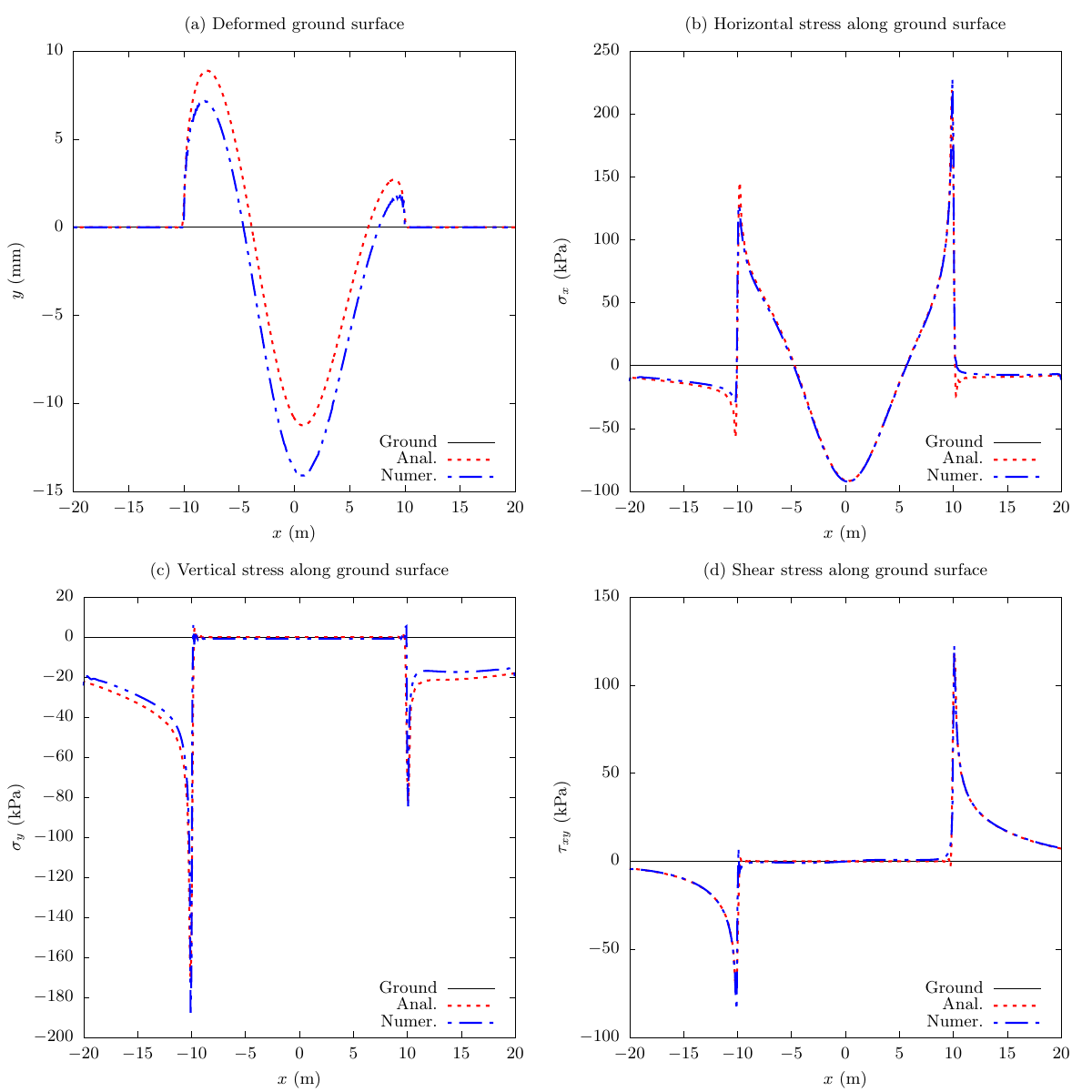}
  \caption{Deformation and stress components along ground surface between analytical and numerical solutions}
  \label{fig:13}
\end{figure}

Substituting all necessary parameters into the analytical solution in Section \ref{sub:Mechanical model and brief solution procedure} and the numerical solution of finite element method in Fig. \ref{fig:11} gives the stress and deformation result comparisons along cavity boundary in Fig. \ref{fig:12} and along ground surface in Fig. \ref{fig:13}.

Figs. \ref{fig:12}a and b show the deformed cavity boundary and the Mises stress $ \sigma_{\theta}^{\rm{Mises}} $ comparisons, which are magnified by $ 10 $ times and reduced by $ 10^{-3} $, respectively, for good visual illustration. We can see that the deformed cavity boundaries and the Mises stress distributions of the analytical and numerical solutions are in good agreements, except for the Mises stress near the arc $ {S_{3}S_{4}} $. Since the seed distributions in Fig. \ref{fig:11}a is dense, the numerical results of finite element solution is expected to be accurate or at least relatively accurate. Thus, we should suspect the accuracy of the analytical solution in Section \ref{sub:Mechanical model and brief solution procedure}. To validate the suspicion, we extract normal and tangential components of residual stress along cavity boundary in both analytical and numerical solutions, which should be zero in theory. The normal residual stress $ \sigma_{n}^{{\rm{res}}} $ and tangential residual stress $ \sigma_{t}^{{\rm{res}}} $ are shown in Figs. \ref{fig:12}c and d with a reduction of $ 10^{-2} $ for good visual illustration. We can see that the residual stress components of both solutions are approximate to zero along cavity boundary segments $ \overset{\LARGE{\frown}}{S_{1}S_{2}} $, $ \overline{S_{2}S_{3}} $, and $ \overline{S_{4}S_{1}} $ as expected. However, both solutions are apparently nonzero along segment $ \overset{\LARGE{\frown}}{S_{3}S_{4}} $, which is the rounded corner. The numerical results are much larger than the analytical one proposed in Section \ref{sub:Mechanical model and brief solution procedure}, indicating that the analytical solution of complex variable solution incorporating the proposed bidirectional conformal mapping has relatively higher accuracy than the numerical one of finite element method.

Fig. \ref{fig:13} shows the stress and displacement component comparisons between the analytical and numerical solutions. Figs. \ref{fig:13}a, c, and d show that the zero-displacement and zero-traction boundary conditions along ground surface in Eqs. (\ref{eq:7.5b}) and (\ref{eq:7.5c}) are well matched, respectively, while all subfigures show that the deformation and stress components of these two solutions are in the same trends with certain discrepancies. The discrepancies are caused by the nonzero residual stress components of finite element solution, as pointed out in Figs. \ref{fig:12}c and d. Therefore, the comparisons in Fig. \ref{fig:13} can only be treated as rough references. The similar trends of stress and displacement components between these two solutions suggest relatively good accuracy of the analytical solution.

\subsection{Second numerical case and comparisons with existing analytical solution}
\label{sec:Second numerical case and comparisons with existing analytical solution}

Certain discrepancies exist in the verification in Section \ref{sub:Results and comparisons with finite element method}, thus, fortified verification of the second numerical case should be conducted to further validate the present analytical solution. As mentioned earlier i Sections \ref{sec:introduction} and \ref{sub:Shallow tunnel}, the bidirectional conformal mapping in this paper is the follow-up research of Ref \cite{lin2024over-under-excavation} (called existing solution or existing one in the text below) with improved bidirectional conformal mapping scheme, and we should compare the complex variable solution incorporating the proposed bidirectional conformal mapping in Section \ref{sub:Mechanical model and brief solution procedure} with the existing one. To be identical, the cavity geometry takes exact the same one in the Section 7.4 of the exisiting solution, which can be analytically expressed as
\begin{equation}
  \label{eq:7.19}
  \left\{
    \begin{aligned}
      & \left( \frac{x}{4} \right)^{2} + \left( \frac{y+10}{5} \right)^{2} = 1, \quad x \leq 0 \\
      & \left( \frac{x}{6} \right)^{2} + \left( \frac{y+10}{5} \right)^{2} = 1, \quad x > 0 \\
    \end{aligned}
  \right.
\end{equation}
The {\emph{collocation points}} along cavity boundary $ \partial{\bm{D}}_{2} $ takes the same distribution in the existing solution. The other mechanical parameters take the same values as those in Table \ref{tab:1}. It should be emphsized that the brief solution procedure in Section \ref{sub:Mechanical model and brief solution procedure} is the same to that in the existing solution, and the total difference between the solution in this paper from the existing one is the reorganized bidirectional conformal mapping proposed in Section \ref{sec:Numerical conformal mapping}.

\begin{figure}[htb]
  \centering
  \includegraphics[width = 0.8\textwidth]{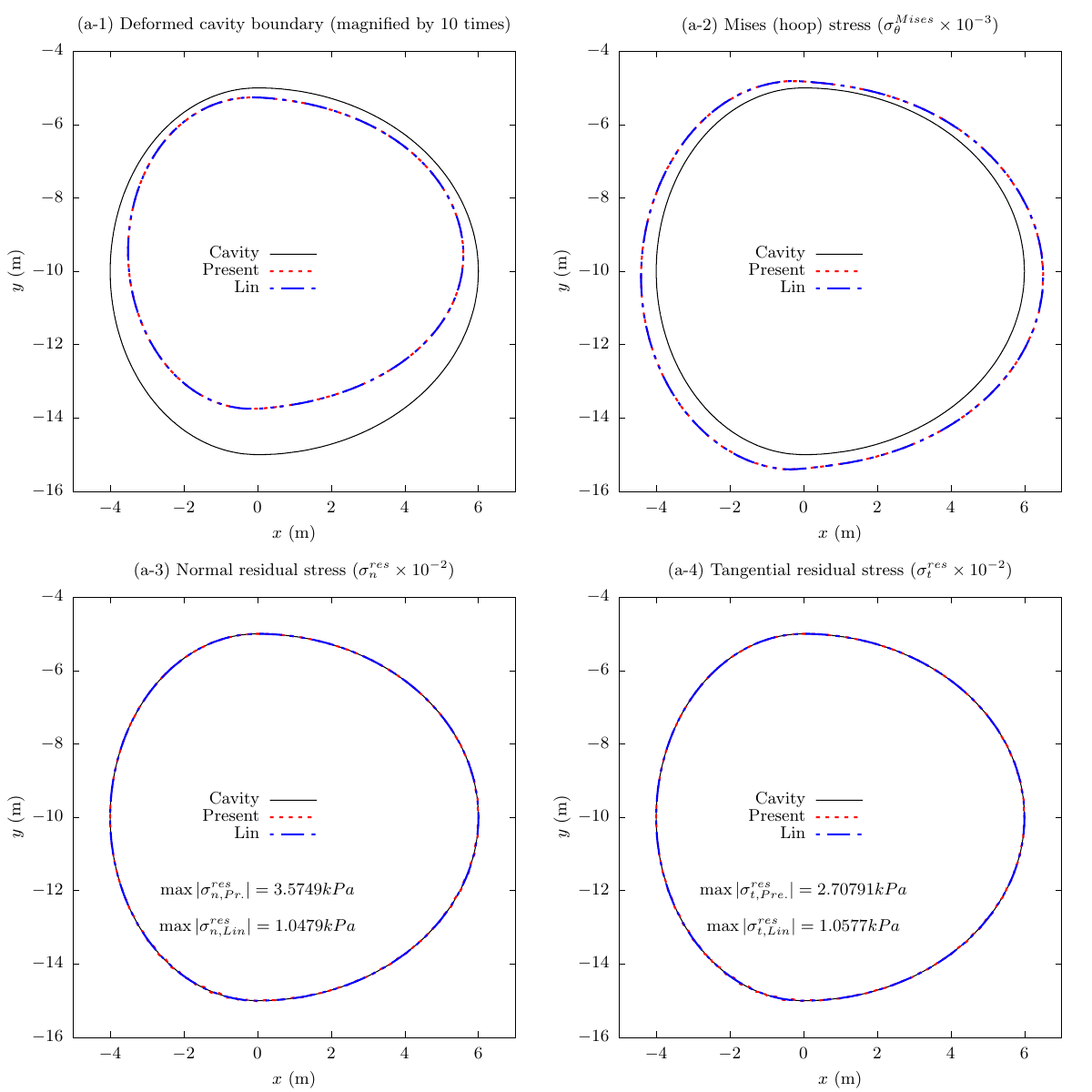}
  \caption{Deformation and stress components along cavity boundary between present solution and Lin's solution \cite{lin2024over-under-excavation}}
  \label{fig:14}
\end{figure}

\begin{figure}[htb]
  \centering
  \includegraphics[width = 0.8\textwidth]{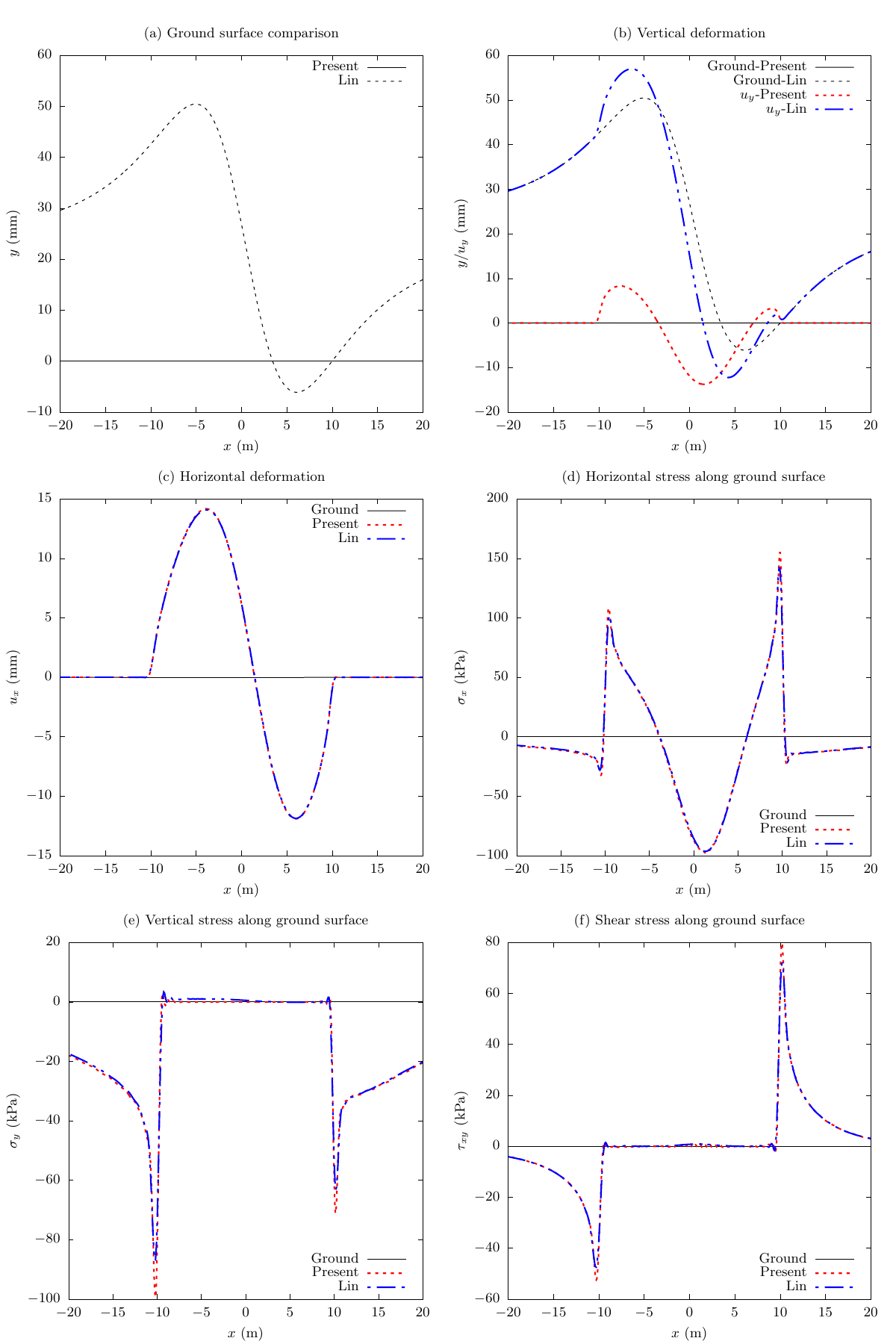}
  \caption{Deformation and stress components along ground surface between present solution and Lin's solution \cite{lin2024over-under-excavation}}
  \label{fig:15}
\end{figure}

Substituting all necessary parameters into present analytical solution and the existing one by Lin et al. gives comparisons of deformation and stress components along cavity boundary in Fig. \ref{fig:14} and along ground surface in Fig. \ref{fig:15}, respectively. Figs. \ref{fig:14}a and b show that the deformation and Mises stress along cavity boundary are almost the same for both solutions, and the residual stress components along cavity boundary in Figs. \ref{fig:14}c and d are very small with approximate maximum values for both solutions. Fig. \ref{fig:14} indicates that the present analytical solution and the existing one can be verified by each other along cavity boundary.

Fig. \ref{fig:15}a shows the comparisons of ground surfacs computed by the two-step backward conformal mappings in this paper and the one in Ref \cite{lin2024over-under-excavation}. Apparently, the backward conformal mapping in this paper restores the ground surface as expected, while certain errors exist in the existing one, which would further cause certain errors in stress and deformation along ground surface. Figs. \ref{fig:15}b-f show the comparisons of stress and deformation components between the solution incorporating the newly proposed bidirectional conformal mapping in Section \ref{sec:Numerical conformal mapping} and the existing one. Fig. \ref{fig:15}b shows the superposition of vertical deformations of present analytical solution and the exsiting one to the ground surface curves in Fig. \ref{fig:15}a, and the vertical deformation solved by the existing solution of Lin et al. is apparently large due to the error of computed ground surface of the backward mapping conformal mapping. Figs. \ref{fig:15}c-f show similar trends of horizontal deformation and stress components, respectively, and certain discrepancies exist due to the error of computed ground surface of the backward conformal mapping in the existing solution. Fig. \ref{fig:15} indicates that the present analytical solution incorporating the newly proposed bidirectional conformal mapping scheme in Section \ref{sec:Numerical conformal mapping} is more reasonable than the existing one with certain improvements, especially in vertical deformation.

\section{Conclusions}
\label{sec:Conclusions}

This paper presents a new bidirectional conformal mapping \added{scheme both geometrically and mechanically applicable} for \replaced{over-break and under-break}{noncircular and asymmetrical} cavitiy in both deep and shallow \replaced{tunnelling}{engineering}\replaced{, which is logically straightforward and computationally simple to conduct by solving only a pair of forward and backward linear systems, comparing to existing mapping schemes.}{ Comparing to existing solution schemes for conformal mappings of deep and shallow tunnels, the new bidirectional conformal mapping is simple to solve and easy to code by only involving a pair of linear systems.} \replaced{The new mapping scheme is validated by several numerical cases of over-break and under-break cavities in both deep and shallow tunnelling, indicating the great potential of the new mapping scheme in geometrical application.}{The new bidirectional conformal mapping is also suitable for cavity with sharp corners with newly developed numerical strategies of corresponding small arc simulation and densified collocation points. Several numerical cases of deep and shallow tunnels are conducted to show the geometrical adaptivity of the new bidirectional conformal mapping.} \replaced{The new mapping scheme is subsequently embeded into the complex variable solutions of under-break shallow tunnelling in gravitational geomaterial with fixed far-field ground surface. The good agreements of result comparisons with corrsponding finite elment solution and existing analytical solution both suggest the great potential of the new mapping scheme in further mechanical application in tunnel engineering.}{Subsequent complex variable solutions embedded by the new bidirectional conformal mapping respectively show good agreements with finite element solution and existing analytical solution, indicating the strong potential of the new mapping in mechanical analyses of tunnel engineering.} \added{Therefore, the new bidirectional conformal mapping scheme would significantly extend the application range of the complex variable method in tunnel engineering by switching the cavity contour from a symmetrical one in theory to an asymmetrical one in practice with consideration of over-break and under-break excavation.}






\bmsection*{Acknowledgments}

This study is financially supported by the Fujian Provincial Natural Science Foundation of China (Grant No. 2022J05190 and 2023J01938), the Scientific Research Foundation of Fujian University of Technology (Grant No. GY-Z20094, GY-Z21026, and GY-H-22050), and the National Natural Science Foundation of China (Grant No. 52178318). The authors would like to thank Ph.D. Yiqun Huang for his suggestion on this study.

\bmsection*{Financial disclosure}

None reported.

\bmsection*{Conflict of interest}

The authors declare no potential conflict of interests.

\clearpage
\bibliography{biblio}

\clearpage
\begin{table}[htb]
  \centering
  \begin{tabular}{cccccc}
    \toprule
    $ E {\rm{(MPa)}} $ & $ \nu $ & $ k_{x} $ & $ \gamma {\rm{(kN/m^{3})}} $ & $ T_{1} $ & $ T_{2} $ \\
    \midrule
    $ 20 $ & $ 0.3 $ & $ 0.8 $ & $ 20 $ & $ (-10,0) $ & $ (10,0) $ \\
    \bottomrule
  \end{tabular}
  \caption{Mechanical parameters in case verification}
  \label{tab:1}
\end{table}

\clearpage

\clearpage

\clearpage

\clearpage

\clearpage

\clearpage

\clearpage

\clearpage

\clearpage

\clearpage

\clearpage

\clearpage

\clearpage

\clearpage

\clearpage

\end{document}